\documentclass{siamart220329}
\usepackage{mathrsfs}
\usepackage{amssymb}
\usepackage{bm}
\usepackage[algo2e]{algorithm2e}
\usepackage{tikz,pgfplots}
\pgfplotsset{compat=1.15}
\usepackage{amsfonts}   
\usepackage{subfig} 
\usepackage{graphicx}       
\usepackage{braket}
\usepackage{booktabs}
\usepackage{multirow}
\usepackage{stmaryrd}
\usepackage{grffile}

\newlength{\tempdima}
\newcommand{\rowname}[1]
{\rotatebox{90}{\makebox[\tempdima][c]{\textbf{#1}}}}

\newtheorem{remark}{Remark}
\numberwithin{equation}{section}

\newcommand\bx{\boldsymbol{x}}

\newcommand\bu{\boldsymbol{u}}
\newcommand\bv{{\boldsymbol{v}}}

\newcommand{\bbf}{\bm{f}}
\newcommand{\br}{\bm{r}}
\newcommand{\bV}{\bm{V}}
\newcommand{\bW}{\bm{W}}

\newcommand\bg{\boldsymbol{g}}

\newcommand\bh{\boldsymbol{h}}
\newcommand\bbR{\mathbb{R}}
\newcommand\bbN{\mathbb{N}}
\newcommand\bbS{\mathbb{S}}

\newcommand{\bbm}{\bm{m}}
\newcommand{\bU}{{\bf U}}

\newcommand\Kn{{\rm Kn}}

\newcommand\mQ{\mathcal{Q}}
\newcommand\mM{\mathcal{M}}

\newcommand{\msF}{\mathscr{F}}
\newcommand{\msW}{\mathscr{W}}
\newcommand{\msM}{\mathscr{M}}
\newcommand{\msK}{\mathscr{K}}
\newcommand{\msT}{\mathscr{T}}

\newcommand\mL{\mathcal{L}}
\newcommand{\mO}{\mathcal{O}}
\newcommand{\bfA}{{\bf A}}
\newcommand{\bfQ}{{\bf Q}}

\newcommand{\bP}{\bm{P}}
\newcommand{\bQ}{\bm{Q}}
\newcommand{\bR}{\bm{R}}
\newcommand{\llb}{\llbracket}
\newcommand{\rrb}{\rrbracket}
\newcommand{\dd}{\; \mathrm{d}}
\newcommand\nnBGK{{\rm NR^{\rm BGK}}}
\newcommand\nnLR{{{\rm NSR^{\rm BGK}_{\rm LR}}}}
\newcommand\nnFSM{{\rm NR^{\rm Quad}}}
\newcommand\nnLA{{{\rm NSR^{\rm Quad}_{\rm LA}}}}

\newcommand{\fneq}{ $f^{\rm neq}$ }
\graphicspath{{media/}}     

\DeclareMathOperator{\Span}{span}

\title{Solving Boltzmann equation with neural sparse representation
}

\author{
  Zhengyi Li\footnote{ School of Mathematical Science,
  Peking University,
  Beijing, China, email: \texttt{lizhengyi@pku.edu.cn}.},
  ~~Yanli Wang\footnote{Beijing Computational Science Research Center,
  Beijing, China, email:
  \texttt{ylwang@csrc.ac.cn}.}
  ~~Hongsheng Liu\footnote{Huawei Technologies Co. Ltd, email:
  \texttt{liuhongsheng4@huawei.com}}
  ~~Zidong Wang\footnote{Huawei Technologies Co. Ltd, email:
  \texttt{wang1@huawei.com}}
  ~~
  Bin Dong\footnote{ Beijing International Center for Mathematical Research \&
  Center for Machine Learning Research,
  Peking University,
  Beijing, China, email:
  \texttt{dongbin@math.pku.edu.cn}.}}

\begin{document}
\maketitle

\begin{abstract}
We consider the neural sparse representation to solve Boltzmann equation with BGK and quadratic collision model, where a network-based ansatz that can approximate the distribution function with extremely high efficiency is proposed. Precisely, fully connected neural networks are employed in the time and spatial space so as to avoid the discretization in space and time. The different low-rank representations are utilized in the microscopic velocity for the BGK and quadratic collision model, resulting in a significant reduction in the degree of freedom. We approximate the discrete velocity distribution in the BGK model using the canonical polyadic decomposition. For the quadratic collision model, a data-driven, SVD-based linear basis is built based on the BGK solution. All these will significantly improve the efficiency of the network when solving Boltzmann equation. Moreover, the specially designed adaptive-weight loss function is proposed with the strategies as multi-scale input and Maxwellian splitting applied to further enhance the approximation efficiency and speed up the learning process. Several numerical experiments, including 1D wave and Sod problems and 2D wave problem, demonstrate the effectiveness of these neural sparse representation methods.

\end{abstract}

\noindent {\bf keyword}: Boltzmann equation, BGK model, quadratic collision, canonical polyadic decomposition, singular value decomposition

\section{Introduction}
People are interested in the simulation of the kinetic theory, due to its extensive applications in the engineering fields, such as aerospace, plasma, and micro-electro-mechanical systems. However, Boltzamnn equation as one of the most important governing equations in kinetic theory, it is quite difficult to solve efficiently and accurately. The main difficulty lies in the high dimensionality, including time, spatial space and microscopic velocity space, and the complex quadratic collision with high dimensional integral and singular collision kernel.

Nowadays, there are several kinds of methods to solve Boltzmann equation. For example, the statistical method as the direct simulation Monte Carlo (DSMC) is brought up in \cite{bird1994molecular}, which directly solves Boltzmann equation with randomness. But it is limited by its low efficiency and the statistical noise. Another kind of method is the deterministic method, such as the discrete velocity methods \cite{liu2020unified}, which solves Boltzmann equation by discreting the distribution function at several discrete velocity points. Fourier spectral method  \cite{mouhot2006fast, wu2013deterministic, gamba2017fast} has also made great progress by approximating the distribution function with trigonometric functions. Recently, Hermite spectral methods are successfully adopted to solve the quadratic collision model \cite{wang2019approximation}. Another important method is the moment method, which is proposed by Grad \cite{grad1949kinetic}, but it is limited by the non-hyperbolicity of the Grad moment equations, even near Maxwellian. The asymptotic-preserving scheme \cite{jin2010micromacro} is also proposed for the Boltzmann equation, and we refer \cite{Dimarco2014} for a comprehensive review of these methods. The low-rank decomposition is applied to numerically solving kinetic equations recently. The adaptive dynamic low-rank method is proposed in  \cite{hu2021adaptive, koch2007dynamicala, lubich2014projectorsplitting} for Boltzmann equation, and a local macroscopic conservation low-rank method is brought up for the Vlasov equation in \cite{guo2022local}. Methods based on higher-order tensor decomposition \cite{ibragimov2009three, boelens2020tensor, chikitkin2021numerical} are also applied for Boltzmann equation. 


Recently, with the development of computers, more and more network-based methods are proposed to solve Boltzmann and other kinetic equations. There are generally two kinds of methods. The first one is combining neural networks and the reduced model of Boltzmann such as the moment models \cite{grad1949kinetic} to learn a closed reduced model \cite{han2019uniformly, huang2021machinea, schotthofer2022neural, li2022learning}. The network-based method is first utilized to learn the moment closure relation for Boltzmann equation in \cite{han2019uniformly}, and the moment closure models which preserve several physical invariances are learned in \cite{li2022learning}. A fast moment closure approximation based on the max entropy method and neural networks is proposed in \cite{schotthofer2022neural}. Other than Boltzmann equation, a neural network-based moment closure model which preserves the hyperbolicity of radiative transfer equation is brought up in \cite{huang2021machinea}. The other kind of method is to solve Boltzmann equation directly in the framework of PINN \cite{raissi2019physics}. PINN was first utilized to solve Boltzmann equation in \cite{lou2021physicsinformed}, but it is only for the BGK model. An asymptotic-preserving network-based method is proposed for linear transport equations in \cite{jin2021asymptoticpreserving}. There is also  some work on the quadratic collision model. In \cite{holloway2021acceleration}, the quadratic collision terms are first approximated by an auto-encoder. Then, with the idea of reduced order models, a new set of basis is learned in \cite{alekseenko2022fast}, where the quadratic collision term is computed with computational effort reduced significantly. However, only the spatially homogeneous problems are tested in both works. In \cite{xiao2021using}, on the other hand, the neural network is adopted to approximate the quadratic collision terms directly, where spatially 1D and 2D examples are studied. DSMC method is utilized in \cite{miller2022neuralnetwork} to generate the training data, based on which the quadratic collision term used at the moment model is obtained.

Generally speaking, to solve Boltzmann equation efficiently is to find an ansatz that can approximate the distribution function well. In this work, we propose a new ansatz, named neural sparse representation (NSR), for the Boltzmann equation based on neural networks. NSR is a promising ansatz for the BGK and quadratic collision model and can express the distribution function with significantly fewer degrees of freedom compared to most traditional methods. In the framework of NSR, the network is built within the framework of the discrete velocity method (DVM), where the input is spatial position and time, and the output is the distribution function at the discrete velocity points. This eliminates the need for discretization in time and space, significantly reducing the degree of freedom and overcoming the curse of dimensionality. This makes the extension of NSR method to high dimensionality without difficulty. To reduce the degree of freedom in the microscopic velocity space, the low-rank property of the distribution function is explored, and different approximation methods are proposed for different collision models. In the BGK model, it is approximated by the canonical polyadic decomposition (CPD), which is widely used in PDE solving \cite{boelens2018parallel, reynolds2016randomized} and accelerated deep learning \cite{lebedev2015speedingup, ji2019survey}. For the quadratic collision model, due to its complexity, the BGK solution is utilized to construct a series of data-driven basis with SVD decomposition to approximate the quadratic collision terms. Moreover, based on prior knowledge of the Boltzmann equation, the multi-scale input is realized by scaling the time and spatial space with parameters at different magnitudes, which can match the multi-scale property of Boltzmann equation. The distribution function is processed using the Maxwellian splitting strategy to capture the behavior of the macroscopic variables much easier. These two structures further improve the approximation efficiency of NSR. An adaptive-weighted loss function is specially designed for the network. Except for the PDE residual loss, which is usually contained in the loss function, the loss from macroscopic variables is added into the loss function to match the property of Boltzmann equation, such as the fact that the density of the PDE residual should be zero. In addition, since the contribution of different microscopic velocity points is different, adaptive weights are added to the error of each microscopic velocity point. All these techniques are employed to accelerate the process and enhance the approximation efficiency of NSR.

The effectiveness of the proposed neural network-based methods is validated through several numerical experiments. These experiments comprise one-dimensional problems with both continuous and discontinuous initial conditions, in addition to a two-dimensional wave problem. To further verify the accuracy and efficiency of the methods, transfer learning is employed to investigate the computational time for the two-dimensional wave problem.

The rest of this paper is organized as follows. In Sec. \ref{sec:Boltzmann}, the Boltzmann and its related properties are introduced. In sec. \ref{sec:net_framework}, the general structure of the network is proposed. The neural sparse representation of the BGK model and the quadratic collision model is discussed in Sec. \ref{sec:lr} and \ref{sec:svd}, respectively. The numerical experiments are presented in Sec. \ref{sec:num}, with some concluding remarks in Sec. \ref{sec:conclusion}. 
\section{Boltzmann equation}
\label{sec:Boltzmann}
In this section, we will introduce the Boltzmann equation, which describes a particle system from a statistical point of view. It has the form below
\begin{equation}
\label{eq:Boltzmann}
\frac{\partial f(\bm x,\bm v,t)}{\partial t}+\bm{v} \cdot\nabla_{\bm x} f(\bm x,\bm v,t)=\mathcal{Q}[f](\bx, \bv,t), \qquad  t \in \bbR^+, \quad \bx \in \bbR^3, \quad \bv \in \bbR^3,
\end{equation}
where $f( \bx, \bv, t)$ is the distribution function. Here, $t$ is the time, $\bx$ is the spatial coordinates, and $\bv$ stands for the microscopic velocity of the particles. $\mQ[f]$ is the collision operator which has a quadratic form 
\begin{equation}
\label{eq:binary-collsion}
\mathcal{Q}[f](\bx, \bv, t) = \mQ(f, f) = \int_{\bbR^3}\int_{\bbS^2}B(\bv - \bv_{\ast}, \sigma)[f(\bm v'_*)f(\bm v')-f(\bm v_*)f(\bm v)] \dd \sigma \dd \bm v_*,
\end{equation}
where $\bm v$ and $\bm v_*$ are the velocities of two particles before the collision, and the velocities of $\bm v'$ and $\bm v'_*$ are the velocities of the particles after the collision. From the conservation of momentum and energy during such a collision, we obtain the following relationship
\begin{equation}
\label{eq:pre_after_col}
\bm v'=\frac{\bm v+\bm v_*}{2}+\frac{\left|\bm v-\bm v_*\right|}{2} \sigma, \qquad\bm v'_*=\frac{\bm v+\bm v_*}{2}-\frac{\left|\bm v-\bm v_*\right|}{2}  \sigma,
\end{equation}
where $\sigma\in \mathbb{S}^2$ is a unit vector. The collision kernel $B(\bv - \bv_{\ast}, \sigma) \geqslant 0$ is a non-negative function that depends on $\bv - \bv_{\ast}$, and  cosine of the derivation angle $\theta$. It is often written as 
\begin{equation}
\label{eq:B}
B\left(\bm v-\bm v_*, \sigma\right)=\Phi\left(\left|\bm v-\bm v_*\right|\right) b(\cos \theta), \quad \cos \theta=\frac{\sigma \cdot\left(\bm v-\bm v_*\right)}{\left|\bm v-\bm v_*\right|}.
\end{equation}
The specific form of $B$, which characterizes the details of the binary interactions, is determined by the mutual potential between the particles. 
One of the commonly used collision kernels is the variable hard sphere (VHS) model proposed by Bird \cite{bird1994molecular} as 
\begin{equation}
    \label{eq:vhs}
    B = C_{\alpha} |\bv - \bv_{\ast}|^{\alpha}, \qquad \alpha = \frac{\eta - 5}{\eta - 1}.
\end{equation}
Here, constant $C_{\alpha}$ is empirically determined \cite{wu2013deterministic}. $\eta$ is the index in the power of distance. The case $\eta > 5$ corresponds to the ``hard potential'' and the case $\eta < 5$ corresponds to the ``soft potential''. The collision kernel $B(\bv - \bv_{\ast}, \sigma)$ is independent of $|\bv - \bv_{\ast}|$ when $\eta = 5$, and is called ``Maxwell molecules'' in this case. 
However, in the VHS model, the differential cross-section $C_{\alpha}|\bv-\bv_{\ast}|^{\alpha-1}$ is independent of the deflection angle. Then, the following anisotropic collision kernel with $\theta$ included in the cross-section is suggested in \cite{mouhot2006fast,wu2013deterministic} as 
\begin{equation}
    \label{eq:vss}
    B = C_{\alpha}' \sin^{\alpha-1}\left(\frac\theta 2 \right) |\bv - \bv_{\ast}|^{\alpha},  
\end{equation}
where $C_{\alpha}'$ is a constant as 
\begin{equation}
    \label{eq:coe_C}
    C_{\alpha}' =\frac{ (\alpha+3)(\alpha+5)}{24}  C_{\alpha}. 
\end{equation}
Model \eqref{eq:vss} is utilized in the numerical simulations in this work.

Due to the complicity of the quadratic collision model \eqref{eq:binary-collsion}, several simplified collision models are proposed, such as the BGK collision model \cite{bhatnagar1954model}
\begin{equation}
    \label{eq:BGK}
\mathcal{Q}^{\rm BGK}[f]=\frac{1}{\tau}(\mathcal{M}-f). 
\end{equation}
Here, $\tau$ is the mean relaxation time, and $\mM$ is the local equilibrium which is called Maxwellian \cite{maxwell1860illustrationsa}, and has  the form below 
\begin{equation}
\label{eq:Maxwellian}
\mathcal{M}=\frac{\rho}{\sqrt{2\pi T}^3}\exp \left( -\frac{|\bm v-\bm u|^2}{2 T} \right).
\end{equation}
Here, $\rho, \bu, T$ is the density, macroscopic velocity, and temperature, respectively.  Their relationship with the distribution function is 
\begin{equation}
\label{eq:macro}
    \begin{gathered}
   \rho(\bm x,t)=\Braket{f}, \qquad 
   \bm{m}(\bx, t) \triangleq  \rho \bm u(\bm x,t)=\Braket{\bm v f},\\
E(\bx, t) \triangleq \frac{3}{2}\rho T(\bx, t) + \frac{1}{2}\rho |\bu|^2 = \frac{1}{2}\Braket{\vert\bm v\vert^2 f},
\end{gathered}
\end{equation}
where $\Braket{\cdot}$ is defined as 
\begin{equation}
    \label{eq:inner_fun}
\Braket{\cdot} =\int_{\mathbb{R}^3} \cdot \dd \bm v.
\end{equation}
The norm for a vector $|\cdot|$ is defined as 
\begin{equation}
    \label{eq:norm_vec}
    |\bm{g}|^2= \sum_{i = 1}^N g_i^2, \qquad \forall \bm{g} \in \bbR^{1\times M} \text{~or~} 
   \bm{g} \in \bbR^{M} , \qquad M \in \bbN^{+}. 
\end{equation}
$\bm{m}(\bm x, t)$ and $E(\bx, t)$ is the momentum and total energy respectively. Moreover, the total energy $E(\bx, t)$ is separated into energy in three directions as 
\begin{equation}
    \label{eq:sep_E}
    E_i(\bx, t) = \frac{1}{2}\Braket{(v^{(i)})^2 f}, \qquad i = 1, 2, 3. 
\end{equation}
The collision terms have some mathematical and physical properties that apply to arbitrary forms of collision terms, such as the conservation of mass, momentum, and energy as 
\begin{equation}\label{eq:Q-cons}
\left\langle
     \left(
    \begin{array}{c}
         1  \\
         \bv \\
         |\bv|^2
    \end{array} 
    \right)
    \mQ[f]
\right\rangle = \bm{0}. 
\end{equation}

For now, we have introduced the Boltzmann equation. There are several classical numerical methods to solve it, such as the fast Fourier spectral method \cite{mouhot2006fast, gamba2017fast}, the moment method \cite{struchtrup2005macroscopic}, the discrete velocity method \cite{liu2020unified}, and DSMC method \cite{bird1994molecular}. But, it is still a challenge to solve Boltzmann equation numerically due to its high dimensionality, quadratic collision operator, etc. In the following sections, a neural sparse representation which is a high-quality ansatz of the distribution function is proposed, and several strategies in the network-based method are brought up to solve Boltzmann equation efficiently.




\section{Neural representation for Boltzmann equation}
\label{sec:net_framework}
In this section, the general framework of the neural representation for the Boltzmann equation will be introduced. When solving Boltzmann equation using neural network, we first discretize Boltzmann equation in the microscopic velocity space to obtain a semi-discrete system which is discussed in Sec. \ref{sec:dis_Bol}. Then, a fully connected neural network is utilized to serve as a parameterized ansatz for the semi-discrete system, and the general network architecture is introduced in Sec. \ref{sec:network}. To optimize the network parameters, a specially designed loss function is proposed, which is discussed in Sec. \ref{sec:loss}.


\subsection{Discretization in the microscopic velocity space}
\label{sec:dis_Bol}
In this section, we will introduce the discretization in the microscopic velocity space, and the semi-discrete system is proposed. The discrete velocity method \cite{liu2020unified} is utilized here to discrete the Boltzmann equation. Assuming the series of points in each direction of the microscopic velocity space are 
\begin{equation}
    \label{eq:dvm_point_dir}
   \bV^{(i)}  = \left[v_1^{(i)}, v_2^{(i)}, \cdots, v_{N_i}^{(i)}\right]\in \bbR^{1 \times N_i}, \qquad i = 1, 2, 3,
\end{equation}
with weights 
\begin{equation}
    \label{eq:point_weight}  
    \bW^{(i)} = [\omega^{(i)}_1, \omega^{(i)}_2, \cdots, \omega^{(i)}_{N_i}] \in \bbR^{1\times N_i}, \qquad i = 1, 2, 3.
\end{equation}
Then, all the points 
\begin{equation}
    \label{eq:single_point}
    \left(v_{l_1}^{(1)}, v_{l_2}^{(2)}, v_{l_3}^{(3)}\right)^T, \qquad 1\leqslant l_i \leqslant N_i, \quad i = 1, 2, 3, 
\end{equation}
make up the full discrete points in the microscopic velocity 
\begin{equation}
\label{eq:dvm_points}
    \bV \triangleq [\bv_1, \bv_2, \cdots, \bv_N] \in \bbR^{3\times N_v}, \qquad \bv_l =  \left(v_{l_1}^{(1)}, v_{l_2}^{(2)}, v_{l_3}^{(3)}\right)^T , 
\end{equation}
with $N_v = \prod_{i = 1}^3 N_i$ and the corresponding weight
\begin{equation}
    \label{eq:weight}
    \bW \triangleq [\omega_1, \omega_2, \cdots, \omega_N] \in \bbR^{1\times N_v}, \qquad \omega_l =\prod\limits_{i=1}^3 \omega_{l_i}^{(i)}.
\end{equation}
We want to emphasize that there exists  a one-to-one mapping between $l$ and $(l_1, l_2, l_3)$, and it will be not listed explicitly. Thus, the discrete distribution functions are 
\begin{equation}
    \label{eq:dis_f}
    \bbf(\bx, t) \triangleq [f_1(\bx, t), f_2(\bx, t), \cdots, f_{N_v}(\bx, t)]\in \bbR^{1\times N_v}, \qquad 
    f_i(\bx, t) = f(\bx, \bv_i, t).
\end{equation}
The macroscopic variables \eqref{eq:macro} and \eqref{eq:sep_E} can be expressed as 
\begin{equation}
    \label{eq:dis_macro_var} 
    \begin{gathered}
        \rho[\bbf] = \bW \bbf^T, \qquad \bm{m}[\bbf] = \bV (\bW \times \bbf)^T, \qquad 
        E[\bbf] = |\bV|^2(\bW \times \bbf)^T, \qquad 
        E_i[\bbf] = |\bV_i|^2(\bW \times \bbf)^T, \quad i = 1, 2, 3. 
    \end{gathered}
\end{equation}
where 
\begin{equation}
    \label{eq:times}
    (\bW \times \bbf)_{ij}=\bW_{ij}\bbf_{ij},
\end{equation}
and $|\bV|^2 = (|\bv_1|^2, |\bv_2|^2, \cdots, |\bv_{N_v}|^2) \in \bbR^{1\times N_v}$. $\bV_i \in \bbR^{1\times N_v}$ is the $i$-th row of $\bV$. For the collision term \eqref{eq:binary-collsion}, the fast Fourier spectral method \cite{mouhot2006fast} is utilized here, and the discrete collision term is labeled as 
\begin{equation}
    \label{eq:dis_Q}
\mQ_i = \mQ[f_1, \cdots, f_{N_v}](\bx, \bv_i, t).
\end{equation}
Then, the Boltzmann equation \eqref{eq:Boltzmann} is reduced into 
\begin{equation}
\label{eq:dis_Boltzmann}
\begin{cases}  
\frac{\partial f_1(\bm x,t)}{\partial t}+\bv_1\cdot \nabla_{\bm x} f_1(\bm x,t)=\mathcal{Q}_1,\\  
\qquad \vdots\\  
\frac{\partial f_{N_v}(\bm x,t)}{\partial t}+\bv_n\cdot \nabla_{\bm x} f_{N_v}(\bm x,t)=\mathcal{Q}_{N_v},
\end{cases}
\end{equation}
Let
\begin{equation}
    \label{eq:vector_Q}
    \mQ[\bbf] \triangleq [\mQ_1, \mQ_2, \cdots, \mQ_{N_v}]^T.
\end{equation}
For the BGK collision model, substituting $\bbf$ into \eqref{eq:BGK}, the discrete collision term is derived directly, which is also denoted as \eqref{eq:vector_Q}. Thus, the system \eqref{eq:dis_Boltzmann} is rewritten as 
\begin{equation} 
\label{eq:vec_dis_Bol}
\frac{\partial \bm{f}(\bm x,t)}{\partial t}+\bV \cdot \nabla_{\bm x} \bm{f}(\bm  x,t)=\mathcal{Q}[\bm{f}].
\end{equation}

In the solving process, a neural network is utilized to approximate the discrete distribution function $\bm f$, whose input is the spatial position $\bx$ and time $t$, and the output is the distribution function value at the fixed discrete velocity. The system \eqref{eq:vec_dis_Bol} is adopted as the governing equation in the loss function, the details of which will be introduced in the following sections.

\subsection{Network architecture}
\label{sec:network}
In this section, the architecture of the network is introduced, where a fully connected neural network is utilized. The general structure of the $L$-layer fully connected neural network is composed of $L$ fully connected layers, each of which consists of a linear transformation $F^{(l)}$ and an activation function $\sigma^{(l)}$. The whole form of the network has the form below
\begin{equation}
\begin{aligned}
\label{eq:full_network}
y(\cdot):=y^{(L)}=\sigma^{(L)}\circ F^{(L)} \circ \sigma^{(L-1)} \circ F^{(L-1)} \circ \cdots \circ \sigma^{(1)} \circ F^{(1)}(y^{(0)}),
\end{aligned}
\end{equation}
where the $l$-th layer is a mapping from $\mathbb{R}^{m_{l-1}}$ to $\mathbb{R}^{m_l}$ as 
\begin{equation}
\label{eq:full_connect_network}
y^{(l)}= \sigma^{(l)}\circ F^{(l)}\left(y^{(l-1)}\right), \qquad  y^{(l-1)} \in \bbR^{m_{l-1}}, \quad y^{(l)} \in \bbR^{m_l}.
\end{equation}
Substituting the specific form of $F^{(l)}$ into \eqref{eq:full_connect_network}, it could be rewritten as  
\begin{equation}
    \label{eq:full_connect_network_entry}
    y^{(l)}_j= \sigma^{(l)}\left(\sum_{i=1}^{m_{l-1}} {\bf W}_{ji}^{(l)}y^{(l-1)}_i + b^{(l)}_j\right), \quad j=1,...,m_{l},
\end{equation}
where $y^{(l)}$ is the output of the $l$-th layer, $y^{(0)}$ is the input and $y^{(L)}$ is the output of the network. $m_l$ is the dimension of the $l$-th layer, $m_0$ is the dimension of the input and  $m_L$ is the dimension of the output. $\sigma^{(l)}, l = 1, \cdots L$ here is the activation function in the $l$-th layer, and we will choose $\sigma^{(L)}$ as the identity and the rest as some nonlinear activation function. Specially, the sine activation function 
\begin{equation}
    \label{eq:sin_activation}
    \sigma(x)=\sin(x)
\end{equation}
is utilized here. This is because sin activation functions are thought to exhibit better representation \cite{sitzmann2020implicit} in implicit neural expression than activation functions such as tanh or softplus. ${\bf W}^{(l)}$, which is a $\mathbb{R}^{m_{l} \times m_{l-1}}$ matrix, is the weight of the $l$-th layer and  $ b^{(l)}\in \mathbb{R}^{m_{l}}$ is the bias of the $l$-th layer. $\{{\bf W}^{(l)},b^{(l)}\}_{l=1}^{L}$ constitutes the parameter sets of the neural network.

Solving PDE using a neural network is the process to determine parameters $\{W^{(l)},b^{(l)}\}_{l=1}^{L}$ using the optimization algorithms with proper loss functions. Here, the input is the spatial space and time as $y^{(0)} = (\bx, t) \in \bbR^4$ and the output is the discrete distribution function $\bbf$ in \eqref{eq:dis_f}. For now, the general form of the neural network to approximate Boltzmann equation is proposed. To improve the approximation efficiency, two strategies as multi-scale input and Maxwellian splitting are utilized.

\paragraph{\bf Multi-scale input} 
Multi-scale neural networks are always utilized when approximating problems with high frequency \cite{liu2020multiscale, huang2021solving}, whose main idea is to convert the learning of data with high frequency to that with low frequency. Related studies show that using multi-scale networks or multi-scale inputs \cite{tancik2020fourier} to change the function frequency can improve the convergence speed of networks. Due to the multi-scale property of Boltzmann equation, the strategy of multi-scale input is adopted here to improve the approximation efficiency of the network. Precisely, a relatively straightforward process is utilized by multiplying the inputs of the network by a sequence of constants as 
\begin{equation}
    y^{(0)}_{\rm multi}=\left(\begin{array}{c}
         c_1  \\
         c_2 \\
         \vdots \\
         c_{n_c}
    \end{array}\right)
     y^{(0)} =\left(\begin{array}{c}
         c_1  \\
         c_2 \\
         \vdots \\
         c_{n_c}
    \end{array}\right)(\bx, t) \in \mathbb{R}^{n_c} \times \bbR^4.
\end{equation}
Here, we want to emphasize that the constants $c_i, i = 1, \cdots n_c$ are problem dependent, and chosen empirically  in the numerical tests. Besides, there is no theoretical evidence right now to show how this strategy will affect the convergence speed of the neural network when approximating the PDE model, it is adopted here following the results in  \cite{liu2020multiscale}. 
 
\paragraph{\bf Maxwellian splitting}
Another strategy we adopt here is Maxwellian splitting, which is also called Micro-Macro decomposition \cite{jin2010micromacro,jang2015high}, where the distribution function is split into the equilibrium whose form is the Maxwellian \eqref{eq:Maxwellian} and the non-equilibrium residual. This splitting corresponds to a first-order Chapman-Enskog expansion \cite{chapman1916law}, which is also  utilized when using a neural network to approximate kinetic equations \cite {lou2021physicsinformed, jin2021asymptoticpreserving}. Precisely, the distribution \eqref{eq:dis_f} is decomposed into two parts as
\begin{equation}
\label{eq:decom-f}
\bbf(\bx, t) = \mathcal{M}^{\rm eq}(\bx, t) + C \bbf^{\rm neq}(\bx, t),
\end{equation}
where $C$ is a problem-dependent constant. Different from the general Micro-Macro decomposition, the equilibrium and non-equilibrium parts in \eqref{eq:decom-f} are chosen as 
\begin{equation}
    \label{eq:decom-f_1}
    \begin{aligned}
 \mM^{\rm eq}_i(\bx,t)=\frac{\tilde \rho(\bm x,t)}{(\pi\tilde{T}(\bm x,t))^{3/2}}\exp\left(\frac{-(v_i-\tilde {\bm u}(\bm x,t))^2}{\tilde T(\bm x,t)}\right), \qquad 
f^{\rm neq}_i(\bm x,t)=\sqrt{\mM_i(\bm x,t)}\tilde {f}_i(\bm x,t). 
    \end{aligned}
\end{equation}
Two neural networks are utilized the equilibrium to learn $\mM^{\rm eq}$ and the non-equilibrium $\bbf^{\rm neq}$, respectively
\begin{align}
\label{eq:decom-f_2}
\tilde \rho(\bm x,t), \tilde {\bm u}(\bm x,t),\tilde T(x,t) & =NN_1(\bm x,t; \cdot), \\
\tilde {\bm f}(\bm x,t) &=NN_2(\bm x,t;\cdot). 
\end{align}
For $\mM^{\rm eq}$, the neural network outputs are only $\tilde{\rho}, \tilde{\bu}, \tilde{T}$, and then generates the distribution function according to \eqref{eq:decom-f_1}. For the distribution function, the outputs are $\tilde{f}_i$. We want to emphasize that this strategy is only inspired by the Micro-Macro decomposition to design such a neural network structure for training purposes, and there is no guarantee that $\tilde \rho(\bm x,t), \tilde {\bm u}(\bm x,t),\tilde T(x,t)$ learned is the exact density, velocity, and temperature of $\bbf$. Moreover, this splitting method is enlightened by the experiments and there is no theoretical evidence right now. We have observed that with the splitting form \eqref{eq:decom-f_1}, the learning process  will be accelerated and a lower error of the loss function will be achieved.

\begin{figure}[!hptb] 
\centering
\includegraphics[width=0.5\textwidth]{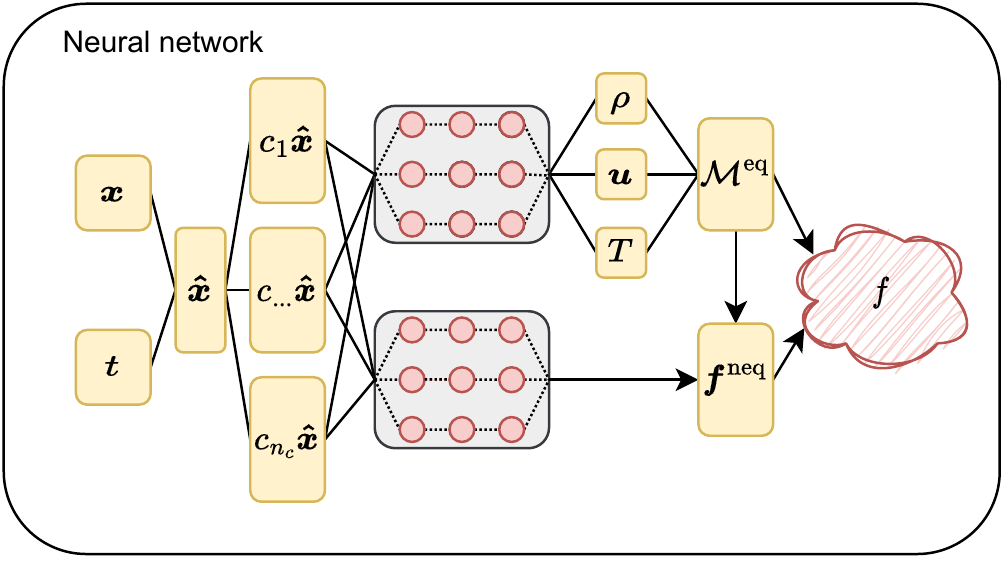} 
\caption{Network architecture. In the neural network, the inputs are the spatial space $\bx$ and time $t$, and the Monte Carlo sampling is utilized to decide the specific points. The multi-scale inputs and the Maxwellian splitting of the distribution function is adopted to improve the approximation efficiency of the network. }
\label{fig:network-main}
\end{figure}

With these two strategies, the efficiency of the neural network to approximate Boltzmann equation can be greatly improved. The total structure of the neural network is shown in Fig. \ref{fig:network-main}. Once the network is set up, the loss function will be constructed to solve the PDE. It is always built by combining the residuals for the equations as well as the initial and boundary conditions, which we will introduce in detail in the next section.

\subsection{Loss function}
\label{sec:loss}
In this section, we will discuss the specially designed loss function to solve the Boltzmann equation. Generally speaking, the loss function for a PDE should contain three parts, namely $L_{\rm IC}$ concerning the initial condition, $L_{\rm BC}$ for the boundary condition, and $L_{\rm PDE}$ to the residual of PDE
\begin{equation}
\label{eq:loss}
L_{\rm loss} =L^f_{\rm IC}+L^f_{\rm BC}+L^f_{\rm PDE}.
\end{equation}
In \eqref{eq:full_network}, the input variables are the spatial position $\bx$, and time $t$, which should be discretized first to derive the loss function. Here, Monte Carlo random sampling  is utilized, and  uniform random points are generated. Assuming the point sizes for the initial, boundary condition, and interior of the computed region are $N_{\rm IC}$, $N_{\rm BC}$, and $N_{\rm PDE}$, then the expression for the loss function can be written as 
\begin{equation}
\label{eq:loss_fun}
\begin{aligned}
&L^f_{\rm IC} = \frac{1}{N_{\rm IC}}\sum_{s=1}^{N_{\rm IC}}\mathcal{L}[\bm{f}(\bm x_s,0)-\bm{f}^0(\bm x_s)],\\
&L^f_{\rm BC} = \frac{1}{N_{\rm BC}}\sum_{s=1}^{N_{\rm BC}}\mathcal{L}[\bm{f}(\bm x_s,t_s)-\bm{f}^b(\bm x_s,t_s)],\\
&L^f_{\rm PDE} = \frac{1}{N_{\rm PDE}} \sum_{s=1}^{N_{\rm PDE}} \mathcal{L}[\bm{r}(\bm x_s,t_s)].
\end{aligned}
\end{equation}
Here, $\bbf^0(\bx_s)$ and $\bbf^b(\bx_s, t_s)$  are the discrete initial and boundary conditions, respectively, and $\br(\bx_s, t_s)$ is the residual of Boltzmann equation with 
\begin{equation}
    \label{eq:residual}
    \br(\bx, t) =\frac{\partial \bm{f}(\bm x,t)}{\partial t}+\bm{v} \cdot \nabla_{\bm x} \bm{f}(\bm  x,t)-\mathcal{Q}[\bm{f}]. 
\end{equation}
In \eqref{eq:loss_fun}, $\mathcal{L}[\cdot]$ is the distance function, and the simplest option is the $l_2$ norm
\begin{equation}
    \label{eq:l2_loss}
    \mathcal{L}_{l_2}[{\bm s}] = \| {\bm s} \|_2^2,
\end{equation}
where $\bm s$ is a vector with any length \cite{lou2021physicsinformed}. 

However, the numerical experiments show that results obtained using $l_2$-norm as a distance function are not satisfactory. When training with distance function \eqref{eq:l2_loss} directly in supervised learning, there is still a relatively large error in the macroscopic variables at the end of the training process. This may be due to that this distance function does not behave well. For example, in this distance function, each point has the same weight, while the distribution function at smaller relative velocity should be more important compared to that at larger because they have a greater impact on the macroscopic variables. Therefore, the specially designed distance function is utilized here, and we will introduce it below. 

\paragraph{\bf Macroscopic variable loss}
The macroscopic variable plays quite an important role in the simulation of Boltzmann equation, and one should try their best to obtain them more exactly. However, the numerical simulations show that with the simple $l_2$-norm distance function  \eqref{eq:l2_loss}, they can not be derived correctly. Therefore, the conserved variable \eqref{eq:dis_macro_var} is added to the distance function. Precisely, since the energy in each direction may vary greatly, the energy in each direction $E_i, i = 1, 2, 3$ instead of the total energy is utilized in the macroscopic variables, and then the macroscopic variables considered are 
\begin{equation}
    \label{eq:macro_vec}
    C[\bbf] = (\rho[\bbf], \bbm[\bbf]^T, E_1[\bbf], E_2[\bbf], E_3[\bbf])^T\in \bbR^{7}.
\end{equation}
The distance function for the conserved variable is defined as
\begin{equation}
    \label{eq:mac_loss}
    \mL^{\rm C}[\bbf] = \|C[\bbf]\|_2^2.
\end{equation}
For the loss function of these conserved variables, the IC, BC, and PDE all have their corresponding conserved variable loss function as 
\begin{equation}
    \label{eq:mac_loss_fun}
    L_S^C = \frac{1}{N_{S}} \sum_{s=1}^{N_{S}} \mathcal{L}^{\rm C}[\bm{g}(\bm x_s,t_s)],
\end{equation}
where $S = {\rm IC, BC}$ and $\rm PDE$, $\bg(\bx, t) = (\bbf(\bx, 0) - \bbf^0(\bx)), (\bbf(\bx,t) - \bbf^b(\bx, t))$ and $\br(\bx, t)$. 
Thus, the loss function is combined by two parts as  
\begin{equation}
    \label{eq:loss_fun_mac}
L_{\rm S} = L^f_{\rm S} + L^C_{\rm S}, \qquad S = {\rm IC,~BC, ~PDE}.
\end{equation}

\begin{figure}[!hptb] 
\centering 
\includegraphics[width=0.6\textwidth]{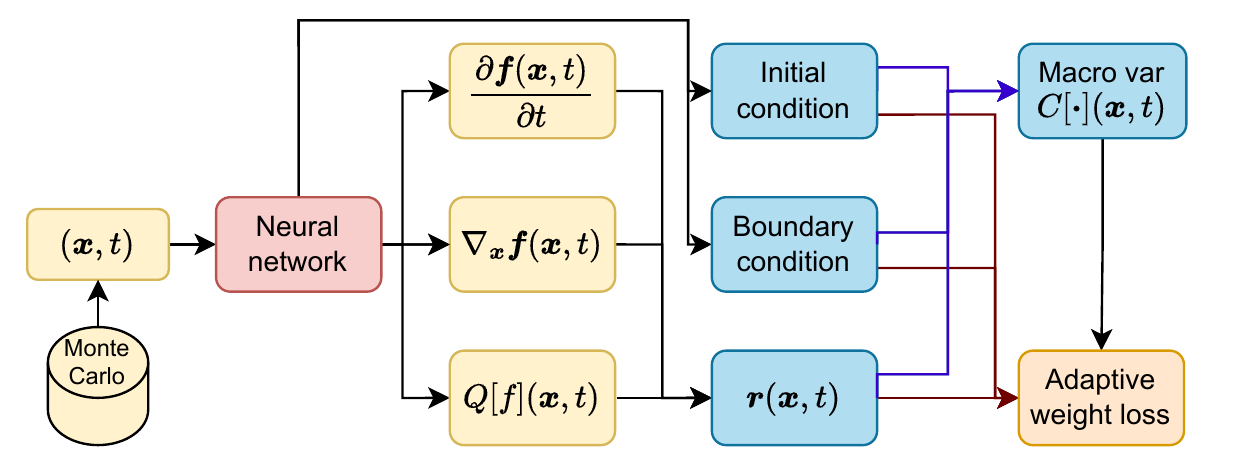} 
\caption{Sketch of the loss function.  Here, the input variables are the spatial space $\bx$ and time $t$, the specific points of which are decided by Monte Carlo. The loss function contains three parts  the initial condition, boundary condition, and residual of PDE, and  each part is composed of the microscopic part and the macroscopic variables part. } 
\label{fig:network-loss} 
\end{figure}

\paragraph{\bf Adaptive weight loss}
In the loss function \eqref{eq:loss_fun}, if the distance function \eqref{eq:l2_loss} is utilized, then each entry of $\bbf$ has the same weight. However, it is obvious that the distance function with a smaller relative velocity is more important. Therefore, how to balance the weight of the distribution function at different velocities is quite important. The lower bound constrained uncertainty weighting \cite{cipolla2018multitask, huang2021solving} is adopted here to assign the weight functions. In particular, instead of simply taking $l_2$ norm, the relative error at each microscopic velocity point $\bv_i$ is considered, and the loss function in \eqref{eq:loss_fun} is changed into 
{\small 
\begin{subequations}
    \label{eq:loss-f}
    \begin{align}
    \label{eq:loss-f_1}
       &\tilde{L}^f_{\rm IC}=\frac{1}{N_{\rm IC}}\sum_{s=1}^{N_{\rm IC}}
\sum_{i=1}^{N_v}\left(\frac{1}{(w_{\rm IC}^{f})_i+\epsilon} \left(f_i(\bm x_s, 0) - f^0_i(\bm x_s)\right)^2
+\log(1+(w_{\rm IC}^{f})_i)\right), \\\label{eq:loss-f_2}
 &\tilde{L}^f_{\rm BC}=\frac{1}{N_{\rm BC}}\sum_{s=1}^{N_{\rm BC}}
\sum_{i=1}^{N_v}\left(\frac{1}{(w_{\rm BC}^{f})_i+\epsilon} \left(f_i(\bm x_s, t_s) - f^b_i(\bm x_s, t_s)\right)^2
+\log(1+(w_{\rm BC}^{f})_i) \right), \\\label{eq:loss-f_3}
    &\tilde{L}^f_{\rm PDE}=\frac{1}{N_{\rm PDE}}\sum_{s=1}^{N_{\rm PDE}}
\sum_{i=1}^{N_v}\left(   \frac{1}{(w_{\rm PDE}^{f})_i+\epsilon} (r_i(\bm x_s,t_s))^2
+\log(1+(w_{\rm PDE}^{f})_i) \right),
    \end{align}
\end{subequations}
} where $(w_{s}^{f})_i \geqslant 0, s = {\rm IC, BC, PDE}, i = 1, \cdots N_v$ are the weights for point $\bv_i$ in different loss functions, and $\epsilon$ is a small positive number preventing division by zero. Here, we want to emphasize that different from the integral weight \eqref{eq:weight}, $(w_{s}^{f})_i$ are the parameters in the neural network, and their values are changing adaptively during the training process. Moreover, the loss function for the macroscopic variables \eqref{eq:mac_loss_fun} is also revised similarly as 
{\small 
\begin{subequations}
\label{eq:loss-prim}
\begin{align}
\label{eq:loss-prim_1}
&\tilde{L}^{C}_{\rm IC}=
\frac{1}{N_{\rm IC}} \sum_{s=1}^{N_{\rm IC}}
\sum_{i=1}^{7}
\left(   \frac{1}{(w_{\rm IC}^C)_i+\epsilon}  (C_i[\bbf(\bx_s, 0) - \bbf^0(\bx_s)])^2
+\log(1+(w_{\rm IC}^C)_i) \right), \\
\label{eq:loss-prim_2}
&\tilde{L}^{C}_{\rm BC}=
\frac{1}{N_{\rm BC}} \sum_{s=1}^{N_{\rm BC}}
\sum_{i=1}^{7}
\left( \frac{1}{(w_{\rm BC}^C)_i+\epsilon}  (C_i[\bbf(\bx_s, t_s) - \bbf^b(\bx_s, t_s)])^2
+\log(1+(w_{\rm BC}^C)_i) \right), \\\label{eq:loss-prim_3}
&\tilde{L}^{C}_{\rm PDE}=
\frac{1}{N_{\rm PDE}} \sum_{s=1}^{N_{\rm PDE}}
\sum_{i=1}^{7} \left(
\frac{1}{(w_{\rm PDE}^C)_i+\epsilon}  (C_i[\bm r(\bm x_s,t_s)])^2
+\log(1+(w_{\rm PDE}^C)_i) \right),
\end{align}
\end{subequations}
} where $(w_s^C)_i \geqslant 0, s = {\rm IC, BC, PDE}, i = 1, \cdots 7$,  are also weight parameters in the neural network. Combining \eqref{eq:loss-f} and \eqref{eq:loss-prim}, we will derive the final loss function for the Boltzmann equation with the sketch shown in Fig. \ref{fig:network-loss}
\begin{equation}
    \label{eq:fin_loss}
    L_{\rm loss} = \tilde{L}_{\rm IC} + \tilde{L}_{\rm BC} + \tilde{L}_{\rm PDE}, \qquad \tilde{L}_{S} = \tilde{L}^f_S + \tilde{L}^C_S, \qquad S = {\rm IC, ~BC, ~PDE}.
\end{equation}

Right now, we have completed the description of this network-based method for Boltzmann equation, where the strategy multi-scale input and Maxwellian splitting are utilized to improve the approximation efficiency, and the specially designed loss in Sec. \ref{sec:loss} are adopted to achieve the final results. This method is then applied to solve Boltzmann equation with BGK and quadratic collision model in Sec. \ref{sec:num}, and we call it $\rm NR$ (${\rm NR^{BGK}}$ and ${\rm NR^{Quad}}$) for short. 

However, due to the high dimensionality of the distribution function, and the complex form of the quadratic collision model, using the NR method directly is still quite expensive. Therefore, NSR to the distribution function is proposed, where the low-rank property of the distribution function is explored, with several sparse expressions introduced to greatly reduce the degree of freedom. 

\section{Neural sparse representation for BGK equation}
\label{sec:lr}
When using a neural network to approximate a discrete distribution function (DDF), most of the parameters in the network are concentrated in the last layer, which increases the computational cost and decreases the approximation efficiency. Moreover, for high-dimensional PDEs, such as Boltzmann equation, the high-order tensor obtained after discretization is one of the important sources for the huge computational cost. To reduce the number of parameters and the computational complexity, the tensor decomposition \cite{boelens2018parallel, reynolds2016randomized}, which can effectively alleviate the curse of dimensionality, is exploited for the discrete distribution function. 

In this section, we will first focus on the BGK model, and the low-rank property of the distribution function is utilized to reduce the number of parameters. We will begin with the introduction to the tensor low-rank decomposition, and then its implication in the neural representation is proposed. 

\subsection{Tensor low-rank decomposition}
\label{sec:low_rank}
Tensor decomposition represents the higher-order tensor with a series of lower-order tensors. For matrix, the second-order tensor, SVD decomposition is the most popular tool. For the higher order, there is Canonical polyadic decomposition (CPD) \cite{carroll1970analysis, harshman1970foundations} and Tucker decomposition (TD) \cite{hitchcock1927expression, tucker1966mathematical} which can be seen as a generalization of SVD from two dimensions to higher dimensions. CDP is utilized here to reduce the degrees of freedom in the discrete distribution function. Precisely, the third-order tensor CPD can be expressed as

\paragraph{\bf Canonical polyadic decomposition (CPD)}
 Assuming $\msF \in \mathbb{R}^{N_1 \times N_2 \times N_3}$ is a third-order tensor, then  $\msT \in \mathbb{R}^{N_1 \times N_2 \times N_3}$ is the CDP factorization of $\msF$, if it is the solution to the optimization problem
 \begin{equation}
    \label{eq:CPD}
    \begin{aligned}
    &\min_{\msT \in \mathbb{R}^{N_1 \times N_2 \times N_3}} \Vert \msT-\msF \Vert_F, \\
    & {\rm s.t.}   \quad  \msT_{ijk}=\sum_{r=1}^K \bU^{(1)}_{ir} \bU^{(2)}_{jr} \bU^{(3)}_{kr},
    \end{aligned}
    \end{equation}
   where $\bU^{(k)}, k = 1, 2 , 3$ are second-order tensors, with $K \in \mathbb{N}^+$. The factorization is called exact if the distance in \eqref{eq:CPD} is zero, or approximated if it is larger than $0$.

It is always expected that $K \ll N_k, k = 1, 2,3 $ in \eqref{eq:CPD} to greatly reduce the computational cost. Generally speaking, it is not easy to solve this optimization problem, and for most  tensors, we can only derive the approximated factorization. There are several classical methods for CPD, such as ALS \cite{carroll1970analysis}, et al. Tensor low-rank decomposition is already widely utilized to solve Boltzmann equation in the framework of discrete velocity method. Precisely, assuming the discrete grid points in the microscopic velocity space are chosen as \eqref{eq:dvm_points}, the obtained discrete distribution function \eqref{eq:dis_f} can be considered as a third-order tensor
\begin{equation}
\label{eq:f_tensor}
\msF(\bx, t) \in \bbR^{N_1\times N_2 \times N_3}, \qquad \msF_{l_1l_2l_3} = f_l(\bx, t), \qquad 1\leqslant l_i \leqslant N_i, \quad i = 1, 2, 3,
\end{equation}
where the one-to-one mapping of $l$ and $l_1,l_2, l_3$ is the same as that in \eqref{eq:dvm_points}. This is an $\mathcal{O}(N_v)$ tensor, which is quite memory-consuming, and also the main reason for the extensive computational effort.  CPD is utilized here to reduce the memory cost. 
The discrete distribution function \eqref{eq:f_tensor}, is a third-order tensor and can be approximated using CPD as 
\begin{equation}\label{eq:low-rank}
\msF_{ijk}(\bx, t):=f(\bx,t,v^{(1)}_i,v^{(2)}_j,v^{(3)}_k) \approx \sum_{r=1}^K \bm{P}_{ir}(\bx, t)\bm{Q}_{jr}(\bx, t)\bm{R}_{kr}(\bx,t),
\end{equation}
where
\begin{equation}
\label{eq:CPD_ele}
  \bP \in \mathbb{R}^{{N_1} \times K}, \qquad \bQ \in \mathbb{R}^{N_2\times K}, \qquad \bR \in \mathbb{R}^{N_3\times K}.
\end{equation}
For convenience, we will write  \eqref{eq:low-rank} as
\begin{equation}
    \label{eq:short_f_CPD}
    \msF = \llb \bP, \bQ, \bR \rrb.
\end{equation}
Here, we want to emphasize that $\bP, \bQ, \bR$ in \eqref{eq:low-rank} can be derived using ALS \cite{carroll1970analysis}, which is also employed in \cite{boelens2020tensor}. Since that in the framework of neural network method, the process of obtaining these components will be completed automatedly in the learning process, we will not focus on this. Moreover, the CPD of the distribution function has several useful properties, and we list them in the following lemma. 
\begin{lemma}
\label{lemma:lowrank-op}
The following operations of a tensor with CPD still keep the form of CPD.

\begin{enumerate}
\item Addition and subtraction

Assuming $\msF^{(1)}, \msF^{(2)} \in \bbR^{N_1 \times N_2 \times N_3}$ are two third-order tensors with the same size, whose CPD are  
\begin{equation}
    \label{eq:FF_CPD}
    \msF^{(1)} =  \llb \bP^{(1)}, \bQ^{(1)}, \bR^{(1)} \rrb, \qquad \msF^{(2)} = \llb \bP^{(2)}, \bQ^{(2)}, \bR^{(2)} \rrb,
\end{equation}
where 
\begin{equation}
     \bP^{(i)} \in \mathbb{R}^{{N_1} \times K^{(i)}}, \qquad \bQ \in \mathbb{R}^{N_2\times K^{(i)}}, \qquad \bR \in \mathbb{R}^{N_3\times K^{(i)}}, \qquad i = 1, 2,
\end{equation}
and there is no requirement that $K^{(1)}$ equals $K^{(2)}$. Then, it holds that 
\begin{equation}
\begin{aligned}
\msF^{(1)} \pm \msF^{(2)} & =\llb \bP^{(1)},\bQ^{(1)},\bR^{(1)} \rrb \pm 
 \llb \bP^{(2)},\bQ^{(2)},\bR^{(2)} \rrb
  \triangleq \llb \bP, \bQ, \bR \rrb,
\end{aligned}
\end{equation}
where  
\begin{equation}
\bP = [\bP^{(1)}, \pm \bP^{(2)}] \in \bbR^{N_1\times K}, \qquad  \bQ = [\bQ^{(1)}, \bQ^{(2)}] \in \bbR^{N_2 \times K}, \qquad \bR = [\bR^{(1)}, \bR^{(2)}] \in \bbR^{N_3\times K}, 
\end{equation}
with $K = K^{(1)} + K^{(2)}$. 
\item Derivatives

Assuming $\msF(x) \in \bbR^{N_1 \times N_2 \times N_3}$ is a  third-order tensor, whos CPD is 
\begin{equation}
    \label{eq:F_CPD}
    \msF(x) = \llb \bP(x), \bQ(x), \bR(x) \rrb,
\end{equation}
where 
\begin{equation}
     \bP(x) \in \mathbb{R}^{{N_1} \times K}, \qquad \bQ(x) \in \mathbb{R}^{N_2\times K}, \qquad \bR(x) \in \mathbb{R}^{N_3\times K}.
\end{equation}
Then, it holds that 
\begin{equation}
(\msF(x))'=(\llb (\bP(x),\bQ(x),\bR(x)\rrb )' \triangleq \llb \overline{\bP}, \overline{\bQ}, \overline{\bR} \rrb,
\end{equation}
where 
\begin{equation}
\begin{aligned}
    \overline{\bP} = [\bP^{\prime}(x), \bP(x), \bP(x)] \in \bbR^{N_1\times 3K},  \\
    \overline{\bQ} = [\bQ(x), \bQ^{\prime}(x), \bQ(x)]\in \bbR^{N_2\times 3K},  \\
    \overline{\bR} = [\bR(x), \bR(x), \bR^{\prime}(x)]\in \bbR^{N_3 \times 3K}.
\end{aligned}
\end{equation}
\end{enumerate}
\end{lemma}
The proof of this lemma is obvious, and we will omit it here. This lemma shows that the summation of two third-order tensors with the same size, whose CPD  ranks are  $K_1$ and $K_2$, respectively, can be written into a tensor with CPD rank $K_1 + K_2$. The derivative of a tensor with CPD rank $K$, can also be written into a tensor with CPD rank $3K$. These two properties are quite important when building loss function and in the learning process. For now, the low-rank decomposition of the distribution function is completed, and the detailed implementation in the $\nnBGK$ method to the BGK model will be introduced in the next section.

\subsection{Implementation of low-rank network}
\label{sec:imp_lr}
With the low-rank decomposition of the discrete distribution function, the degrees of freedom will be greatly reduced, and then the computational cost will decrease. The calculation of the macroscopic variables \eqref{eq:dis_macro_var} and the loss function \eqref{eq:loss-prim} will all become much easier based on CPD \eqref{eq:low-rank}.

\paragraph{Computation of the macroscopic variables}
Once the discrete distribution function is written into CDP form \eqref{eq:low-rank}, the computation of the macroscopic variables can be calculated more easily.  Define the moments as 
\begin{equation}
    \label{eq:moment}
    M_{i_1i_2i_3}[f] = \int_{\bbR^3}(v^{(1)})^{i_1}(v^{(2)})^{i_2}(v^{(3)})^{i_3} f \dd \bv. 
\end{equation}
Then, the macroscopic variables can be express by $M_{i_1i_2i_3}$. For example, the density $\rho$ corresponds to $(i_1, i_2, i_3) = (0,0,0)$, and $E_i$ in \eqref{eq:sep_E} corresponds to $2e_i$. 

Since the full discrete points $\bV$ defined in \eqref{eq:dvm_points} and the weight \eqref{eq:weight} can be treated as tensors with CPD rank one,  when approximating the discrete distribution function with CPD factorization, the moments \eqref{eq:moment} can be calculated as 
\begin{equation}
    \label{eq:CDP_macro}
    \begin{aligned}
    M_{i_1i_2i_3}[\msF] &= \sum_{l_1}^{N_1}\sum_{l_2}^{N_2}\sum_{l_3}^{N_3}
    \omega_{l_1}^{(1)} \omega_{l_2}^{(2)} \omega_{l_3}^{(3)}
(v_{l_1}^{(1)})^{i_1}(v_{l_2}^{(2)})^{i_2}(v_{l_3}^{(3)})^{i_3}
\msF_{l_1l_2l_3} \\
    & = \sum_{r=1}^K \big( \sum_{l_1}^{N_1} \omega_{l_1}^{(1)}(v_{l_1}^{(1)})^{i_1}\bm{P}_{l_1r} \big)\big( \sum_{l_2}^{N_2} \omega_{l_2}^{(2)}(v_{l_2}^{(2)})^{i_2}\bm{Q}_{l_2r} \big)\big( \sum_{l_3}^{N_3} \omega_{l_3}^{(3)}(v_{l_3}^{(3)})^{i_3}\bm{R}_{l_3r} \big).
\end{aligned}    
\end{equation}
We can find that the computational cost for the moments is reduced from $\mathcal{O}({N_v})$ to $\mathcal{O}(K{N_v}^{\frac13})$ with the assumption $N_1, N_2, N_3$ at the same order, which is greatly reduced especially when $K$ is small.

\paragraph{\bf Approximation of BGK collision term}
To approximate the BGK collision term, we mainly need to approximate the equilibrium $\mM$ in the low-rank form. Its CPD factorization is 
\begin{equation}
    \label{eq:CPD_M}
    \msM_{ijk}(\bx, t) :=\mM(\bx, t, v_i^{(1)}, v_j^{(2)}, v_k^{(3)}) = {\bf M}^{(1)}_i {\bf M}^{(2)}_j{\bf M}^{(3)}_k,
\end{equation}
where
\begin{equation}
    \label{eq:CPD_M_ele}
  {\bf M}^{(l)} \in \bbR^{N_l}, \qquad   {\bf M}^{(l)}_i = \frac{\rho^{1/3}}{\sqrt{2 \pi T}}\exp\left(-\frac{
(v_i^{(l)} - u^{(l)})^2 }{2T}\right), \qquad i = 1, \cdots, N_l, \quad 
l = 1, 2,3.
\end{equation}
It shows that the CPD rank of the Maxwellian is one, which will make the approximation to the BGK collision term much easier.

\paragraph{\bf Loss function}
In this part, we will introduce how to build the adaptive weight loss function \eqref{eq:fin_loss} with the CPD of the discrete distribution function. For the discrete distribution function part \eqref{eq:loss-f}, assume the adaptive weight $(w_s^f)_l, l = 1, \cdots N_v, s = {\rm IC, BC, PDE}$  in \eqref{eq:loss-f} is rank one, which means that 
\begin{equation}
    \label{eq:cpd_weight}
    (w_s^f)_l = (w_s^f)^{(1)}_{l_1}(w_s^f)^{(2)}_{l_2}(w_s^f)^{(3)}_{l_3},
\end{equation}
where the mapping of $l$ and $l_1,l_2,l_3$ is the same as that in \eqref{eq:weight}. Based on this, the loss function for the distribution function can be written into the matrix form. Let \eqref{eq:loss-f_3} as an example. With the properties in Lemma \ref{lemma:lowrank-op}, $\bm{r}(x, t)$ in \eqref{eq:residual} can be CP decomposed. Suppose its CPD factorization is 
\begin{equation}
    \label{eq:CPr}
    \mathscr{R}(\bx, t) = \llb {\bf R}(\bx,t), {\bf S}(\bx,t), {\bf T}(\bx,t)\rrb.
\end{equation}
Let the first part of \eqref{eq:loss-f_3} be 
\begin{equation}
    \label{eq:CPD_loss_f_0}
    F(\bx_s, t_s) \triangleq \sum_{i = 1}^{N_v} \frac{r_i(\bx_s, t_s)^2}{(w_{\rm PDE}^f)^2 + \epsilon},
 \end{equation}
with \eqref{eq:cpd_weight} and \eqref{eq:CPr}, it can be expressed into the CPD form as 
\begin{equation}
    \label{eq:CPD_loss_f_1}
    F(\bx_s,t_s) \approx  \sum_{l_1l_2l_3}\left(\frac{\mathscr{R}_{l_1l_2l_3}(\bx_s, t_s)}{\msW_{l_1l_2l_3} + \epsilon}\right)^2
     \approx \sum_{l_1l_2l_3} \left(
    \sum_{r = 1}^{\overline{K}} \left( \frac{{\bf R}_{l_1r}(\bx_s, t_s)}{w_{l_1}^{(1)} + \epsilon}\right)
    \left(\frac{{\bf S}_{l_2r}(\bx_s, t_s)}{w_{l_2}^{(2)} + \epsilon}\right)
    \left(\frac{{\bf T}_{l_3r}(\bx_s, t_s)}{w_{l_3}^{(3)} + \epsilon}\right)
    \right)^2,
\end{equation}
where $\sum\limits_{l_1l_2l_3} = \sum\limits_{l_1=1}^{N_1}\sum\limits_{l_2=1}^{N_2}\sum\limits_{l_3=1}^{N_3}$, and $\overline{K}$ is the CPD rank of $\mathscr{R}$. Here, the subscripts $f$ and $\rm PDE$ of $(w_s^f)^{(i)}_{l_i}$ are omitted. 
Assuming there exists a third-order tensor $\mathscr{G}$, whose CPD factorization satisfies that 
\begin{equation}
    \label{eq:tensor_G}
    \mathscr{G}_{l_1l_2l_3}(\bx_s, t_s) :=  \sum_{r = 1}^{\overline{K}} \left( \frac{{\bf R}_{l_1r}(\bx_s, t_s)}{w_{l_1}^{(1)} + \epsilon}\right)
    \left(\frac{{\bf S}_{l_2r}(\bx_s, t_s)}{w_{l_2}^{(2)} + \epsilon}\right)
    \left(\frac{{\bf T}_{l_3r}(\bx_s, t_s)}{w_{l_3}^{(3)} + \epsilon}\right),
\end{equation}
the adaptive weight loss function  \eqref{eq:loss-f_3} is changed into
\begin{equation}
    \label{eq:loss-f_3_CP}
    \hat{L}^f_{\rm PDE}=\frac{1}{N_{\rm PDE}}\sum_{s=1}^{N_{\rm PDE}} (\|\mathscr{G}(\bx_s, t_s)\|_{F})^2
    + \sum_{l_1l_2l_3}\log(1+(w_{l_1}^{(1)} w_{l_2}^{(2)}w_{l_3}^{(3)})),
\end{equation}
where $\|\mathscr{G}\|_{F}$ is the Frobenius norm of $\mathscr{G}$. Moreover, the loss function \eqref{eq:loss-f_1} and \eqref{eq:loss-f_2} can be revised similarly. 
To calculate the F-norm of a third-order tensor, we present the lemma below. 
\begin{lemma}\label{lemma:lowrank-square}
Assuming $\msF \in \bbR^{N_1\times N_2 \times N_3}$ is a third-order tensor, which has a CPD factorization of rank $K$ as 
\begin{equation}
    \label{eq:CPD_F}
    \msF = \llb {\bf P}_{l_1r}, {\bf Q}_{l_2r}, {\bf R}_{l_3r}\rrb,
\end{equation}
then it holds for its F-norm that 
\begin{equation}
    \label{eq:Fnorm_F}
    \|\msF\|_F^2 = \sum_{r =1}^{K}\sum_{r' = 1}^K {\bf H}_{rr'},
\end{equation}
where ${\bf H} =({\bf P}^T {\bf P})\times ({\bf Q}^T {\bf Q}) \times ({\bf R}^T {\bf R})$, and $\times$ is defined in \eqref{eq:times}.
\end{lemma}
{\newcommand{\proofname}{Proof of Lemma \ref{lemma:lowrank-square}}
\begin{proof}
From the definition of the F-norm and the CPD factorization \eqref{eq:CPD_F}, it holds that 
\begin{equation}
\begin{aligned}
 \Vert \msF \Vert _F^2&=\sum_{l_1l_2l_3} \left(\sum_{r=1}^K {\bf P}_{l_1r}{\bf Q}_{l_2r}{\bf R}_{l_3r}\right)^2= \sum_{l_1l_2l_3} \left[\big(\sum_{r=1}^K {\bf P}_{l_1r}{\bf Q}_{l_2r}{\bf R}_{l_3r}\big)  \big(\sum_{r'=1}^K {\bf P}_{l_1r'}{\bf Q}_{l_2r'}{\bf R}_{l_3r'}\big)\right]. 
\end{aligned}
\end{equation}
Changing the order of summation, we derive \eqref{eq:Fnorm_F}, and the proof is completed. 
\end{proof}
}

For the macroscopic variables part of the adaptive weight function \eqref{eq:loss-prim}, the calculation of the macroscopic variables in \eqref{eq:macro_vec} can be derived using the CPD form of the moment \eqref{eq:CDP_macro}, where the computational cost will be reduced from $\mathcal{O}(N_v)$ to $\mathcal{O}(KN_v^{\frac13})$. For now, the implementation of CPD on the neural representation for the BGK model is completed, and we call it neural sparse representation method ($\nnLR$  for short). When solving BGK model using $\nnLR$, the CPD factorization of the distribution function \eqref{eq:low-rank} is utilized instead of the discrete distribution function \eqref{eq:dis_f} as the outputs, in which case the parameter number of outputs will be greatly reduced. Moreover, compared to $\nnBGK$, whose computational cost of the loss function is $\mO(N_sN_v)$, where $N_s = N_{\rm IC} + N_{\rm BC} + N_{\rm PDE}$ is the total grid number in the spatial space, that in $\nnLR$ is $\mO(N_sN_v^{\frac13}K^2+N_v)$. Moreover, memory usage can be decreased by a large amount as well.

However, $\nnLR$ can not be extended to the full Boltzmann equation, since the quadratic collision term can not be expressed in the CPD form. A data-driven quadratic collision model is proposed to reduce the computational cost of the full Boltzmann equation, which we will introduce in the next section. 
\section{Neural sparse representation for Boltzmann equation with quadratic collision}
\label{sec:svd}
Due to the complex form of the quadratic collision model,  $\nnLR$ can not be utilized directly for the full Boltzmann equation. To build the neural sparse representation of the quadratic collision model, a series of basis vectors for the collision term are obtained from the BGK solution in the framework of the discrete system \eqref{eq:vec_dis_Bol}, with which the degree of freedom for the collision term can be  greatly reduced. 

\subsection{Approximating the quadratic collision term}
\label{sec:cal_coe}
Inspired by the spectral methods, we want to find a series of basis vectors to approximate the quadratic collision term in the framework of the discrete system \eqref{eq:vec_dis_Bol}. We begin from the continuous form. Supposing there are two sets of standard basis functions $\{g_r(\bm v)\}_{r=1}^{N_a},\{h_r(\bm v)\}_{r=1}^{N_b}$ for the distribution function and the collision term, respectively. Thus, it holds that 
\begin{subequations}
\label{eq:quad_basis}
\begin{align}
 f(\bm x,\bm v,t) &\approx \tilde{f}(\bx, \bv, t) = \sum_{r=1}^{N_a} \tilde{f}_r(\bm 
 x,t)g_r(\bm v), \label{eq:exp_f} \\  \mathcal{Q}[f](\bm x,\bm v,t) &\approx \tilde{\mQ}(\bx, \bv, t) =  \sum_{r=1}^{N_b} \tilde{Q}_r(\bm 
 x,t)h_r(\bm v). \label{eq:exp_Q}
\end{align}
\end{subequations}
By the orthogonality  of the basis functions, the expansion coefficients are obtained as 
\begin{equation}
    \label{eq:exp_coe}
    \tilde{f}_r(\bx, t) = \Braket{f g_r}, \qquad \tilde{Q}_r(\bx, t) = \Braket{\mQ[f], h_r} = \Braket{\mQ(f, f) h_r},
\end{equation}
where $\Braket{ \cdot}$ is defined in \eqref{eq:inner_fun}. Substituting \eqref{eq:exp_f} into \eqref{eq:exp_coe},  the coefficient $Q_r$ can be obtained as 
\begin{equation}
    \label{eq:exp_Q_1}
\tilde{Q}_r = \Braket{\mQ(\tilde{f} \tilde{f})h_r} = \sum_{i=1}^{N_a}\sum_{j = 1}^{N_a} \tilde{f}_i \tilde{f}_j \Braket{\mQ(g_i, g_j) h_r}.
\end{equation}
Defining the third-order tensor $\msK \in \bbR^{N_a \times N_a \times N_b}$ as 
\begin{equation}
    \label{eq:coe_ker}
    \msK_{ijr} = \Braket{\mQ(g_i, g_j) h_r}, 
\end{equation}
the quadratic collision term can be approximated as 
\begin{equation}
    \label{eq:exp_Q_2}
    \tilde{\mQ}(\bx, t) = \sum_{i = 1}^{N_a}\sum_{j =1}^{N_a}\sum_{r = 1}^{N_b}\msK_{ijr}(\bx, t) \tilde{f}_i(\bx, t) \tilde{f}_j(\bx, t) h_r(\bv).
\end{equation}
Here, $\msK_{ijr}$ can be treated as the expansion coefficients of the collision kernel, which is independent of the distribution function, and can be pre-computed offline once the basis functions are decided. 
\begin{remark}
    In the Hermite spectral method \cite{wang2019approximation}, the distribution function and the quadratic collision term are also approximated in the form \eqref{eq:quad_basis}, where both basis functions $g_r, h_r$ are the Hermite functions, and the expansion order is the same as $N_a = N_b$. For the Fourier spectral method \cite{mouhot2006fast}, the basis functions are the trigonometric functions. 
\end{remark}

This inspires us that if we can find a series of basis vectors in the framework of discrete system \eqref{eq:vec_dis_Bol} which can approximate the collision term with high efficiency, then the computational cost for the quadratic collision model can be greatly reduced. Assume we already have a series of basis vectors for the discrete distribution function and the quadratic collision term as 
$\bg_r \in \bbR^{N_v}, r = 1, \cdots N_a$, and $\bh_r\in\bbR^{N_v}, r = 1, \cdots, N_b$, in the whole spatial and time space. Then, the discrete distribution function $\bbf$ and discrete collision term $\mQ[\bbf]$ can be approximated as 
\begin{equation}
    \label{eq:exp_dis}
   \bbf_{s} \triangleq \bbf(\bx_s, t_s) \approx \sum_{r = 1}^{N_a} \hat{f}_r \bg_r, \qquad \mQ[\bbf_s]\approx \sum_{r = 1}^{N_b} \hat{q}_r \bh_r, \qquad \hat{f}_r, \hat{q}_r \in \bbR,
 \end{equation}
with
\begin{equation}
    \label{eq:exp_dis_coe}
    \hat{f}_r = \braket{\bbf_s, \bg_r}, \qquad \hat{q}_r = \braket{\mQ[\bbf_s], \bh_r},
\end{equation}
where $\braket{\cdot, \cdot}$ is the inner-product of two vectors. Without loss of generality, we scale the weight \eqref{eq:weight} to $1$ here. With the similar process of \eqref{eq:exp_Q_1}, the expansion coefficient $q_r$ is calculated as 
\begin{equation}
    \label{eq:exp_dis_coe_1}
    \hat{q}_r = \braket{\mQ(\bbf_s, \bbf_s), \bh_r} = \sum_{i = 1}^{N_a}\sum_{j = 1}^{N_a}\hat{f}_i \hat{f}_j\braket{\mQ(\bg_i, \bg_j), \bh_r} 
    \triangleq \sum_{i = 1}^{N_a}\sum_{j = 1}^{N_a}\hat{f}_i \hat{f}_j \hat{\msK}_{ijr},
\end{equation}
with 
\begin{equation}
    \label{eq:kerner_dis}
    \hat{\msK}_{ijr} = \braket{\mQ(\bg_i, \bg_j), \bh_r} 
\end{equation}
the expansion coefficients of the collision kernel under the discrete basis function. Similarly, $\hat{\msK}_{ijr}$ depends only on the basis vectors $\bg_r$ and $\bh_r$, which only need to compute once. The discrete collision term can be approximated as 
\begin{equation}
    \label{eq:exp_dis_Q}
    \mQ[\bbf_s] \approx \sum_{i = 1}^{N_a}\sum_{j = 1}^{N_a}\sum_{r = 1}^{N_b}  \hat{\msK}_{ijr} \hat{f}_i \hat{f}_j \bh_r.
\end{equation}
This shows that the more efficient the basis vectors are, the smaller the degree of freedom $N_a$ and $N_b$ are. Therefore, the last problem we have is finding the proper basis vectors.

\subsection{Choosing data-driven basis vectors}
\label{sec:choose_basis}
 Rather than producing a large amount of data and then learning the bases from them, it seems a shortcut to derive the bases from rough solutions to the same problem. To find proper basis vectors for the discrete Boltzmann equation, the main idea here is to build the data-driven-based basis function using the BGK solution, with which the discrete distribution function is well represented and the quadratic collision term is constructed based on it.

 The solution to BGK model can be treated as  an approximation to the full Boltzmann equation, we will derive a series of the basis vectors from the discrete distribution function of the BGK solution. Precisely, supposing the discrete distribution function derived for the BGK model \eqref{eq:vec_dis_Bol} is $\bbf^{\rm BGK}(\bx_s, t_s), s = 1, \cdots N_s$, we write it into the matrix form as 
\begin{equation}
    \label{eq:sol_BGK}
    \bfA^{\rm BGK} = [\bbf_1, \cdots, \bbf_{N_s}] \in \bbR^{N_v \times N_s}, \qquad \bbf_s = \bbf^{\rm BGK}(\bx_s, t_s).
\end{equation}
The corresponding discrete quadratic collision model of these distribution functions can be written as 
\begin{equation}
    \label{eq:sol_quad}
    \bfQ^{\rm BGK} = [\mQ[\bbf]_1, \cdots, \mQ[\bbf]_{N_s}] \in \bbR^{N_v \times N_s}, \qquad \mQ[\bbf]_s = \mQ[\bbf](\bx_s, t_s).
\end{equation}
Then, performing truncated singular value decomposition (SVD) on $\bfA^{\rm BGK}$ and $\bfQ^{\rm BGK}$, it holds that 
\begin{equation}
    \label{eq:svd}
    \begin{aligned}
&    \bfA^{\rm BGK} \approx {\bf G} {\bf S}_1 {\bf V}_1^T, \qquad {\bf G} \in \bbR^{N_v \times N_a}, \quad {\bf S}_1 \in \bbR^{N_a \times N_a}, \quad {\bf V}_1^T \in \bbR^{N_a \times N_s}, \\
 &    \bfQ^{\rm BGK} \approx {\bf H} {\bf S}_2 {\bf V}_2^T, \qquad {\bf H} \in \bbR^{N_v \times N_b}, \quad {\bf S}_2 \in \bbR^{N_b \times N_b}, \quad {\bf V}_2^T \in \bbR^{N_b \times N_s},
    \end{aligned}    
\end{equation}
where $N_a$ and $N_b$ are the truncated order, which can be unequal. 

\begin{remark}
It is obvious that the larger $N_a, N_b$ are, the more accurate the approximation will be, but also the greater the computational consumption is. The truncation error utilized here is 
\begin{equation}
    \label{eq:svd_error}
e=\frac{\lVert \bfA^{\rm BGK} -{\bf G} {\bf S}_1 {\bf V}_1^T \rVert_F}{\lVert \bfA^{\rm BGK} \rVert_F}.    
\end{equation}
A moderately sized quantity is chosen to balance the truncation error and the computational cost. In the simulation, the number is set as $N_a = N_b = 40$, in which case the truncation error \eqref{eq:svd_error} is at the magnitude of $\mO(10^{-3})$. 
\end{remark}
Then, the orthogonal matrices
\begin{equation}
    \label{eq:basis}
    {\bf G} = [\bg_1, \cdots \bg_{N_a}], \qquad {\bf H} = [\bh_1, \cdots \bh_{N_b}] 
\end{equation}
satisfy that 
\begin{equation}
    \label{eq:orth}
    \langle \bg_i, \bg_j \rangle = \delta_{ij}, \quad i, j = 1, \cdots N_a, \qquad \langle \bh_i, \bh_j \rangle = \delta_{ij}, \quad i,j = 1, \cdots, N_b. 
\end{equation}
Consequently, we choose $\bg_i, i = 1, \cdots, N_a$ and $\bh_j, j = 1, \cdots N_b$ as the basis vectors of the discrete distribution function, and the discrete quadratic collision term, respectively. 
\paragraph{\bf Conservation basis}
Since the quadratic collision term keeps the conservation of density, momentum, and energy, it is expected that the reduced collision model still maintains this property. To achieve this, we let each basis vector of the collision term have the property below
\begin{equation}
    \label{eq:con_basis}
    \rho[\bh_j] = 0, \qquad \bm{m}[\bh_j] = 0, \qquad E[\bh_j] = 0, \qquad \forall~j = 1, \cdots, N_b.
\end{equation}
Precisely, define the matrix $\bf M$ as 
\begin{equation}
    \label{eqdis_macro}
    {\bf M} = [1, \bV, |\bV|^2]^T \in \bbR^{N_v\times 5}.
\end{equation}
First, orthogonalize $\bf M$ into $\tilde{\bf M}$, which satisfies 
\begin{equation}
    \label{eq:M}
    (\tilde{\bf M})^T {\bf M} = {\bf I}, \qquad \Span(\bf M) = \Span(\bf \tilde{M}). 
\end{equation}
Let 
\begin{equation}
    \label{eq:bar_H}
    \overline{\bf H} = {\bf H} - \tilde{\bf M} \tilde{\bf M}^T {\bf H}, 
\end{equation}
and it is easy to verify that $\overline{\bf H}$ satisfies \eqref{eq:con_basis}. Finally, re-orthogonalize $\overline{\bf H}$, and we obtain the final set of basis vectors 
\begin{equation}
    \label{eq:fin_H}
    \tilde{\bf H} = \mathrm{Ortho}(\overline{\bf H}).
\end{equation}
\begin{remark}
    The process to obtain $\tilde{\bf H}$ from $\bf H$ is completed by the scipy.linalg.orth function in scipy \cite{2020SciPy-NMeth}.
    
 In the learning process, since the collision kernel \eqref{eq:kerner_dis} is pre-computed, the main computational cost is to obtain the collision term \eqref{eq:exp_dis_Q}, which is $\mO(N_a^2N_b+(N_a+N_b)N_v)$. Compared to the initial cost of the DVM method $\mO(MN_v^2)$ where $M$ is the number of points on the unit-sphere, this is greatly reduced. 
\end{remark}

In the simulation, the method proposed in Sec. \ref{sec:lr} is first adopted to obtain the BGK solution, which is then utilized to form the information matrix \eqref{eq:sol_BGK}. With this matrix, the series of basis vectors for the discrete distribution function and the quadratic collision term are derived. These two sets of basis vectors only need to derive once and are fixed during the learning process, which is also independent of the spatial and time variables. 

This method presents a neural sparse representation for the quadratic collision term, and we call it $\nnLA$ for short. For completeness, the algorithms to derive \eqref{eq:kerner_dis} and \eqref{eq:exp_dis_Q} are presented in the next section, where all the operations are finished by simple calculation between matrix and vectors, which can also speed up this method. 

\subsection{Algorithm for the quadratic collision model}
\label{sec:algorithm}

In the implementation, we use $\nnBGK$ to obtain the BGK solution. Once the BGK solution $\bbf^{\rm BGK}$ is obtained, Alg. \ref{alg:svd-kernel} will be utilized to derive the two sets of basis vectors and to obtain the collision kernel \eqref{eq:kerner_dis}. Then, Alg. \ref{alg:svd-collision} will be used in $\nnLA$ to compute the quadratic collision term \eqref{eq:exp_dis_Q}.

\begin{algorithm}[!hptb]
\caption{Obtain the data-driven collision kernel}\label{alg:svd-kernel}
\KwData{${\bf A}^{\rm BGK}=[\bbf_1,...,\bbf_{N_s}]\in \mathbb{R}^{N_v\times N_s }, {\bf M}=[\bm{1},\bV,|\bV|^2]^T\in \mathbb{R}^{N_v\times 5}$}
\KwResult{$\hat{\msK}\in \mathbb{R}^{N_a\times N_a\times N_b }$, ${\bf G}=[\bm{g}_1,...,\bm{g}_{N_a}]\in \mathbb{R}^{N_v\times N_a }$, ${\bf H}=[\bm{h}_1,...,\bm{h}_{N_b}]\in \mathbb{R}^{N_v\times N_b}$}
\For{$i=1:N_s$}
{$\bm{q}_i=\mathcal{Q}[\bbf_i]$}
${\bf Q}^{\rm BGK} = [\bm{q}_1,..., \bm{q}_{N_s}]$\;

${\bf G},{\bf S_1}, {\bf V}_1^T= \rm{truncatedSVD}({\bf A}^{\rm BGK})$\; 

${\bf H},{\bf S}_2, {\bf V}_2^T= \rm{truncatedSVD}({\bf Q}^{\rm BGK})$\;

Obtain $\overline{\bf H}$ through \eqref{eq:M}, \eqref{eq:bar_H}, and \eqref{eq:fin_H}\;

\For{$k=1:N_b$}{
\For{$i=1:N_a$}{
\For{$j=1:N_a$}{
{$\hat{\msK}_{ijk}=\Braket{\mathcal{Q}(\bm{g}_i, \bm{g}_j), \bm{\bar{h}}_k }$}
}
}
}
\end{algorithm}

\begin{algorithm}[!hptb]
\caption{Collision term for a given distribution function}\label{alg:svd-collision}
\KwData{$\hat{\msK}\in \mathbb{R}^{N_a \times N_a\times N_b }$, ${\bf G}\in \mathbb{R}^{N_v\times N_a}$, ${\bf H}\in\mathbb{R}^{N_v\times N_b}$,$\bbf\in \mathbb{R}^{N_v}$}
\KwResult{${\bf Q}\in \mathbb{R}^{N_v}$}
$\bm{a}={\bf G}^T \bbf$\;
\For{$k=1:N_b$}{
$q_k=\sum\limits_{i=1}^{N_a}\sum\limits_{j=1}^{N_a} \hat{\msK}_{ijk} \bm{a}_i \bm{a}_j $
}
${\bf Q}={\bf H} \bm{b}$\;
\end{algorithm}

\section{Numerical Experiment}
\label{sec:num}
In this section, several numerical examples are studied to validate the numerical methods proposed in this work, where the spatially 1-dimensional wave and Sod tube problem, and spatially 2-dimensional periodic problem are tested. In each example, 
the Adam optimizer with learning rate $\eta_0 = 0.005$ and the cosine annealing learning rate decay algorithm \cite{loshchilov2017sgdr} are utilized. For the cosine annealing learning rate decay algorithm, the learning rate at $i$-th step is 
\begin{equation}
    \label{eq:cosine}
\eta_{i}=\frac{1}{2} \eta_{0}\left (1+ \cos\left(\frac{i}{T_{\max}}\right)\right ).    
\end{equation}
This optimizer is widely adopted in learning PDEs with neural network \cite{lou2021physicsinformed,jin2021asymptoticpreserving}, and we refer \cite{kingma2014adam} for more details.  All the tests are performed on a machine with Intel(R) Xeon(R) Gold 6240 and 4 Tesla V100 SXM2 16GB. Unless otherwise specified, all the experiments are conducted under the MindSpore\footnote{https://www.mindspore.cn} and code will be available online\footnote{https://gitee.com/mindspore/mindscience/tree/master/MindFlow/applications/physics\_driven/boltzmann}. At the end of this section, some tentative work on the efficiency of the NSR method ($\nnLR$ and $\nnLA$) with transfer learning is studied.

\subsection{1D wave problem}
\label{subsec:wave1d}
In this section, the 1D3V wave problem with periodic boundary condition is studied. The initial condition is Maxwellian with the macroscopic variables as below 
\begin{equation}
\rho(x)=1+0.5\sin(2\pi x),\qquad 
\bm{u}(x)=0,\qquad 
T(x)=1+0.5\sin(2\pi x+0.2),
\end{equation}
with the spatial space $x \in [-0.5, 0.5]$. This is a problem with the smooth initial condition, and a similar one is also studied in \cite{e2018deep,li2022learning}. 

\begin{table}[!htb]
    \centering
    \def\arraystretch{1.5}
\scalebox{0.85}{
{\footnotesize
\begin{tabular}{c|cc}
\multirow{3}{*}{neural network} & layer number         & $5$          \\ & neurons              & $80$ \\
& steps                & $10,000$   \\ \hline
 \multirow{3}{*}{ optimizer}  & method            & Adam         \\ 
  & learning rate & 0.005        \\
  &  decay algorithm   & cosine annealing  \\ \hline 
 \multirow{3}{*}{sampling points}  & $N_{\rm IC}$         & $100$  \\
 & $N_{\rm BC}$         & $200$                        \\
 & $N_{\rm PDE}$        & $500$                \\ 
 \hline 
\multirow{3}{*}{computational parameters}  
 & time  & $t \in [0, 0.1]$ \\
 & Knudsen number ($\Kn$) & $0.01, 0.1, 1.0$ \\
 & microscopic velocity space     & $\bv \in [-10,10]^3$  \\
 & grid number & $24\times 24 \times 24$
\end{tabular}
 }
}
\caption{(1D wave problem in Sec. \ref{subsec:wave1d}) Parameters of the NR/NSR methods.} 
\label{tab:ex1_para}
\end{table}
In the simulation, the BGK and  quadratic collision term with Knudsen number $\Kn = 0.01, 0.1$ and $1.0$ is considered. 
The time region is $t \in [0, 0.1]$. The computational domain in the microscopic velocity space is $[-10, 10]^3$, with the grid number $24 \times 24 \times 24$. In the neural representation, both networks of $\mM^{\rm eq}$ and \fneq consist of a $5$ layer fully connected network. Each layer has $80$ neurons. The parameters in \eqref{eq:loss-prim} are set as $N_{\rm IC} = 100$, $N_{\rm BC} = 200$, and $N_{\rm PDE} = 500$. This means that we randomly sample $100$ points in $x \in [-0.5, 0.5]$ at $t = 0$ for the initial condition, $200$ points in $x = 0$ and $t \in [0, 0.1]$ for the boundary condition, and $500$ points in the spatial space and time region $x\times t \in [-0.5, 0.5] \times [0, 0.1]$. The total training step is $10, 000$. The detailed parameters are listed in Tab. \ref{tab:ex1_para}.

\begin{figure}[!hptb] 
  \centering 
  \subfloat[$\Kn = 0.01, t = 0$]{
      \includegraphics[width=0.25\textwidth]{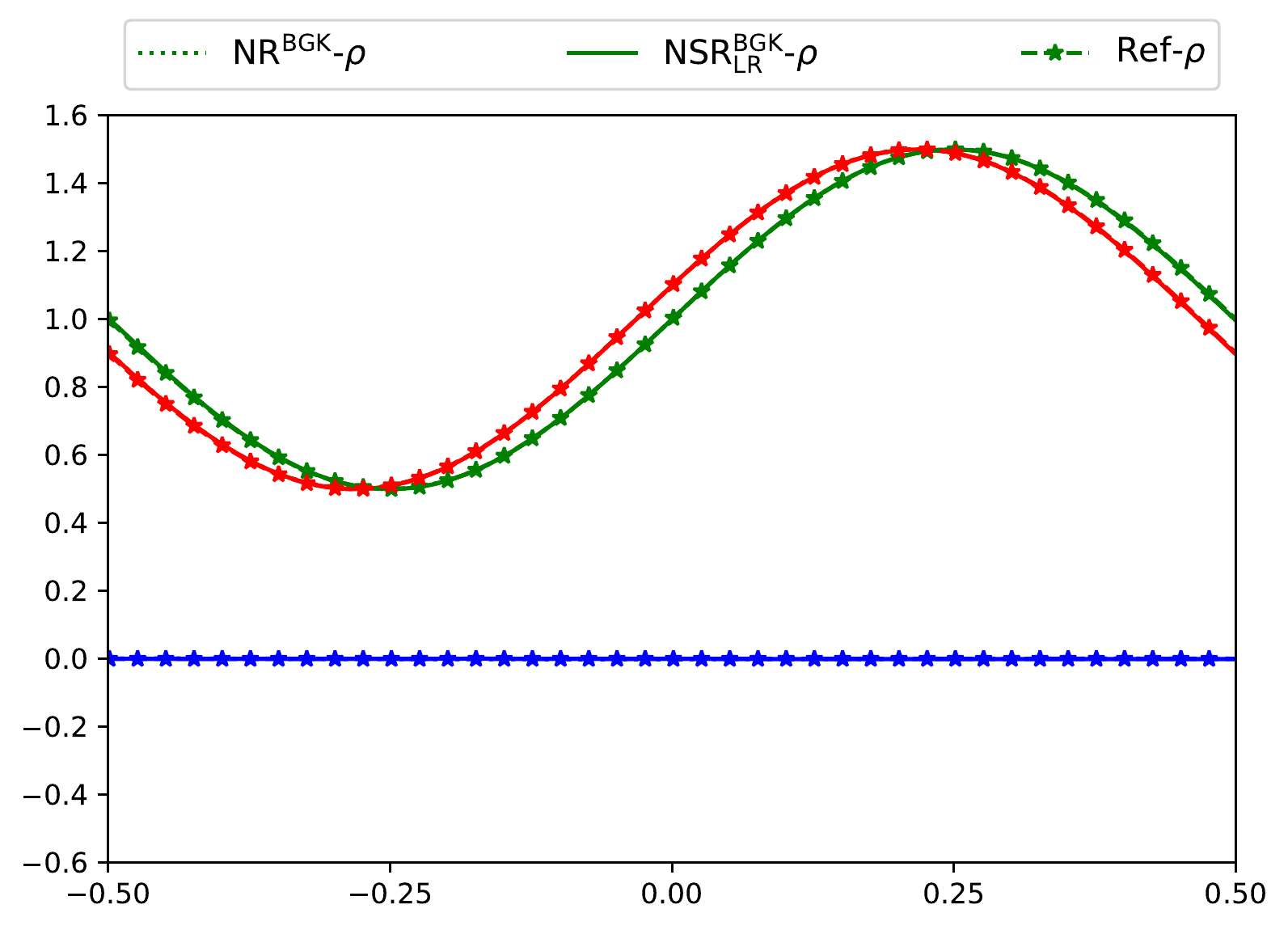} 
  } \hfill
  \subfloat[$\Kn = 0.1, t = 0$]{
      \includegraphics[width=0.25\textwidth]{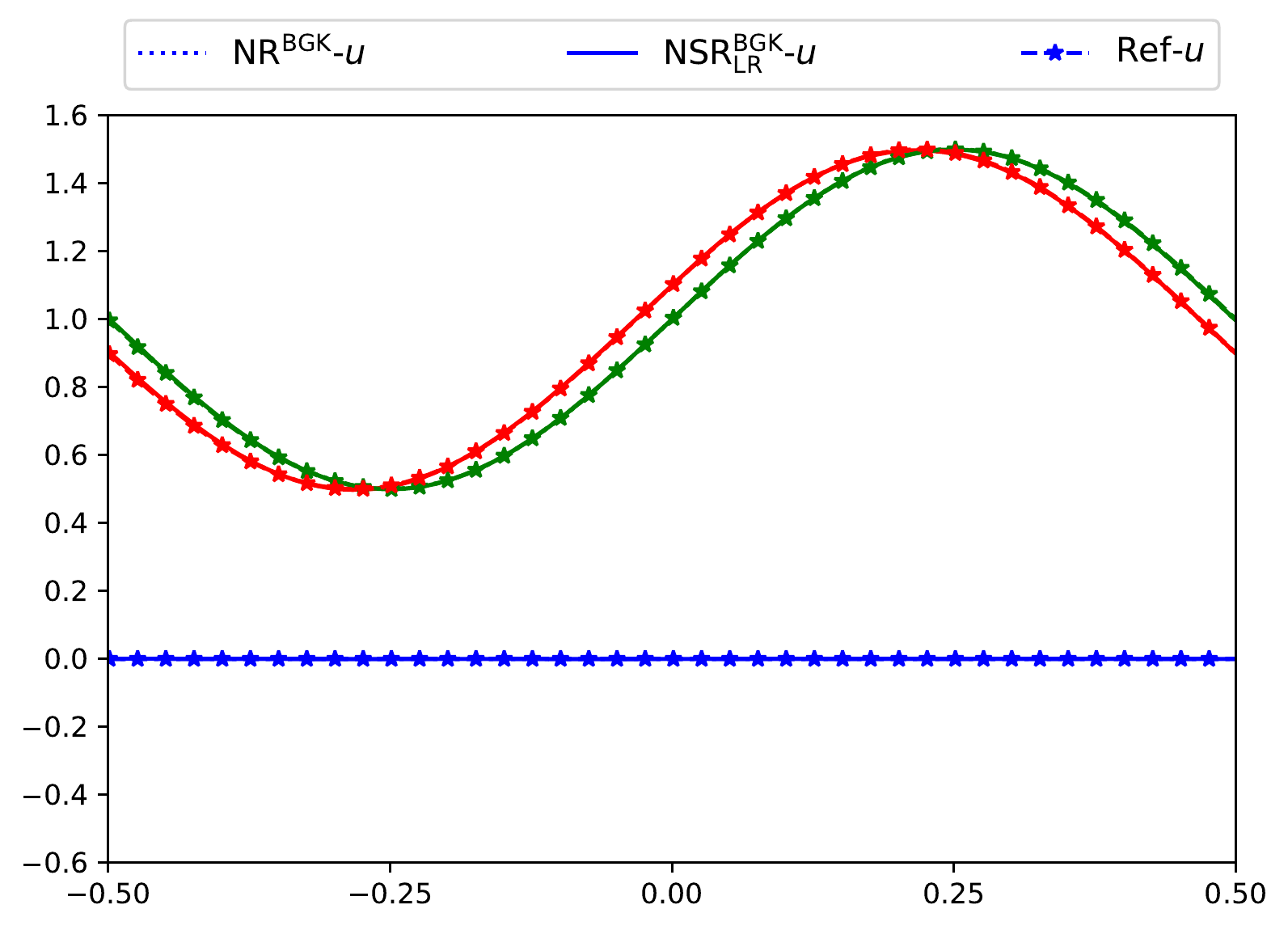} 
  } \hfill
  \subfloat[$\Kn = 1.0, t = 0$]{
      \includegraphics[width=0.25\textwidth]{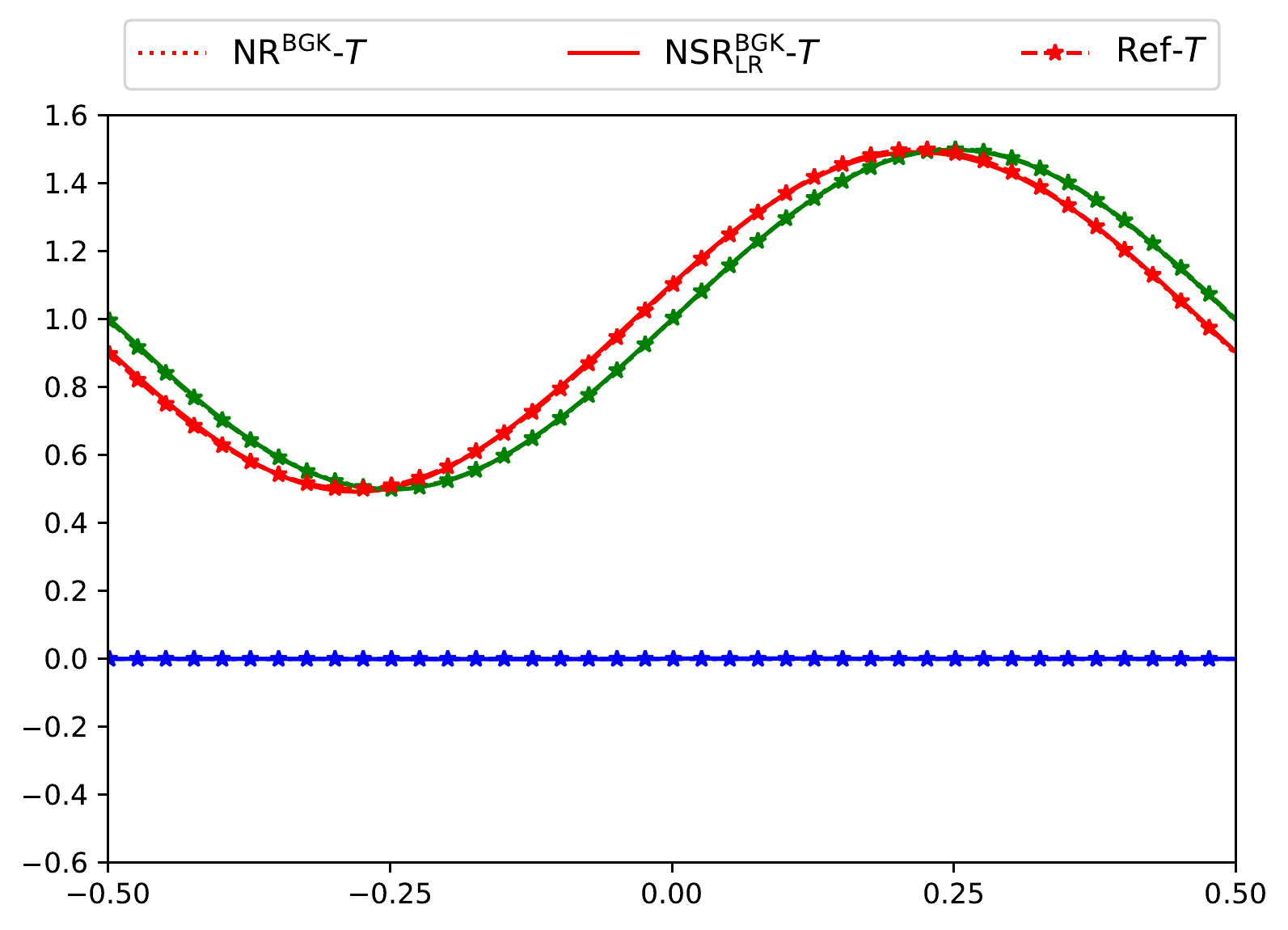} 
  } \\
   \subfloat[$\Kn = 0.01, t = 0.1$]{
      \includegraphics[width=0.25\textwidth]{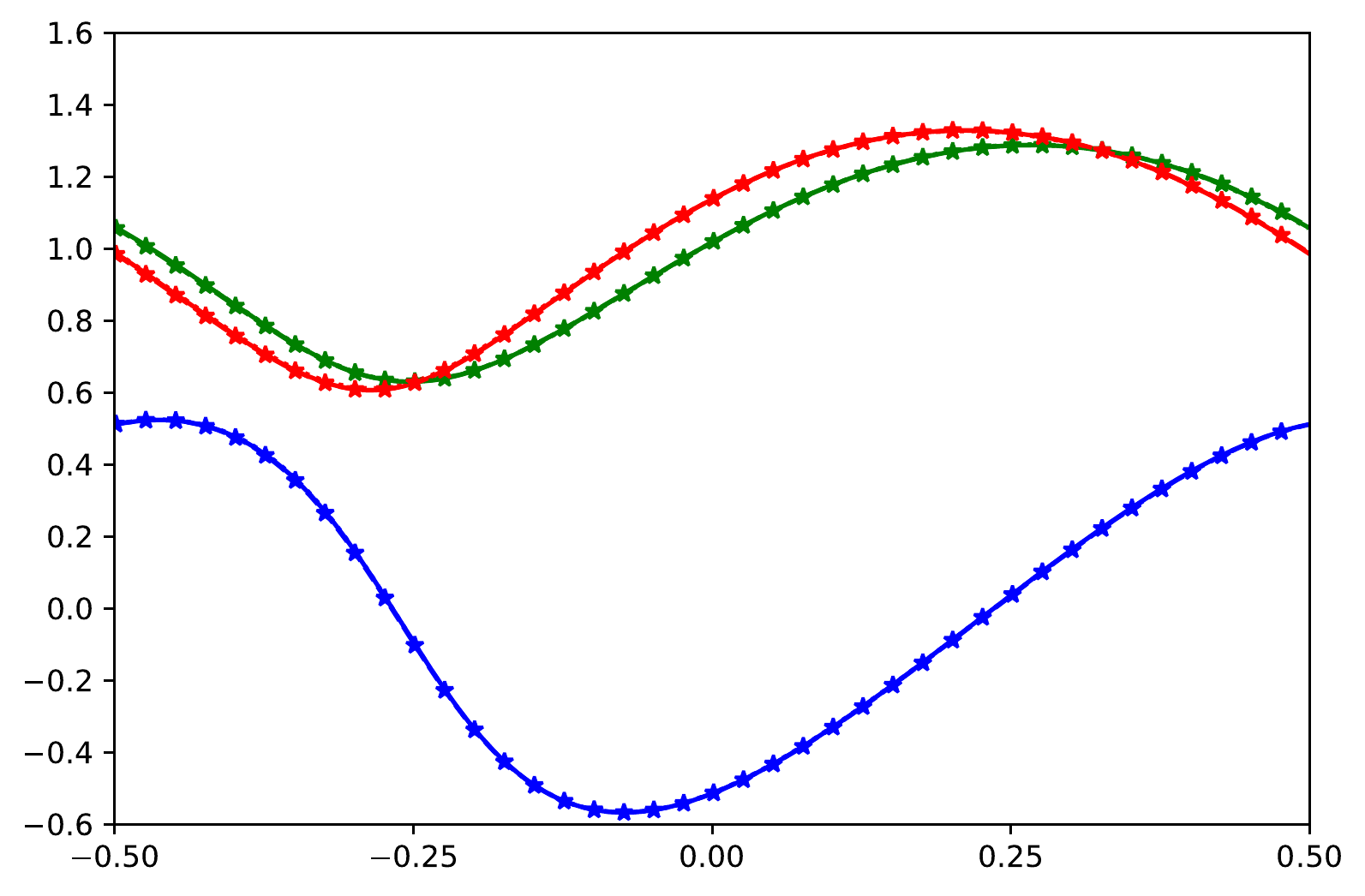} 
  } \hfill
  \subfloat[$\Kn = 0.1, t = 0.1$]{
      \includegraphics[width=0.25\textwidth]{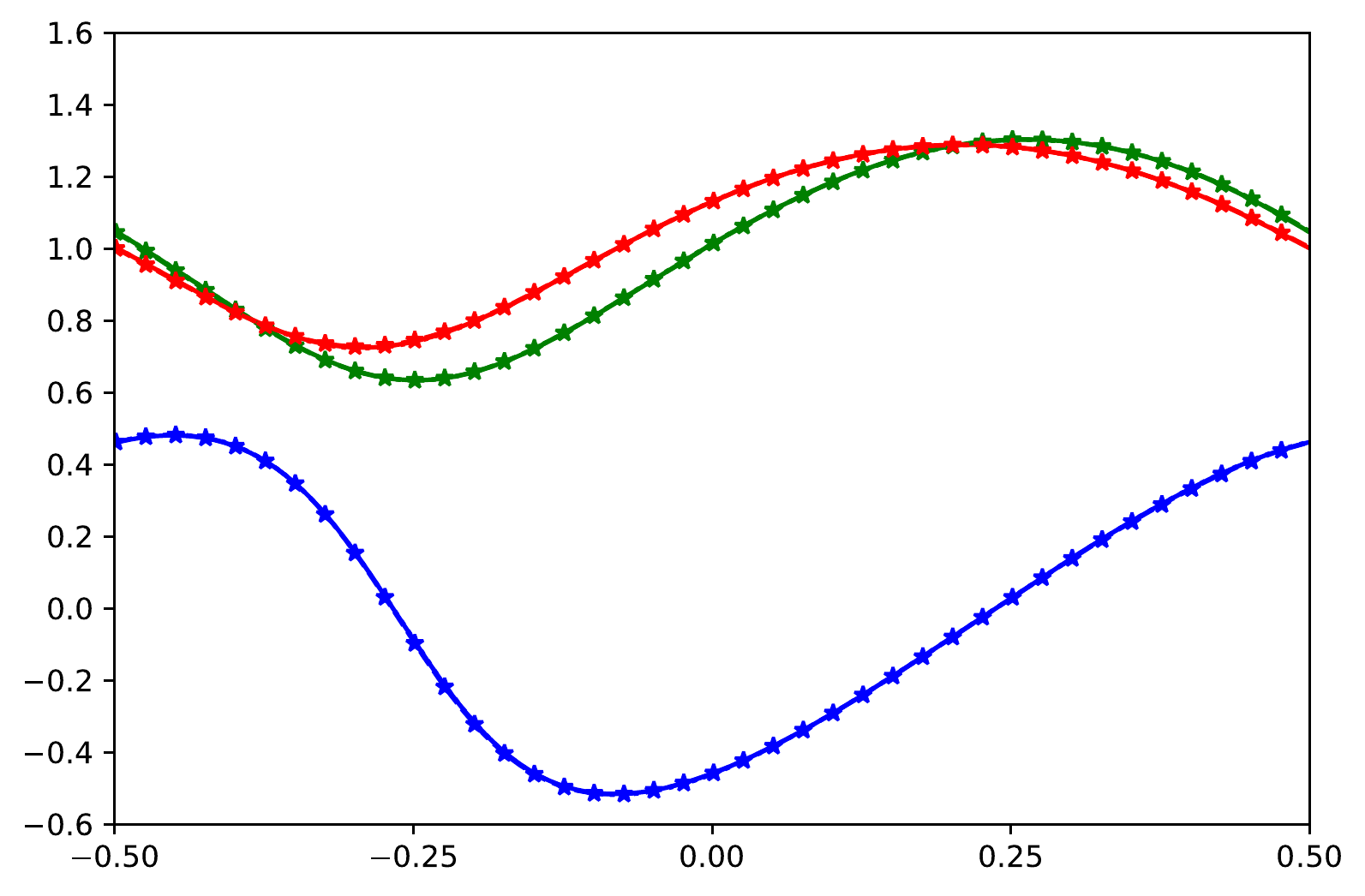} 
  }\hfill
  \subfloat[$\Kn = 1.0, t = 0.1$]{
      \includegraphics[width=0.25\textwidth]{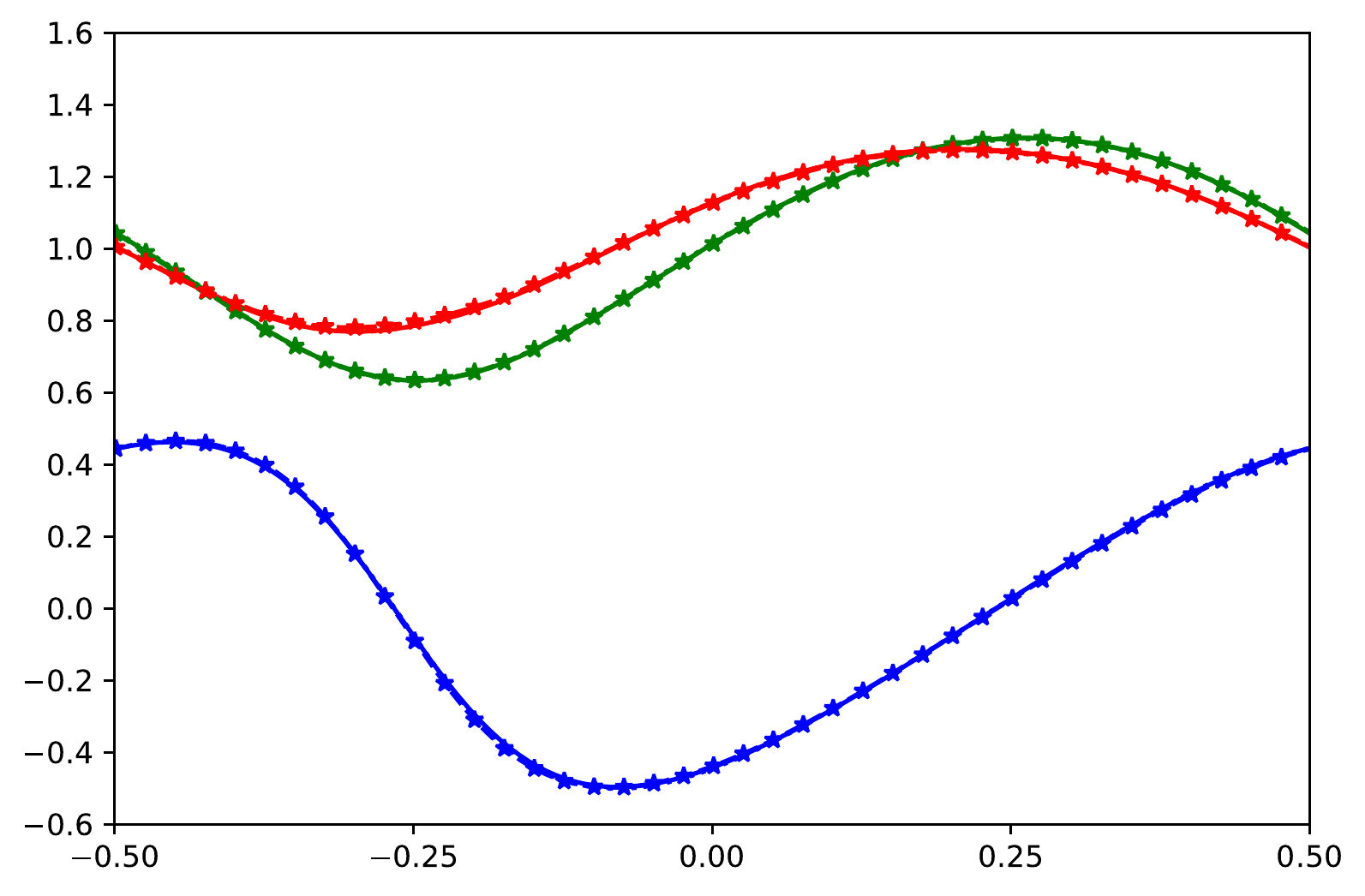} 
  }
  \caption{(1D wave problem in Sec. \ref{subsec:wave1d}) Numerical solution of 1D3V wave problem with BGK model by $\nnBGK$ and $\nnLR$. The first row corresponds to $t=0.0$, and  the second row corresponds to $t=0.1$. The dashed line is the numerical solution of ${\rm NR^{\rm BGK}}$, the solid line is that of ${\rm NSR}_{\rm LR}^{\rm BGK}$, and the dot-dash line is the reference solution by DVM. 
}
  \label{fig:wave-bgk}
\end{figure}
The numerical results of BGK model using $\nnBGK$ and $\rm NSR_{\rm LR}^{\rm BGK}$ are plotted in Fig. \ref{fig:wave-bgk}, where the density $\rho$, macroscopic velocity $u_1$, and temperature $T$ are studied. Since the initial condition should also be approximated by the neural network, there can be a small deviation between the numerical solution and the reference solution. The numerical solution at $t = 0$ and $t = 0.1$ are listed in Fig. \ref{fig:wave-bgk}. It illustrates that the numerical solution matches well with the reference solution at both time for the three Knudsen numbers. Here, the reference solution is obtained by DVM method, where the grid number in spatial space $N_x = 400$ with linear reconstruction  and the upwind numerical flux utilized. The computational domain in the microscopic velocity space is $[-10, 10]^3$ with $24$ grids in each direction.  

\begin{figure}[!hptb] 
  \centering 
  \subfloat[$\Kn=0.01, t = 0$]{
      \includegraphics[width=0.25\textwidth]{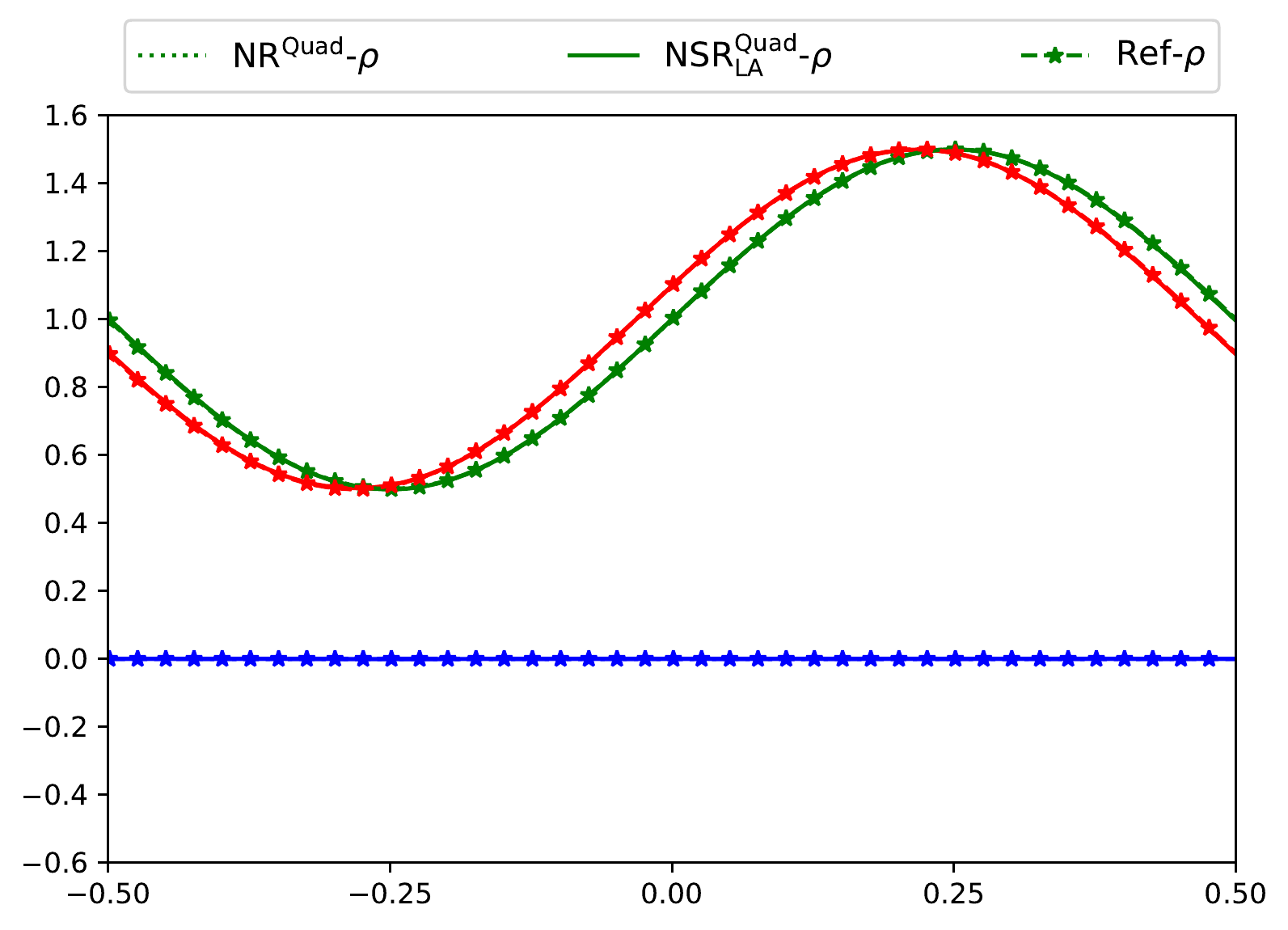} 
  }\hfill
  \subfloat[$\Kn=0.1, t = 0$]{
      \includegraphics[width=0.25\textwidth]{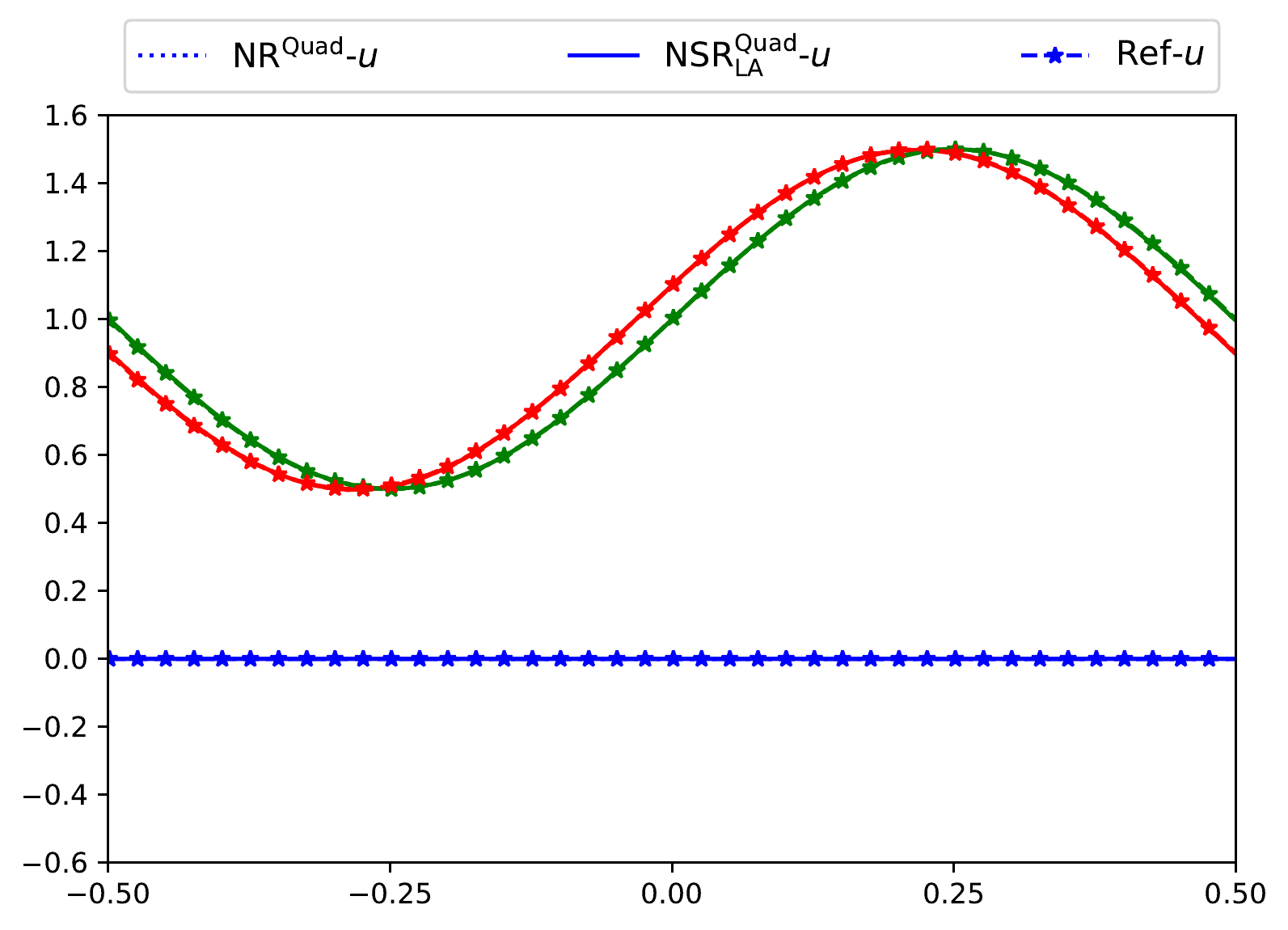} 
  }\hfill
  \subfloat[$\Kn=1.0, t = 0$]{
      \includegraphics[width=0.25\textwidth]{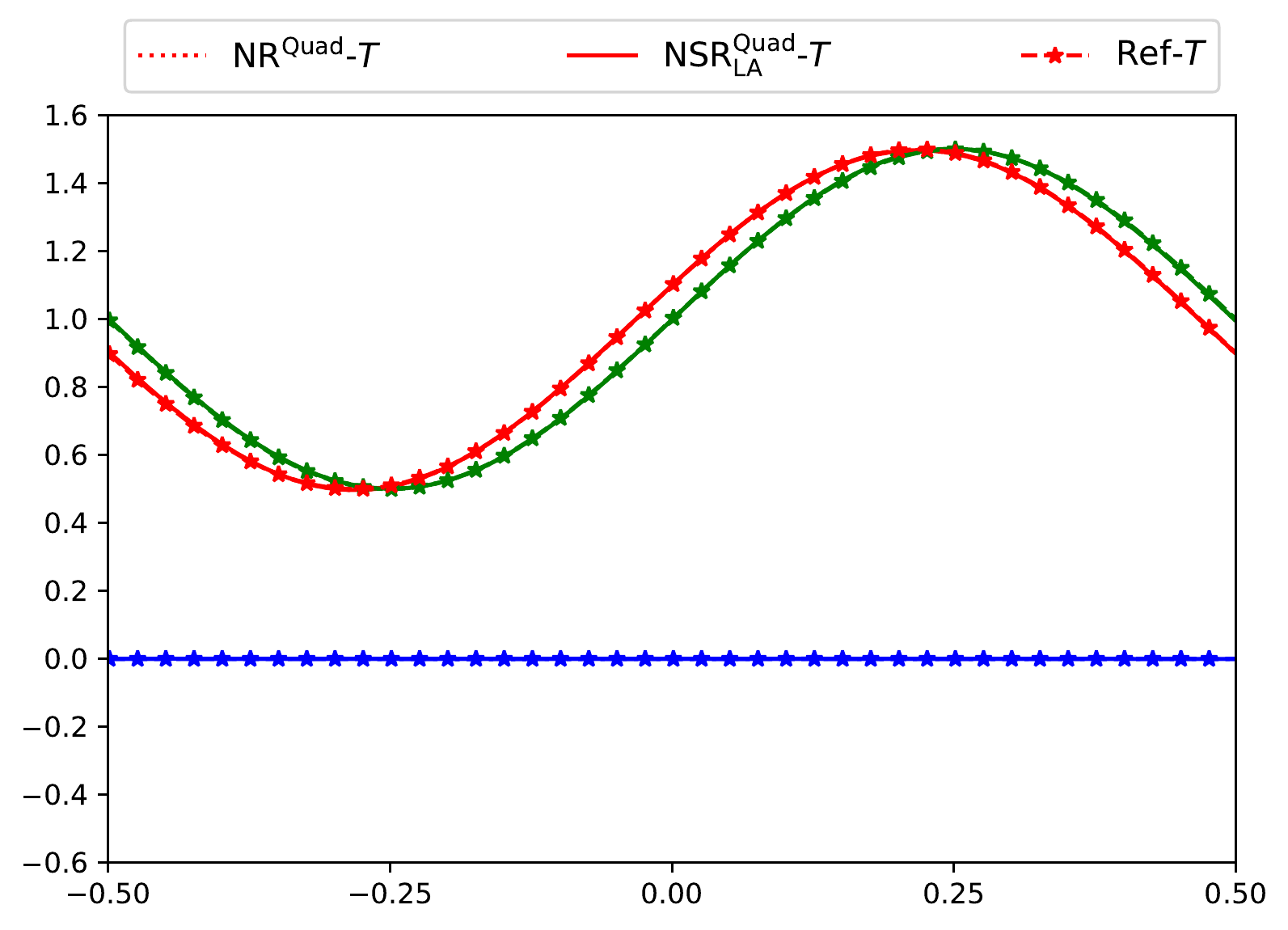} 
  } \\
   \subfloat[$\Kn=0.01, t = 0.1$]{
      \includegraphics[width=0.25\textwidth]{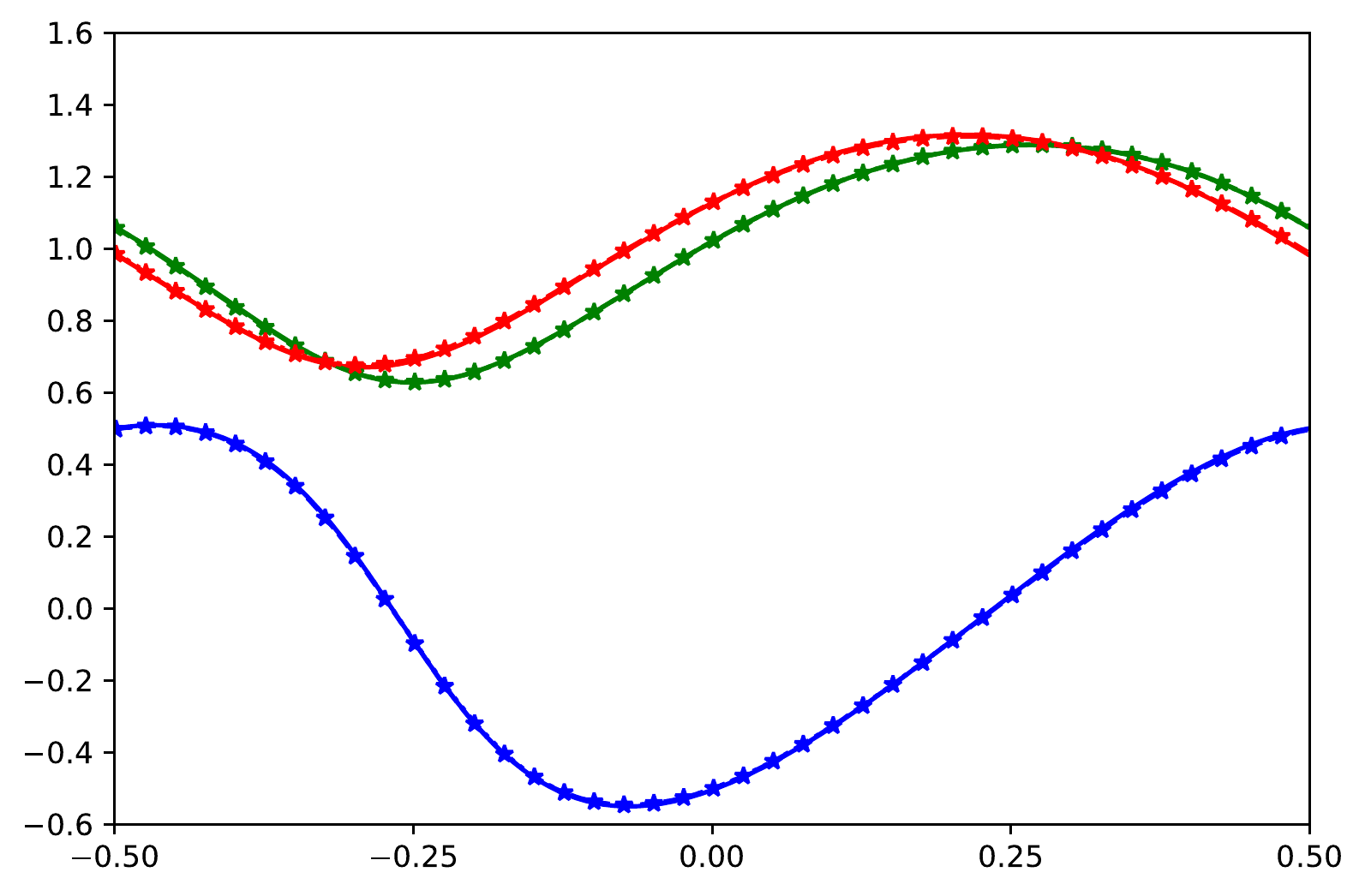} 
  }\hfill
  \subfloat[$\Kn=0.1, t = 0.1$]{
      \includegraphics[width=0.25\textwidth]{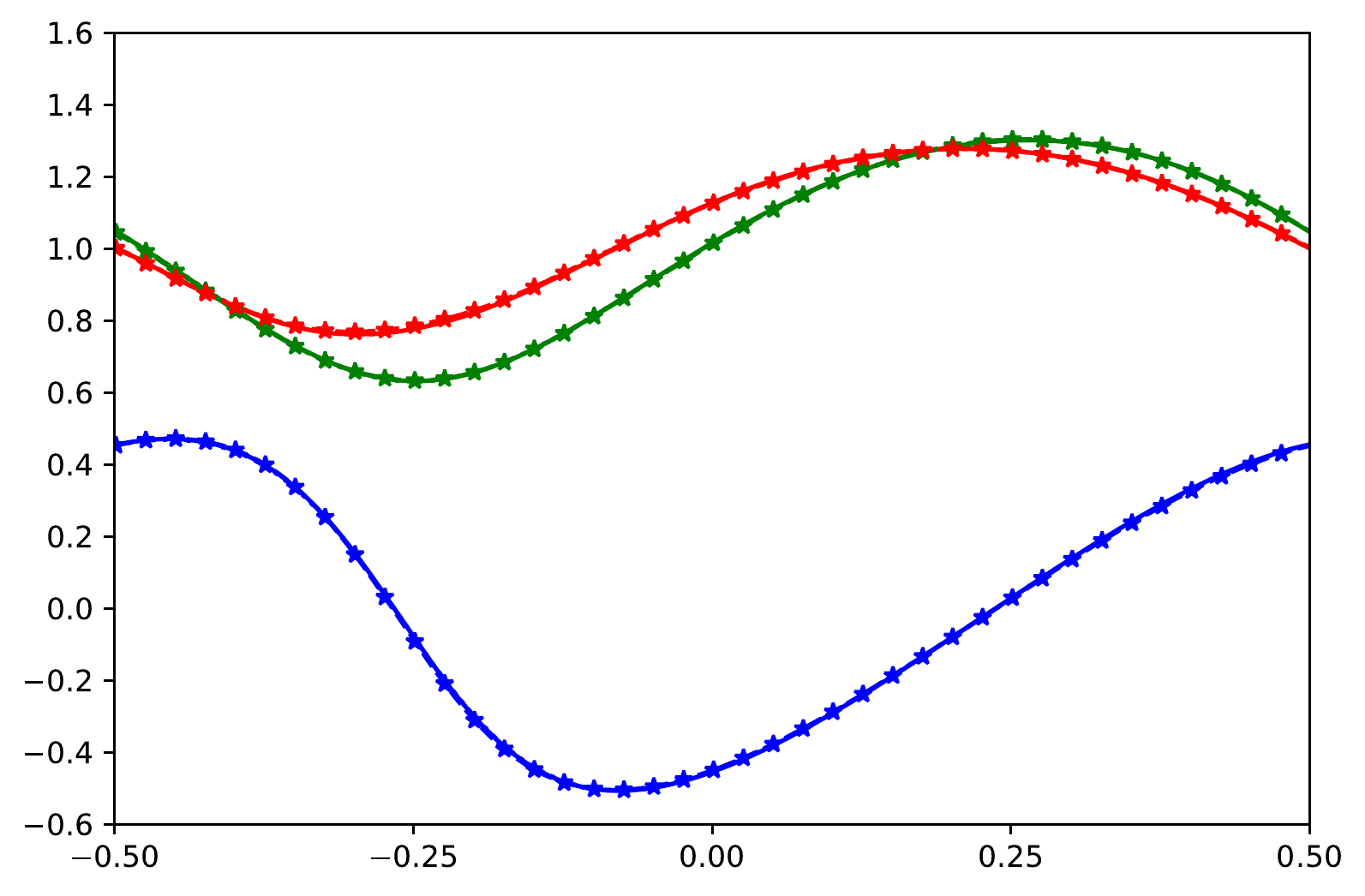} 
  }\hfill
  \subfloat[$\Kn=1.0, t = 0.1$]{
      \includegraphics[width=0.3\textwidth]{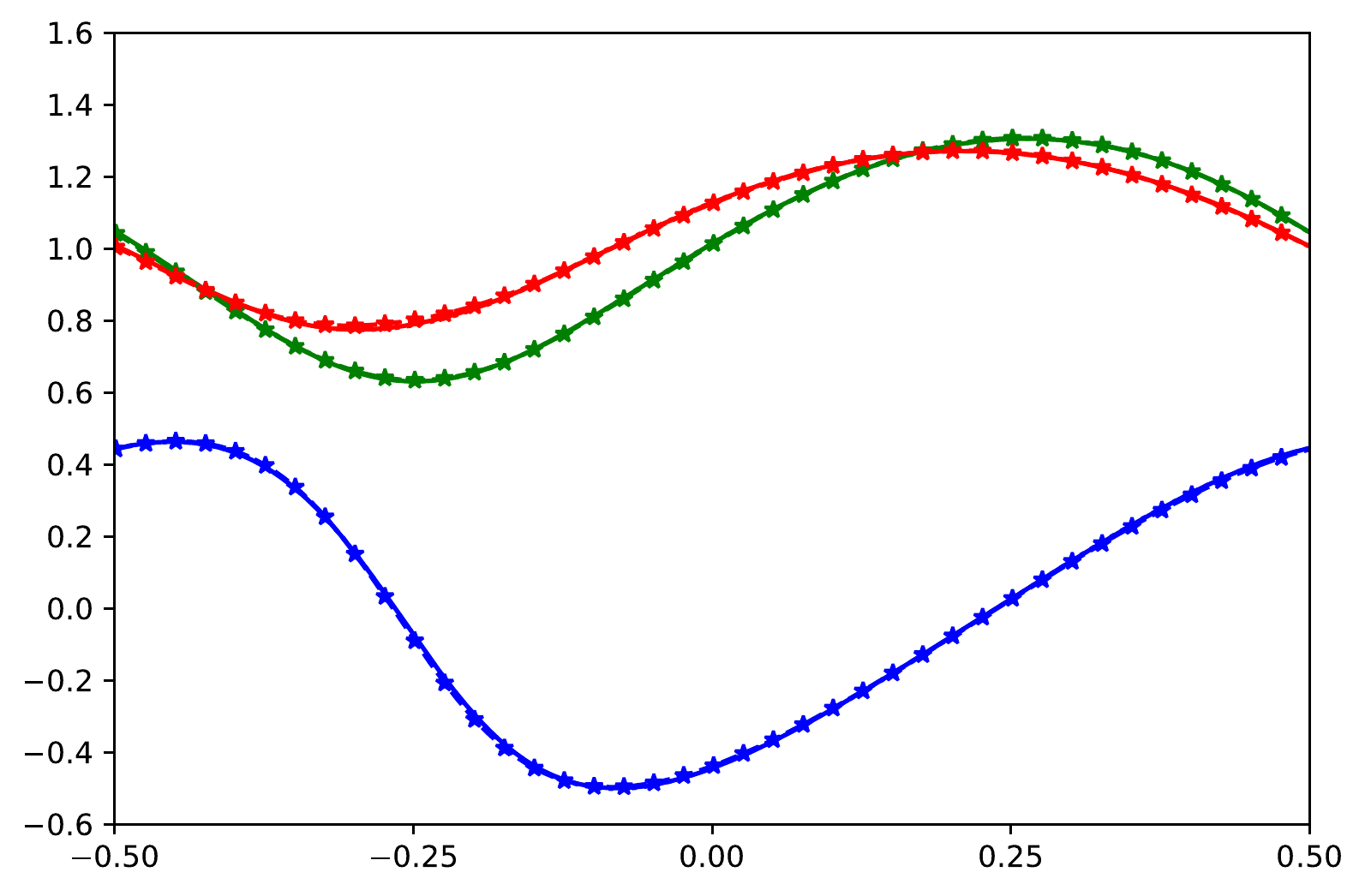} 
  }
   \caption{(1D wave problem in Sec. \ref{subsec:wave1d}) Numerical solution of 1D3V wave problem with quadratic collision model by $\nnFSM$ and $\nnLA$. The first row corresponds to $t=0.0$, and  the second row corresponds to $t=0.1$. The dashed line is the numerical solution of ${\rm NR^{\rm Quad}}$, the solid line is that of ${\rm NSR}_{\rm LA}^{\rm Quad}$, and the dot-dash line is the reference solution by fast Fourier method. }
  \label{fig:wave-fsm}
\end{figure}

The numerical results of the quadratic collision model using $\nnFSM$ and $\rm NSR_{\rm LA}^{\rm Quad}$ is provided in Fig. \ref{fig:wave-fsm}, where the same macroscopic variables $\rho, u_1$ and $T$ at $t = 0$ and $0.1$ with $\Kn = 0.01, 0.1$ and $1.0$ are shown. It is clear that for the quadratic collision model, the numerical solution also agrees well with the reference solution, where the reference solution is obtained by fast Fourier method with $24$ modes in each microscopic velocity direction.

Define the relative error between the numerical solution and the reference solution as 
\begin{equation}
    \label{eq:error1}
{\rm error} = \frac{\Vert s_{\rm num}- s_{\rm ref} \Vert_2}{\Vert s_{\rm num}\Vert_2}, \qquad s = \rho, T,
\end{equation}
for the density and temperature and 
\begin{equation}
    \label{eq:error2}
{\rm error} = \frac{\Vert u_{\rm num}-u_{\rm ref} \Vert_2}{1 + \Vert u_{\rm num}\Vert_2}, 
\end{equation}
for the macroscopic velocity to avoid the case $u = 0$. The relative error of the four methods with different Knudsen numbers at $t = 0$, and $0.1$ are shown in Tab. \ref{tab:wave}. It illustrates that the error of all the neural representation methods for both BGK and quadratic collision model reach the magnitude $\mO(10^{-3})$. In particular, the error of the sparse representation-based methods $\nnLR$ and $\nnLA$ is also at the same level compared to the direct neural representation methods $\nnBGK$ and $\nnFSM$, which means that $\nnLR$ and $\nnLA$ can reach a similar accuracy of $\nnBGK$ and $\nnFSM$.

\begin{table}[hpt!]
\centering
\def\arraystretch{1.5}
\scalebox{0.85}{
{\footnotesize
\begin{tabular}{@{}lclllllllll@{}}
\toprule
\multicolumn{1}{l}{Kn} & \multicolumn{1}{l}{} & \multicolumn{3}{c}{0.01}                            & \multicolumn{3}{c}{0.1}                             & \multicolumn{3}{c}{1.0}        \\ \midrule
  &  $t$   & \multicolumn{1}{c}{$\rho$} & \multicolumn{1}{c}{$u_1$} & \multicolumn{1}{c|}{$T$}      & \multicolumn{1}{c}{$\rho$} & \multicolumn{1}{c}{$u_1$} & \multicolumn{1}{c|}{$T$}      & \multicolumn{1}{c}{$\rho$} & \multicolumn{1}{c}{$u_1$} & \multicolumn{1}{c}{$T$} \\
\multirow{2}{*}{$\nnBGK$} & 0.0                    & 1.64e-03 & 6.23e-04 & \multicolumn{1}{l|}{1.60e-03} & 1.80e-03 & 1.06e-04 & \multicolumn{1}{l|}{1.64e-03} & 2.08e-03 & 4.09e-04 & 1.90e-03 \\
                     & 0.1                    & 1.63e-03 & 1.98e-03 & \multicolumn{1}{l|}{2.97e-03} & 1.08e-03 & 2.17e-03 & \multicolumn{1}{l|}{1.89e-03} & 2.30e-03 & 4.62e-03 & 4.14e-03 \\
\multirow{2}{*}{$\nnLR$}  & 0.0                    & 1.63e-03 & 3.67e-04 & \multicolumn{1}{l|}{1.67e-03} & 1.66e-03 & 2.78e-04 & \multicolumn{1}{l|}{2.26e-03} & 1.66e-03 & 8.44e-04 & 4.92e-03 \\
                     & 0.1                    & 1.13e-03 & 1.57e-03 & \multicolumn{1}{l|}{1.24e-03} & 1.13e-03 & 1.79e-03 & \multicolumn{1}{l|}{1.35e-03} & 1.26e-03 & 4.98e-03 & 4.85e-03 \\
\multirow{2}{*}{$\nnFSM$} & 0.0                    & 1.61e-03 & 2.08e-04 & \multicolumn{1}{l|}{1.75e-03} & 1.82e-03 & 2.20e-04 & \multicolumn{1}{l|}{2.07e-03} & 2.01e-03 & 3.72e-04 & 2.04e-03 \\
                     & 0.1                    & 1.07e-03 & 1.73e-03 & \multicolumn{1}{l|}{1.46e-03} & 1.58e-03 & 3.84e-03 & \multicolumn{1}{l|}{2.59e-03} & 2.30e-03 & 4.74e-03 & 4.71e-03 \\
\multirow{2}{*}{$\nnLA$}  & 0.0                    & 1.60e-03 & 1.84e-04 & \multicolumn{1}{l|}{1.84e-03} & 1.86e-03 & 1.48e-04 & \multicolumn{1}{l|}{1.80e-03} & 1.98e-03 & 2.02e-04 & 2.06e-03 \\
                     & 0.1                    & 1.33e-03 & 2.99e-03 & \multicolumn{1}{l|}{3.48e-03} & 1.57e-03 & 3.79e-03 & \multicolumn{1}{l|}{2.99e-03} & 2.23e-03 & 5.11e-03 & 4.10e-03 \\ \bottomrule
\end{tabular}}}
\caption{(1D wave problem in Sec. \ref{subsec:wave1d})  The relative error between the numerical solution by NR/NSR and the reference solution for the density $\rho$, macroscopic velocity $u_1$ and the temperature $T$ with $\Kn = 0.01, 0,1$ and $1$ at $t = 0$ and $0.1$.}
\label{tab:wave}
\end{table}

\subsection{1D Sod tube problem}
\label{sec:sod}

\begin{figure}[!hptb]
  \centering 
  \subfloat[$\Kn=0.01, t = 0$]{
      \includegraphics[width=0.25\textwidth]{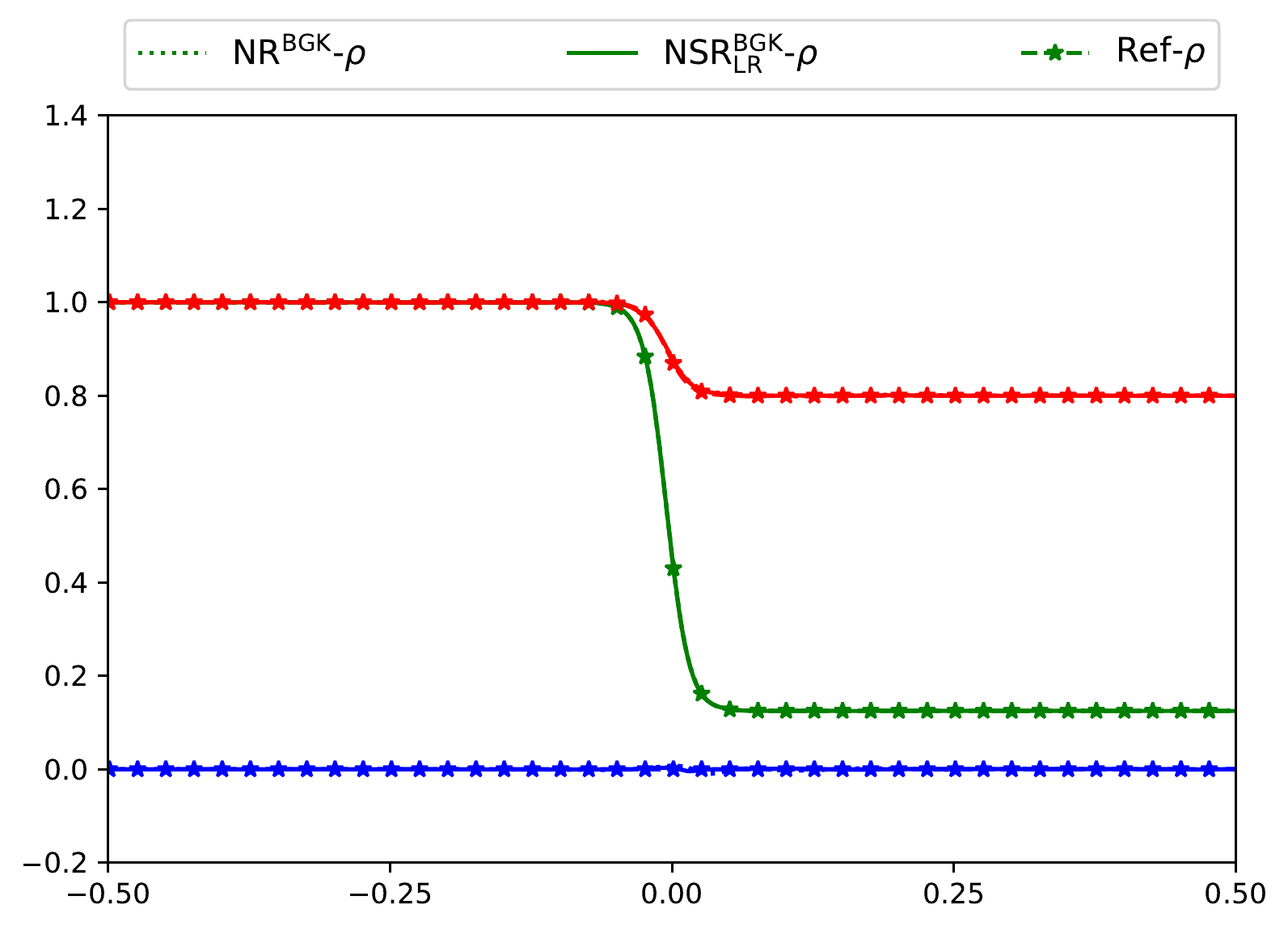} 
  } \hfill
  \subfloat[$\Kn=0.1, t = 0$]{
      \includegraphics[width=0.25\textwidth]{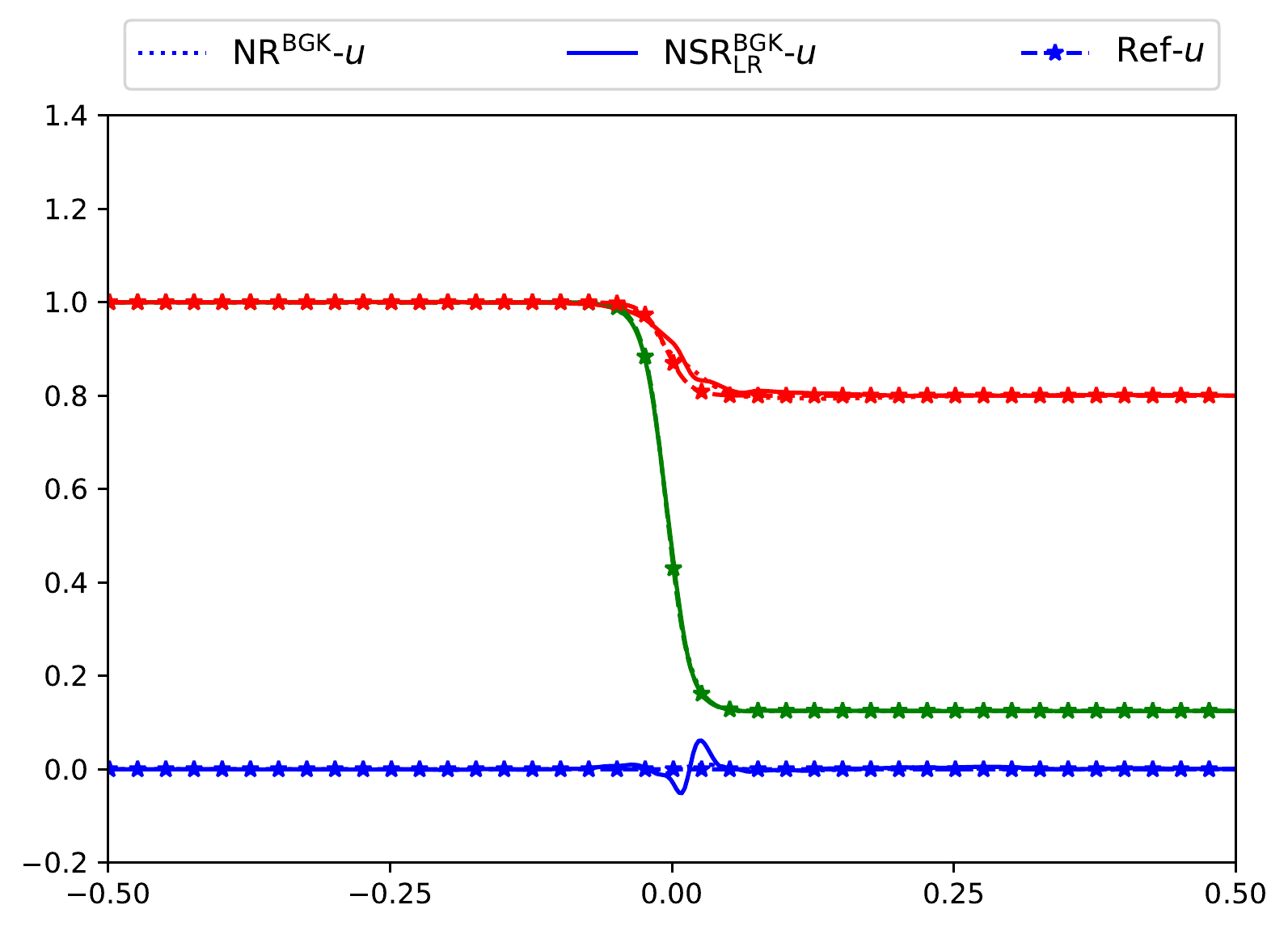} 
  } \hfill
  \subfloat[$\Kn=1.0, t = 0$]{
      \includegraphics[width=0.25\textwidth]{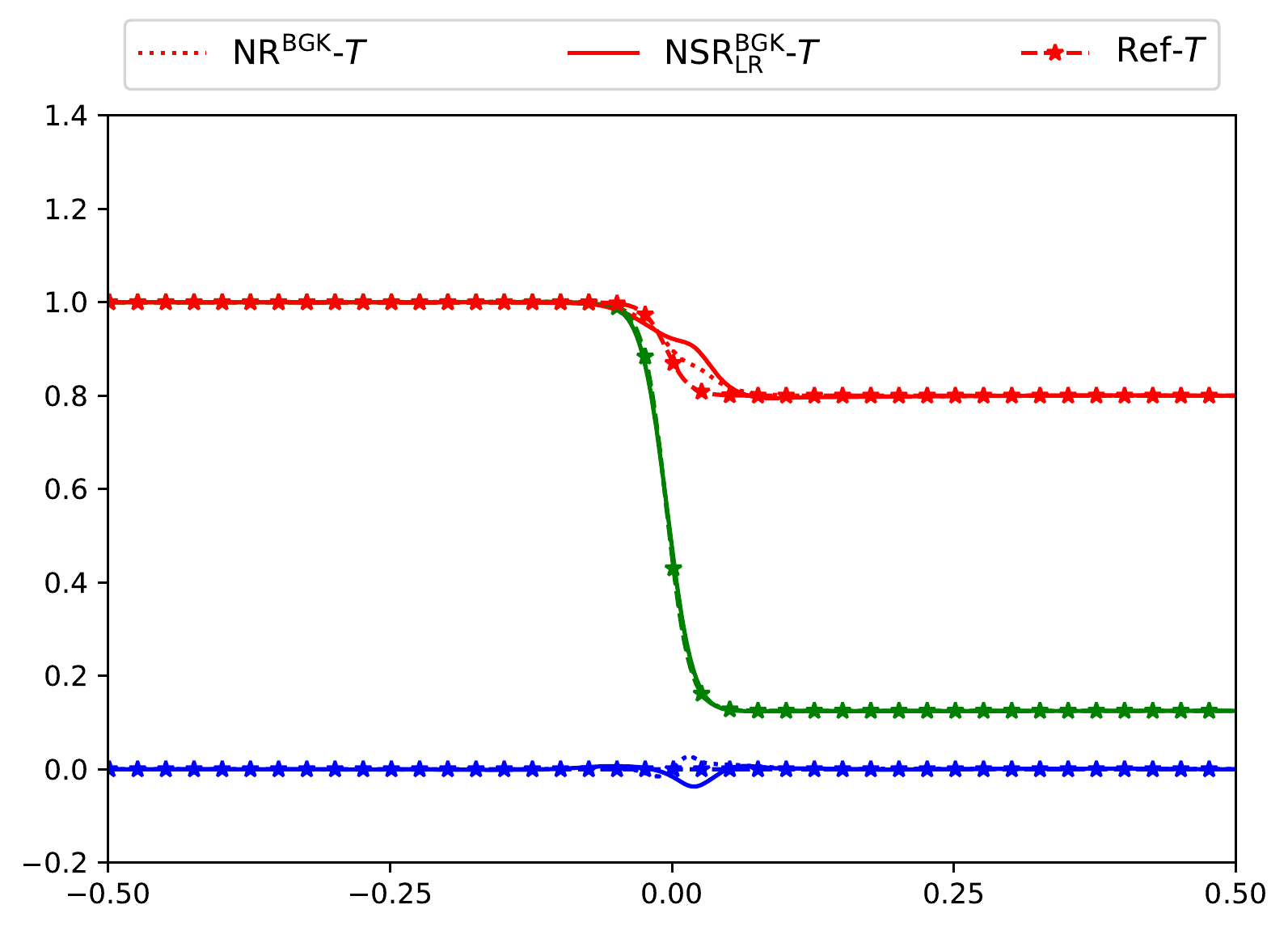} 
  } \\
  \subfloat[$\Kn=0.01, t = 0$]{
      \includegraphics[width=0.25\textwidth]{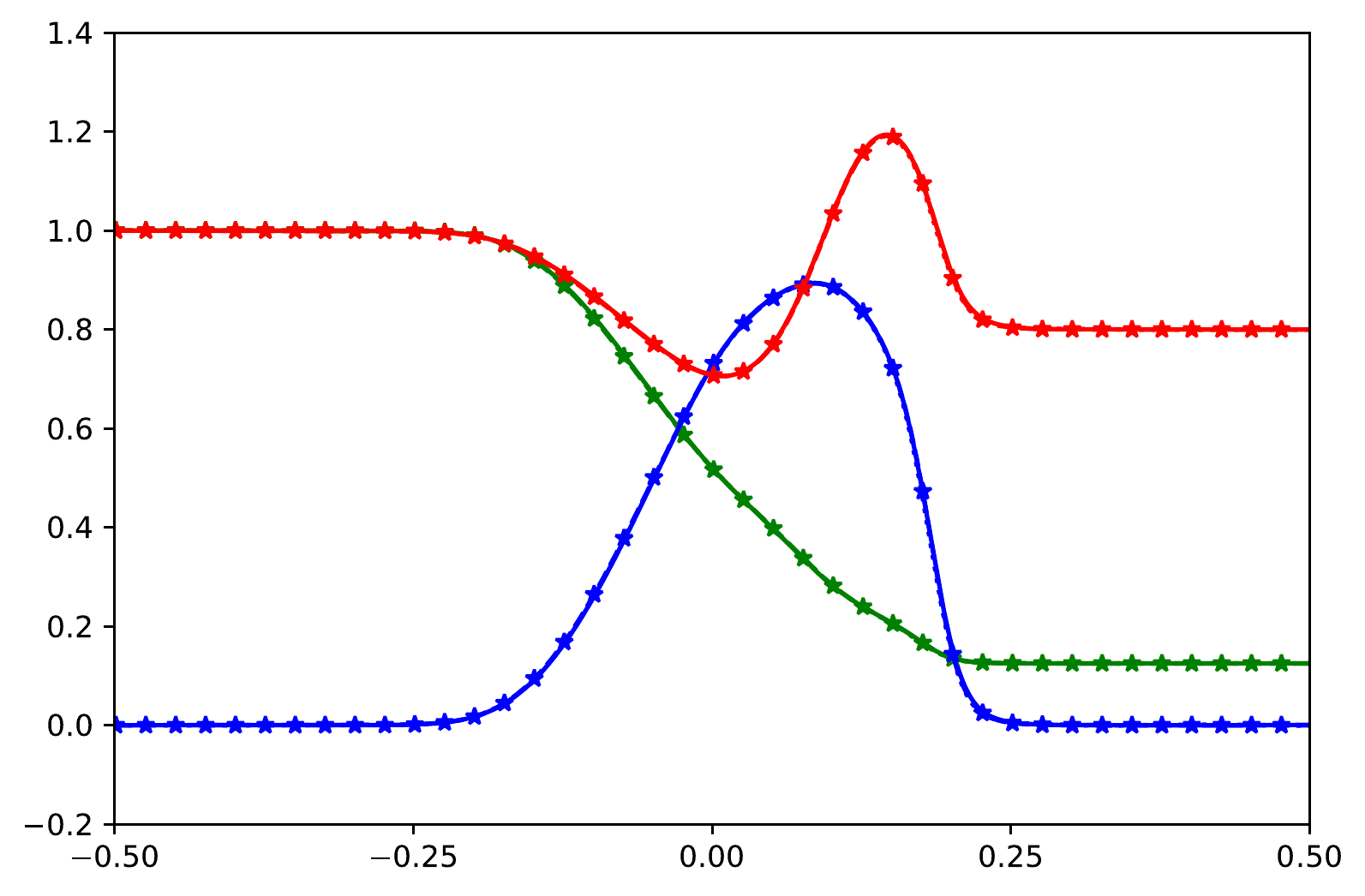} 
  } \hfill
  \subfloat[$\Kn=0.1, t = 0$]{
      \includegraphics[width=0.25\textwidth]{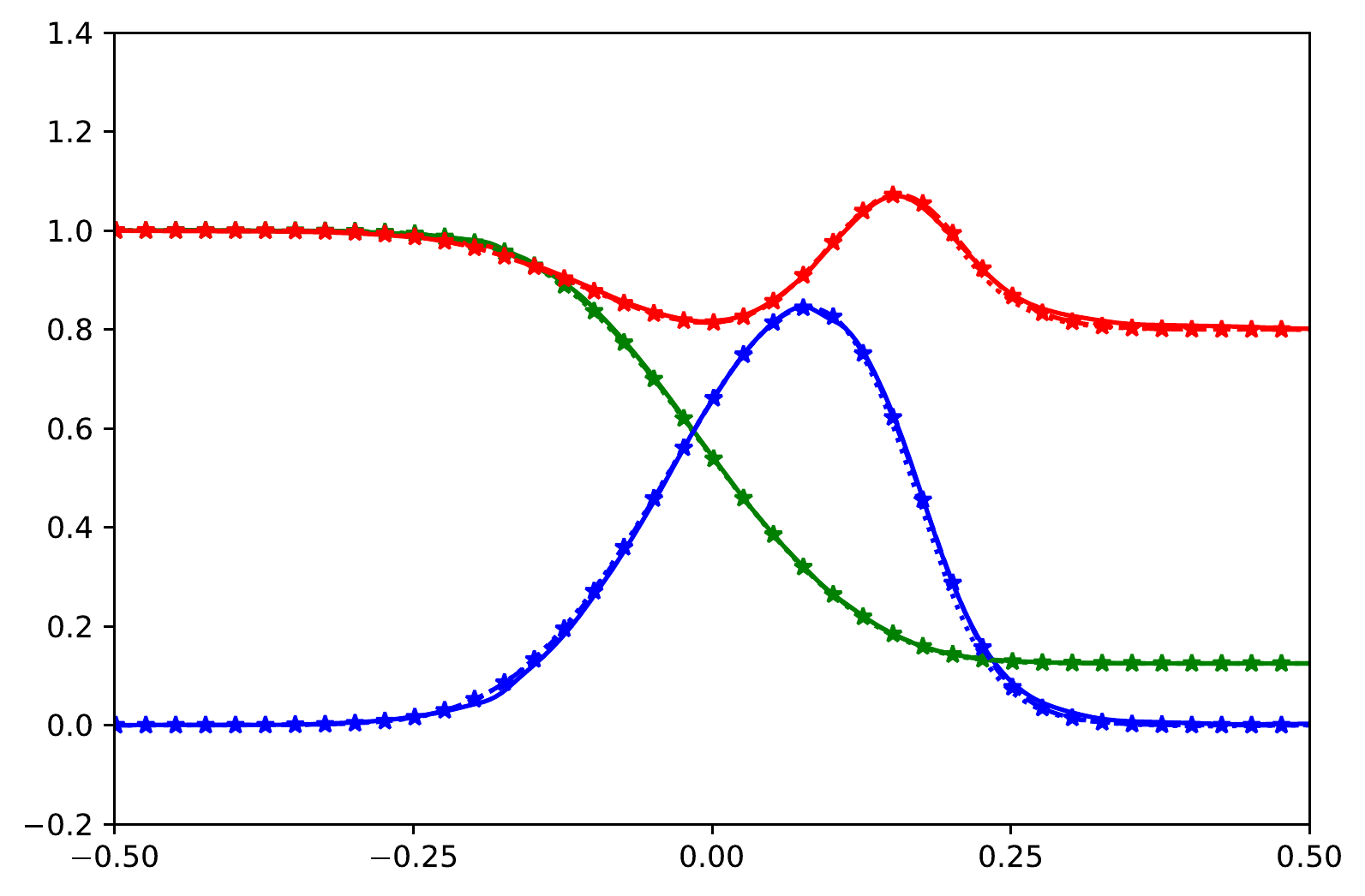} 
  } \hfill
  \subfloat[$\Kn=1.0, t = 0$]{
      \includegraphics[width=0.25\textwidth]{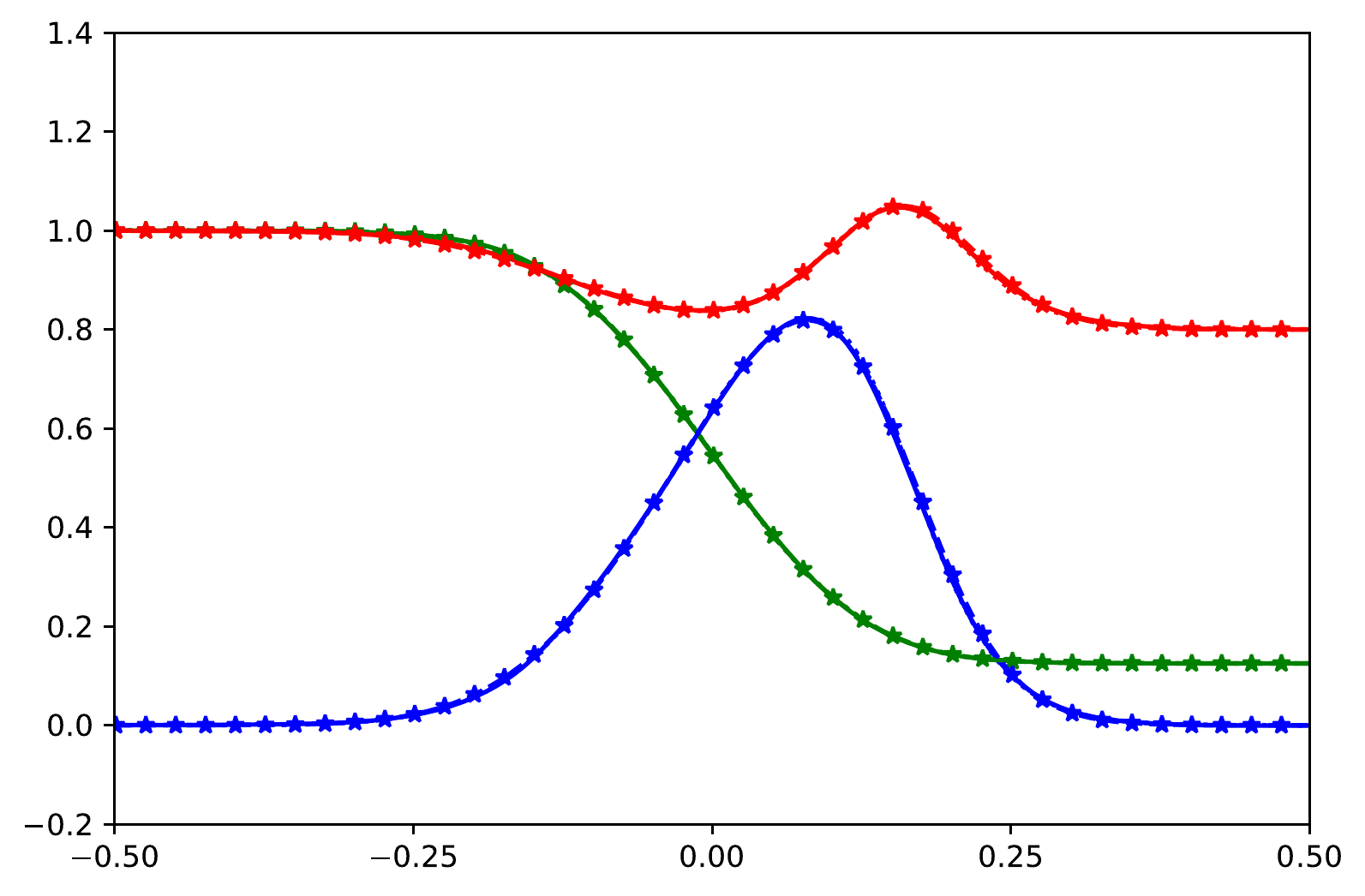} 
  } 
  \caption{(1D Sod tube problem in Sec. \ref{sec:sod}) Numerical solution of 1D3V Sod tube problem with BGK model by $\nnBGK$ and $\nnLR$. The first row corresponds to $t=0.0$, and  the second row corresponds to $t=0.1$. The dashed line is the numerical solution of ${\rm NR^{\rm BGK}}$, solid line is that of ${\rm NSR}_{\rm LR}^{\rm BGK}$, and the dot-dash line is the reference solution by DVM.}
\label{fig:sod-bgk}
\end{figure}

In this section, the classical 1D3V Sod tube problem is studied, where the initial condition is Maxwellian, with the macroscopic variables as below
\begin{equation}
\label{eq:ex2_ini_1}
\left[\rho(x), u_1(x), T(x)\right] = \left\{
\begin{aligned}
&[\rho_l, u_l,  T_l],  &x < 0, \\
&[\rho_r, u_r, T_r],  &x \geqslant  0. 
\end{aligned}
\right. 
\end{equation}
with 
\begin{equation}
    \label{eq:ex2_ini_2}
 \rho_l = 1.0, \qquad \rho_r = 0.125, \qquad u_l = u_r = 0, \qquad  T_l = 1.0, \qquad  T_r = 0.8.
\end{equation}
This is a problem with discontinuous initial condition, and a similar one is also studied in \cite{li2022learning,patel2022thermodynamically}. Since the neural network can not represent the discontinuous functions as well as the smooth functions, which is a common problem for network-based methods when solving PDEs, and there did not exist a generally effective method to solve this yet \cite{fuks2020limitationsa,lv2021hybrida}, the smoothing technique is utilized for this discontinuous initial condition problem. Precisely, the smoothing function is 
\begin{equation}
\label{eq:ex2_smooth}
\begin{aligned}
\rho(x)=\frac{\rho_r-\rho_l}{1+e^{-x/b}}+\rho_l, \qquad 
u(x)=\frac{u_r-u_l}{1+e^{-x/b}}+u_l, \qquad 
T(x)=\frac{T_r-T_l}{1+e^{-x/b}}+T_l, 
\end{aligned}
\end{equation}
with the smoothing factor  $b=0.005$. 

\begin{figure}[!hptb]
  \centering 
   \subfloat[$\Kn=0.01, t = 0$]{
      \includegraphics[width=0.25\textwidth]{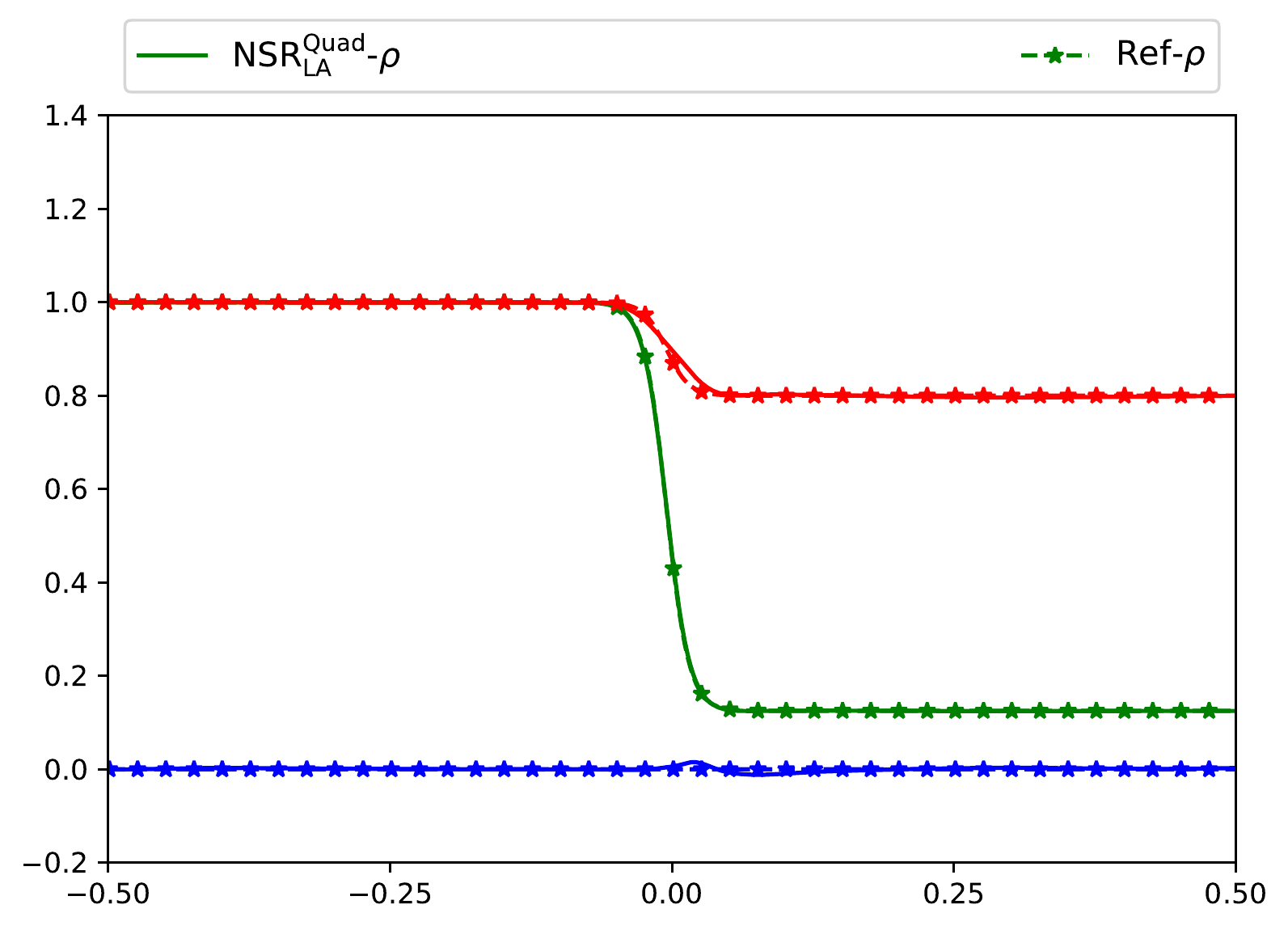} 
  }\hfill
  \subfloat[$\Kn=0.1, t = 0$]{
      \includegraphics[width=0.25\textwidth]{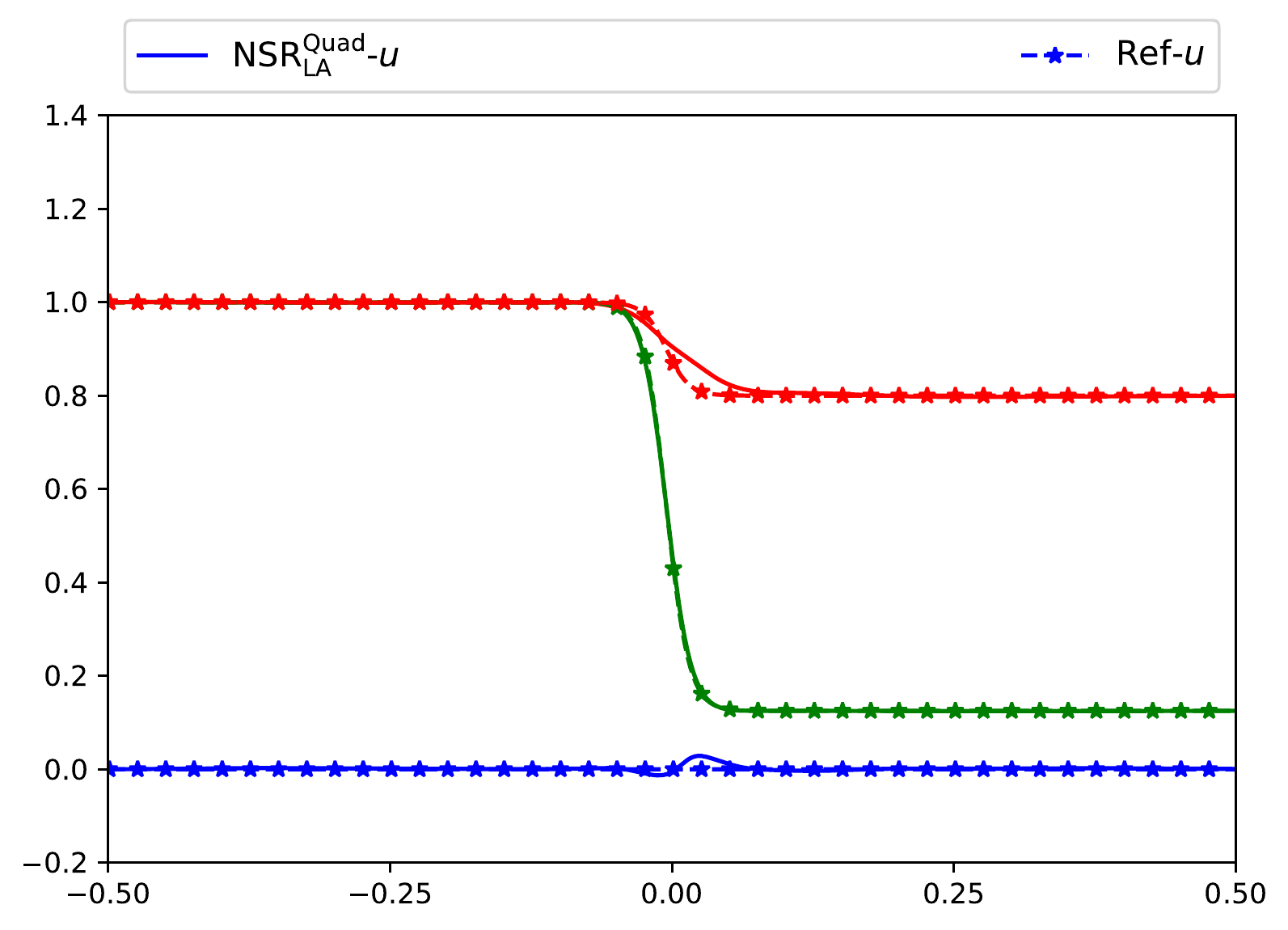} 
  }\hfill
  \subfloat[$\Kn=1.0, t = 0$]{
      \includegraphics[width=0.25\textwidth]{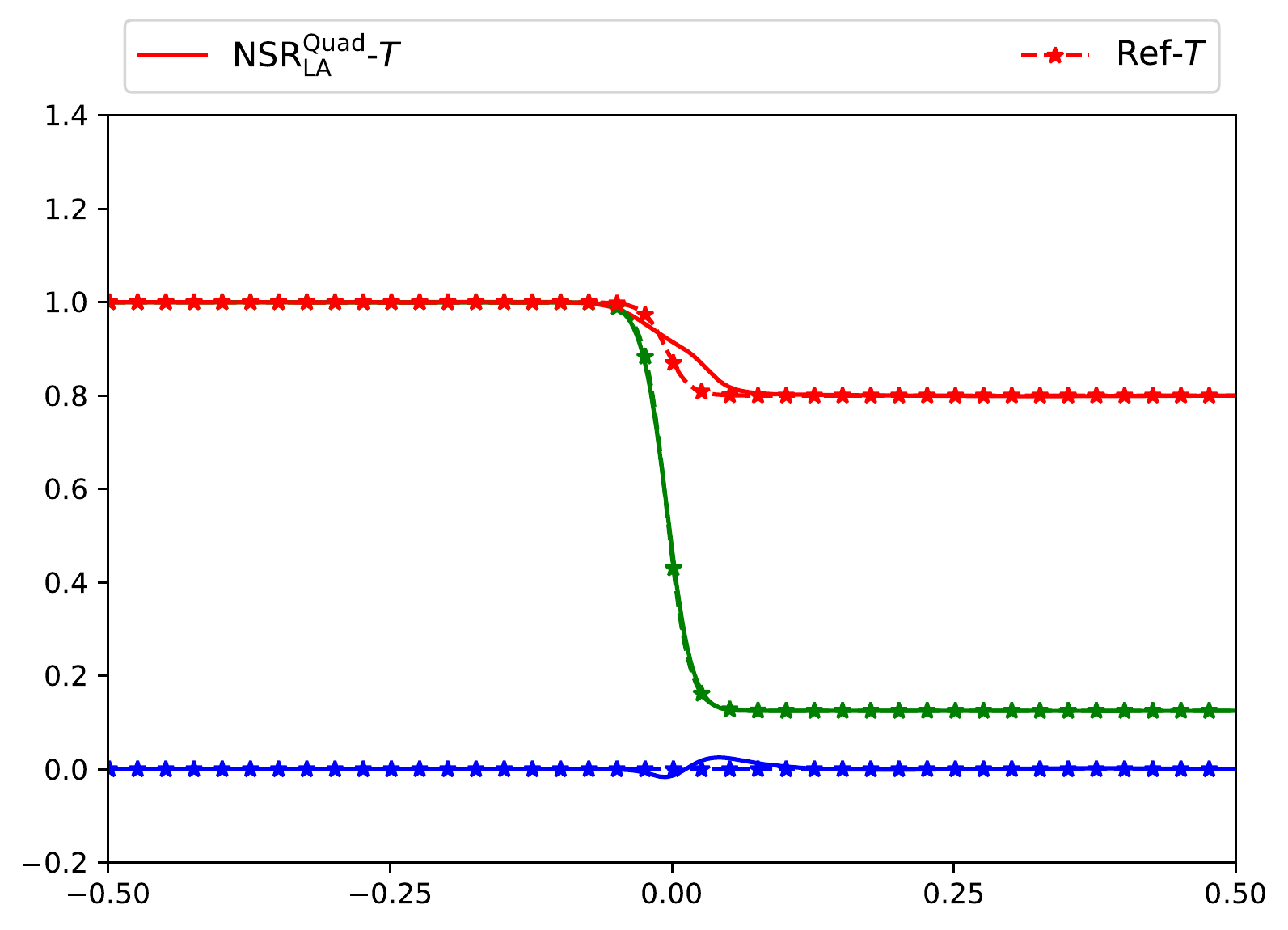} 
  } \\
  \subfloat[$\Kn=0.01, t = 0$]{
      \includegraphics[width=0.25\textwidth]{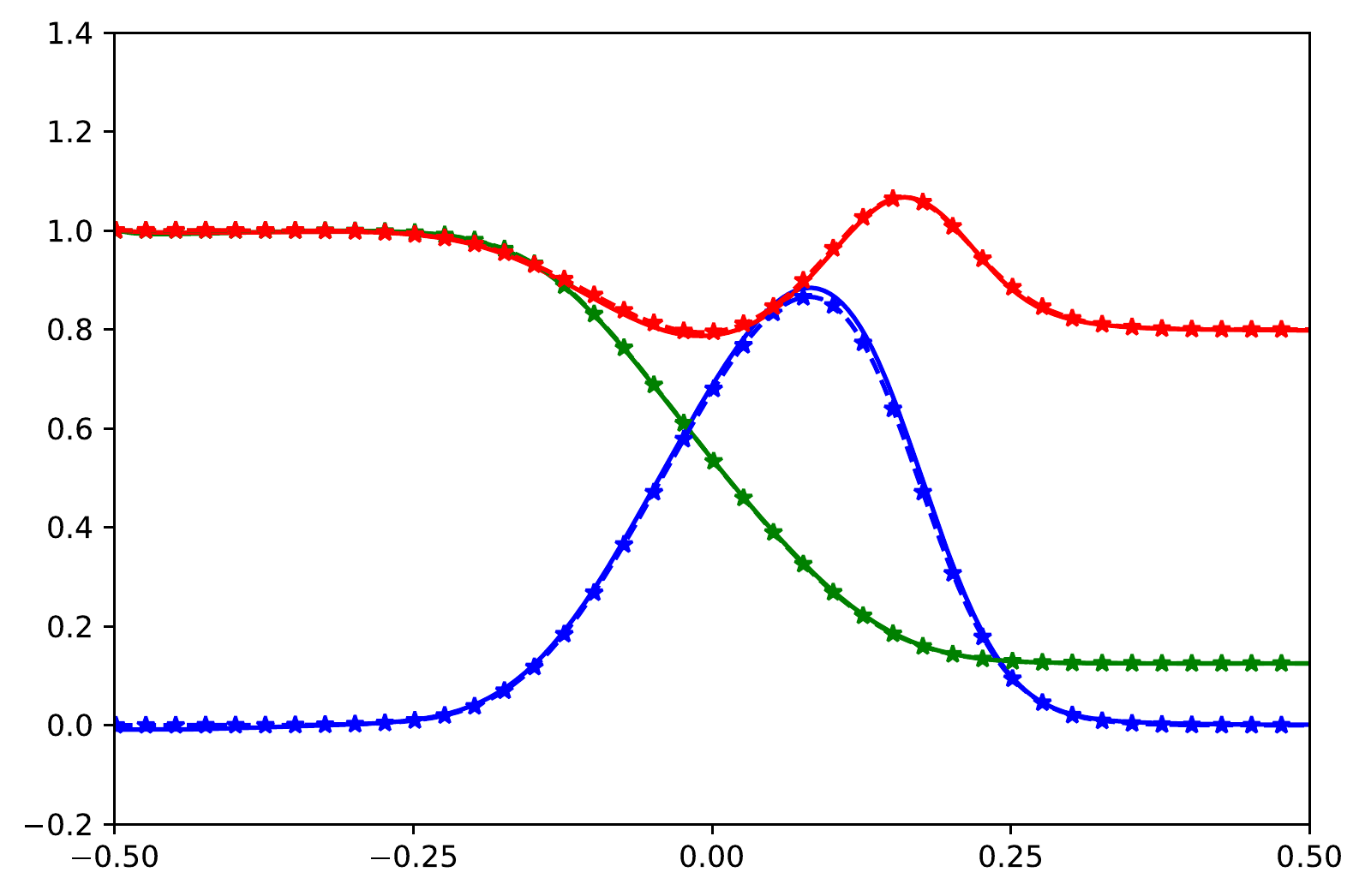} 
  }\hfill
  \subfloat[$\Kn=0.1, t = 0$]{
      \includegraphics[width=0.25\textwidth]{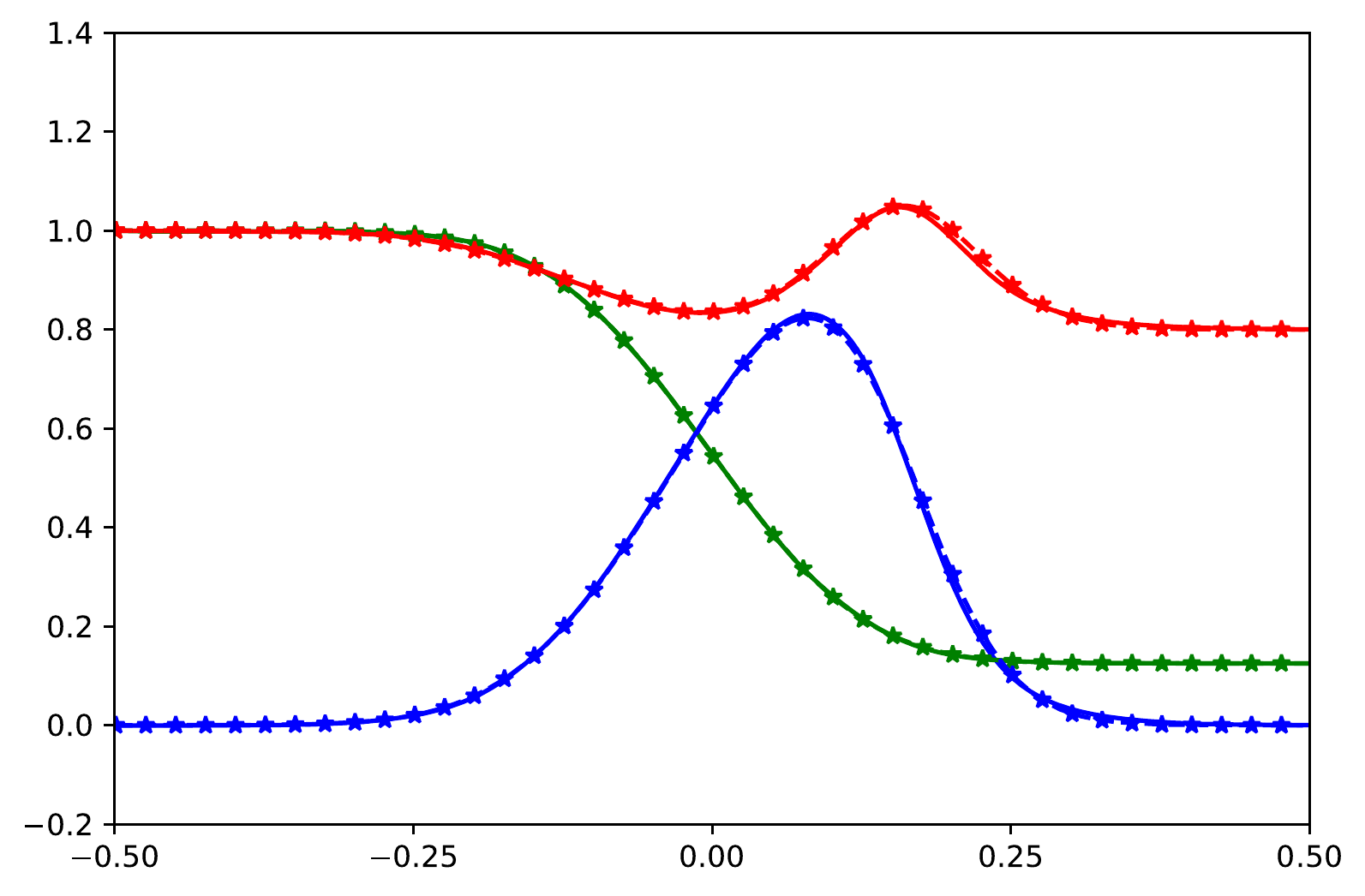} 
  } \hfill
  \subfloat[$\Kn=1.0, t = 0$]{
      \includegraphics[width=0.3\textwidth]{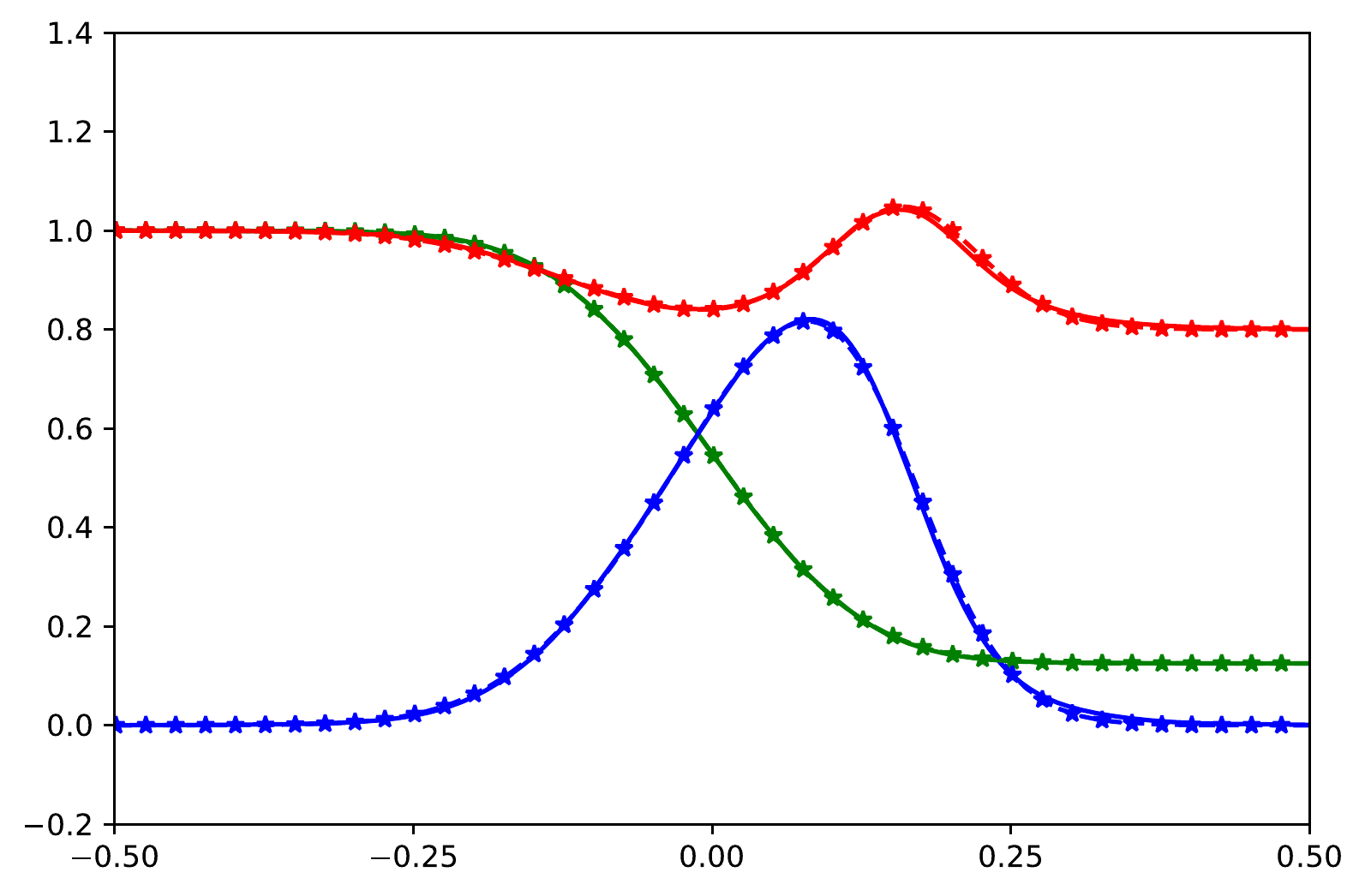} 
  } 
  \caption{(1D Sod tube problem in Sec. \ref{sec:sod}) Numerical solution of 1D3V Sod tube problem with quadratic collision model by $\nnLA$. The first row corresponds to $t=0.0$, and  the second row corresponds to $t=0.1$. The solid line is the numerical solution of ${\rm NSR}_{\rm LA}^{\rm Quad}$, and the dot-dash line is the reference solution by fast Fourier spectral method. 
  }
    \label{fig:sod-svd}
\end{figure}

The BGK and quadratic collision model are studied, where the parameters here are the same as in Tab. \ref{tab:ex1_para}. The numerical solution by  $\nnBGK$, $\nnLR$ is shown in Fig. \ref{fig:sod-bgk}, with the reference solution obtained  by DVM. Due to the discontinuous property here, the grid points in the microscopic velocity space of  DVM are increased to $[96, 24, 24]$. In Fig. \ref{fig:sod-bgk}, we can see that there is a little discrepancy for the initial condition, even though it  has been smoothed. Compared to the initial condition, the numerical solution and reference solution agree well with each other at $t = 0.1$. The numerical solution of the quadratic collision model by $\nnLA$ is plotted in Fig. \ref{fig:sod-svd}, while that of $\nnFSM$ is not presented due to memory limitation. The behavior of the quadratic collision model is similar to that of BGK model, where there is a relatively larger error in the initial condition while they are almost on top of each other at $t = 0.1$.

This example shows that even though the neural network can not approximate the discontinuous initial condition well,  the error does not rise monotonically with time increasing. This phenomenon also appeared in \cite{mattey2021physics, schafer2022generalization}. A possible reason for this is that the solution is gradually getting smoother, which causes the error to decrease as well.

Tab. \ref{tab:sod} shows the relative error defined in \eqref{eq:error1} and \eqref{eq:error2} between the numerical solution and the reference solution. As stated, the neural network can not approximate the discontinuous functions well, and the error at $t = 0$ for the temperature is increased to $\mO(10^{-2})$, though the smoothing technique is utilized. Moreover, the numerical solution is becoming smoother with time increasing, and the error is decreasing, which is different from the traditional numerical method.

\begin{table}[hpt!]
\centering
\def\arraystretch{1.5}
\scalebox{0.85}{
{\footnotesize
\begin{tabular}{@{}lllllllllll@{}}
\toprule
\multicolumn{1}{c}{Kn}           &    & \multicolumn{3}{c}{0.01}                                                             & \multicolumn{3}{c}{0.1}                                                              & \multicolumn{3}{c}{1.0}                                                        \\ \midrule
                                 &$t$ & \multicolumn{1}{c}{$\rho$} & \multicolumn{1}{c}{$u$} & \multicolumn{1}{c|}{$T$}      & \multicolumn{1}{c}{$\rho$} & \multicolumn{1}{c}{$u$} & \multicolumn{1}{c|}{$T$}      & \multicolumn{1}{c}{$\rho$} & \multicolumn{1}{c}{$u$} & \multicolumn{1}{c}{$T$} \\
\multirow{2}{*}{$\nnBGK$} & 0.0 & 4.55e-04                   & 1.87e-03                & \multicolumn{1}{l|}{1.79e-03} & 1.56e-03                   & 2.42e-03                & \multicolumn{1}{l|}{7.15e-03} & 2.90e-03                   & 4.56e-03                & 1.02e-02                \\
                                 & 0.1 & 1.12e-03                   & 2.31e-03                & \multicolumn{1}{l|}{1.93e-03} & 1.95e-03                   & 6.02e-03                & \multicolumn{1}{l|}{4.42e-03} & 1.15e-03                   & 3.53e-03                & 4.02e-03                \\
\multirow{2}{*}{$\nnLR$}  & 0.0 & 2.38e-04                   & 7.79e-04                & \multicolumn{1}{l|}{1.19e-03} & 3.05e-03                   & 8.97e-03                & \multicolumn{1}{l|}{8.76e-03} & 5.92e-03                   & 6.18e-03                & 1.75e-02                \\
                                 & 0.1 & 1.42e-03                   & 1.53e-03                & \multicolumn{1}{l|}{1.40e-03} & 3.07e-03                   & 4.83e-03                & \multicolumn{1}{l|}{5.47e-03} & 1.36e-03                   & 3.64e-03                & 3.09e-03                \\
\multirow{2}{*}{$\nnLA$} & 0.0 & 1.73e-03                   & 3.83e-03                & \multicolumn{1}{l|}{6.12e-03} & 3.61e-03                   & 5.30e-03                & \multicolumn{1}{l|}{1.22e-02} & 5.00e-03                   & 5.77e-03                & 1.41e-02                \\
                                 & 0.1 & 3.31e-03                   & 8.16e-03                & \multicolumn{1}{l|}{4.76e-03} & 1.23e-03                   & 5.03e-03                & \multicolumn{1}{l|}{5.92e-03} & 1.17e-03                   & 4.78e-03                & 5.23e-03                \\ \bottomrule
\end{tabular}}
}
\caption{(Sod tube problem in Sec. \ref{sec:sod})
 The relative error between the numerical solution by NR/NSR and the reference solution for the density $\rho$, macroscopic velocity $u_1$ and the temperature $T$ with $\Kn = 0.01, 0,1$ and $1$ at $t = 0$ and $0.1$.}
\label{tab:sod}
\end{table}

\subsection{Two-dimensional case}
\label{sec:2D}

\begin{figure}[!hptb]
\centering 
\includegraphics[width=0.25\linewidth]{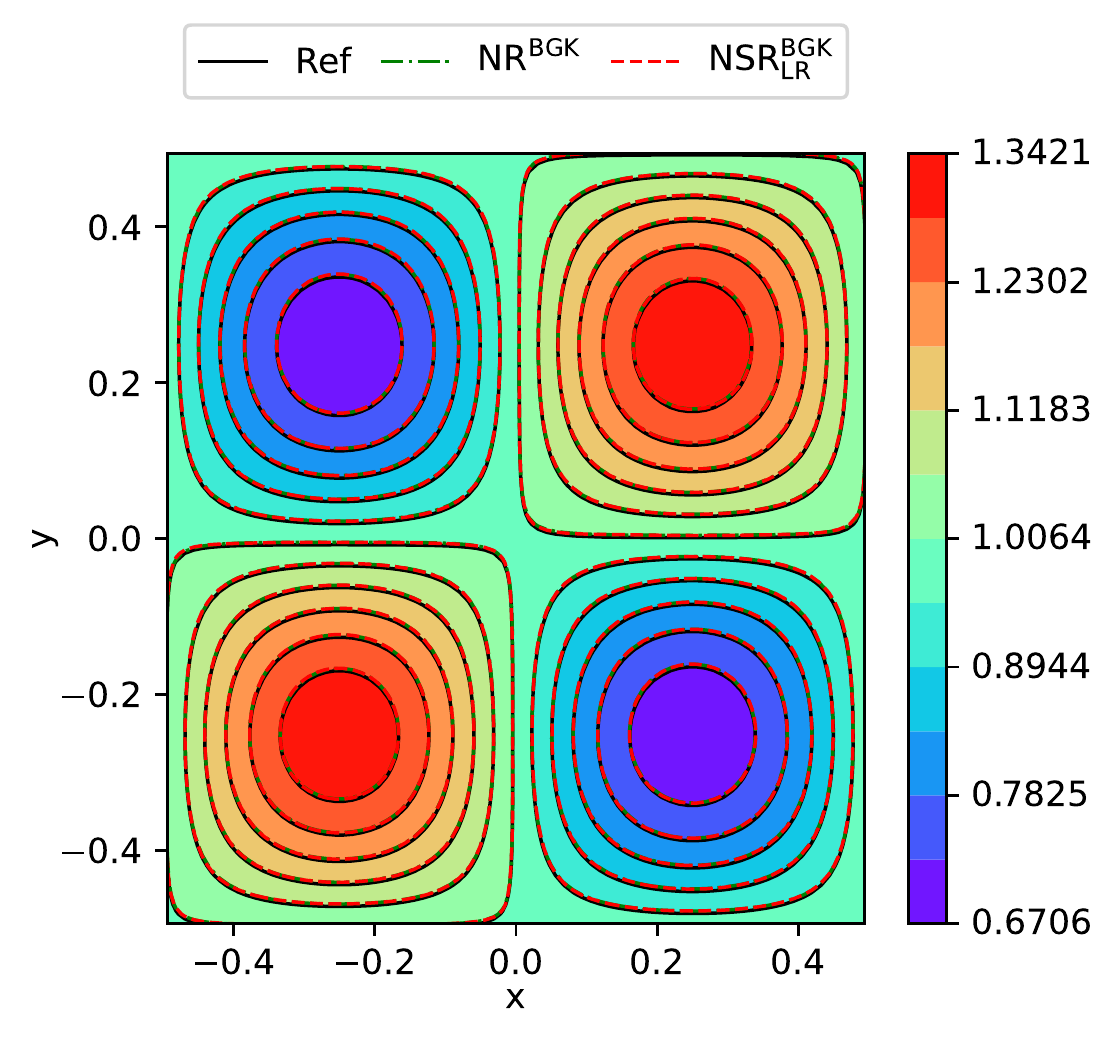}
\includegraphics[width=0.25\linewidth]{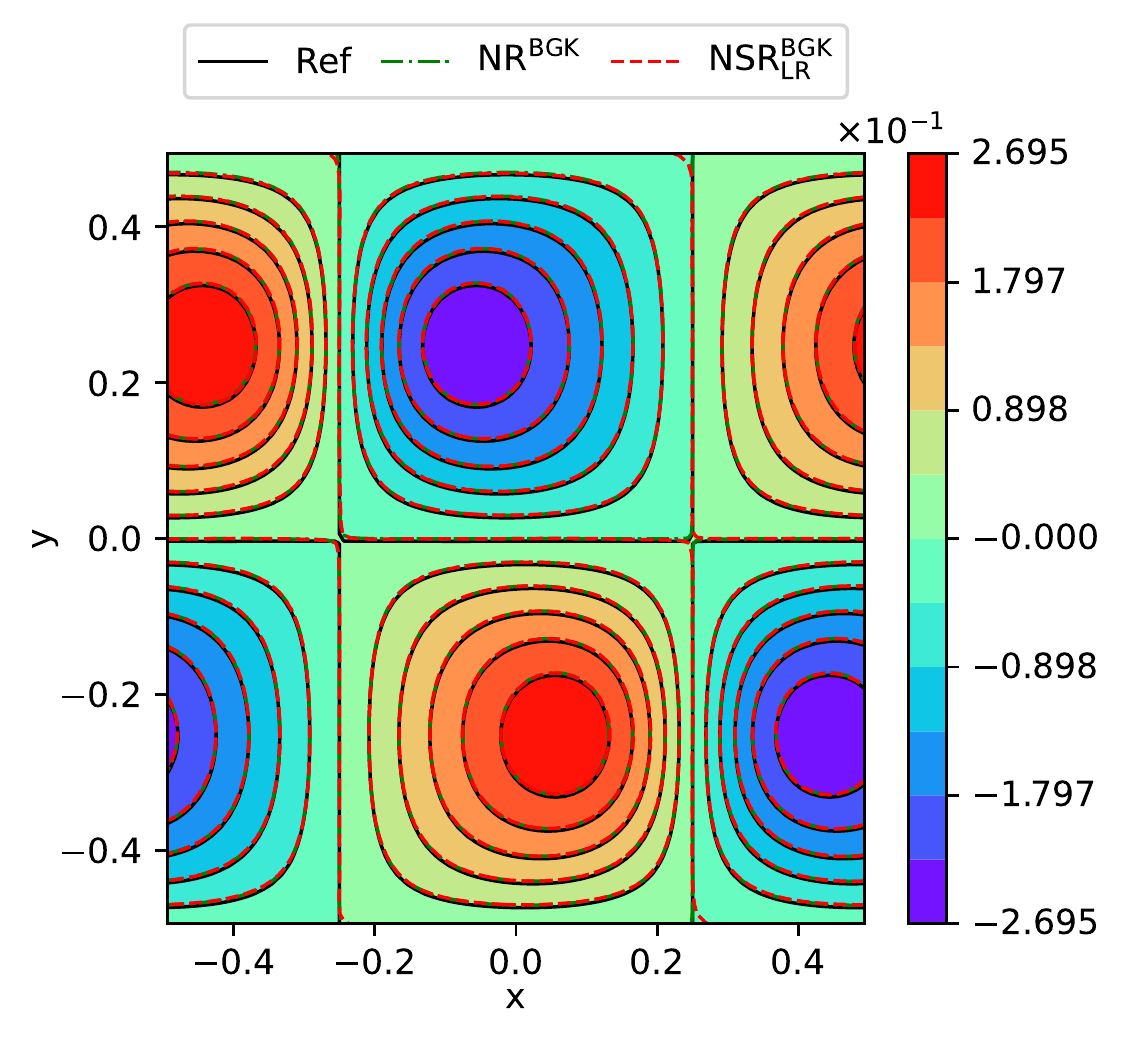}
\includegraphics[width=0.25\linewidth]{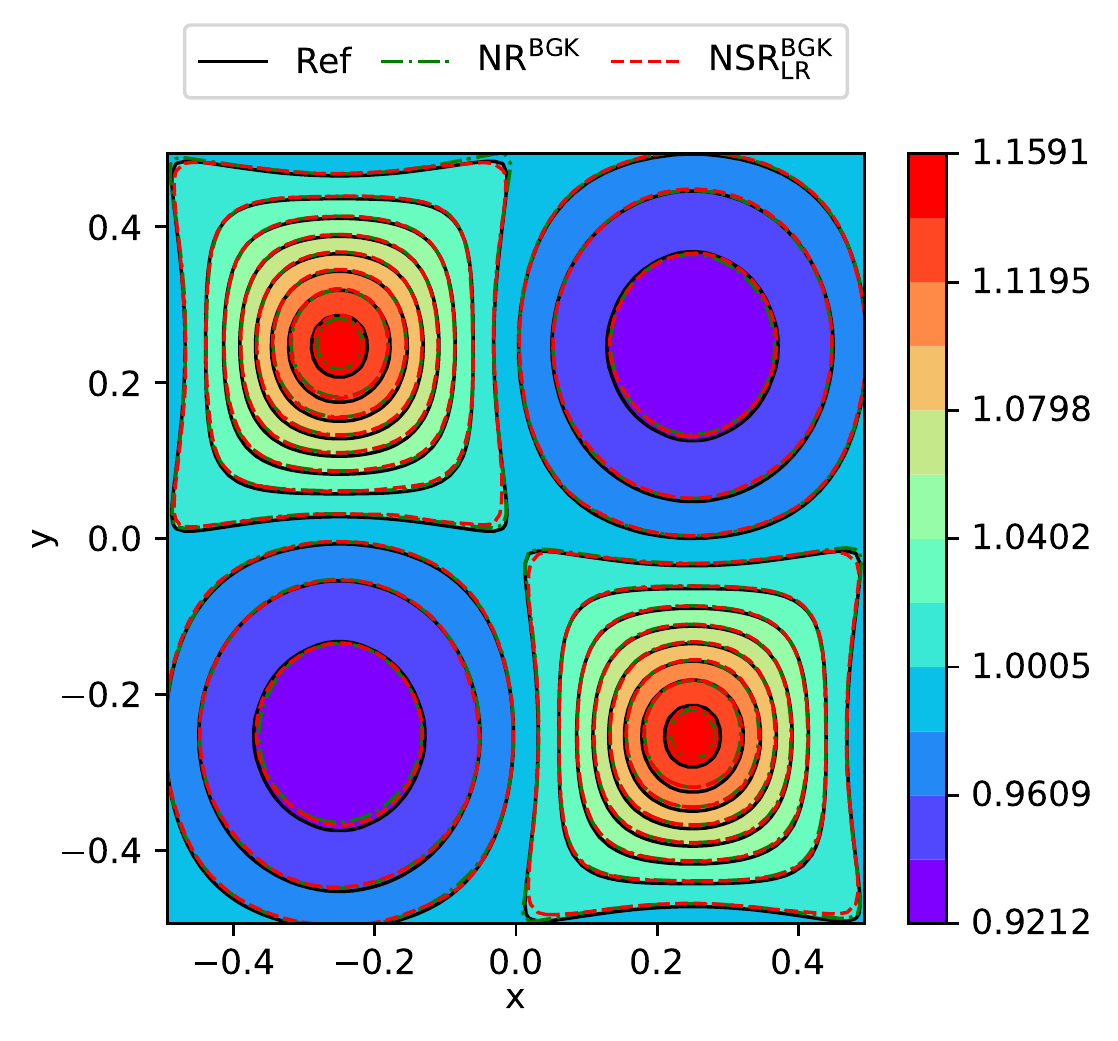}
\includegraphics[width=0.25\linewidth]{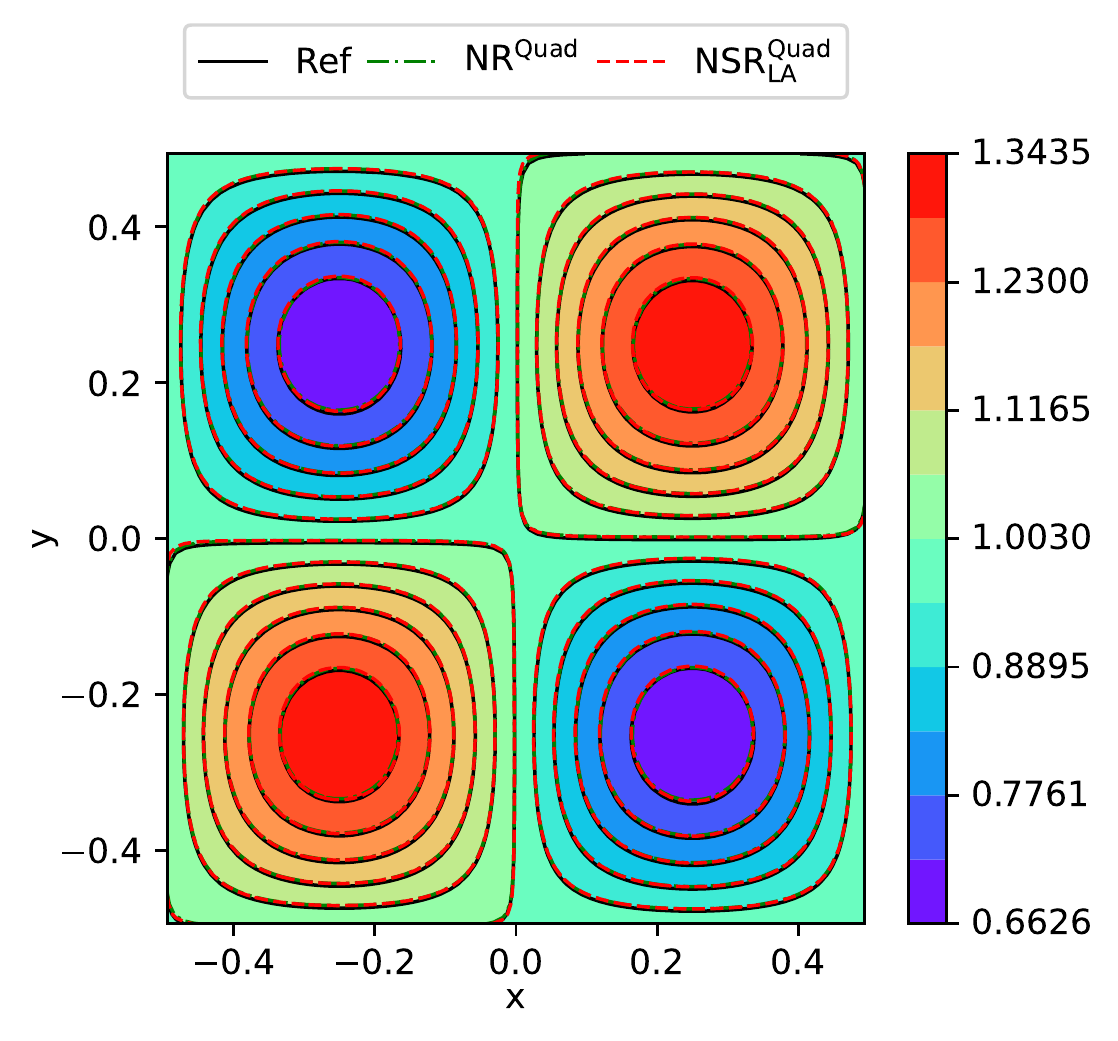}
\includegraphics[width=0.25\linewidth]{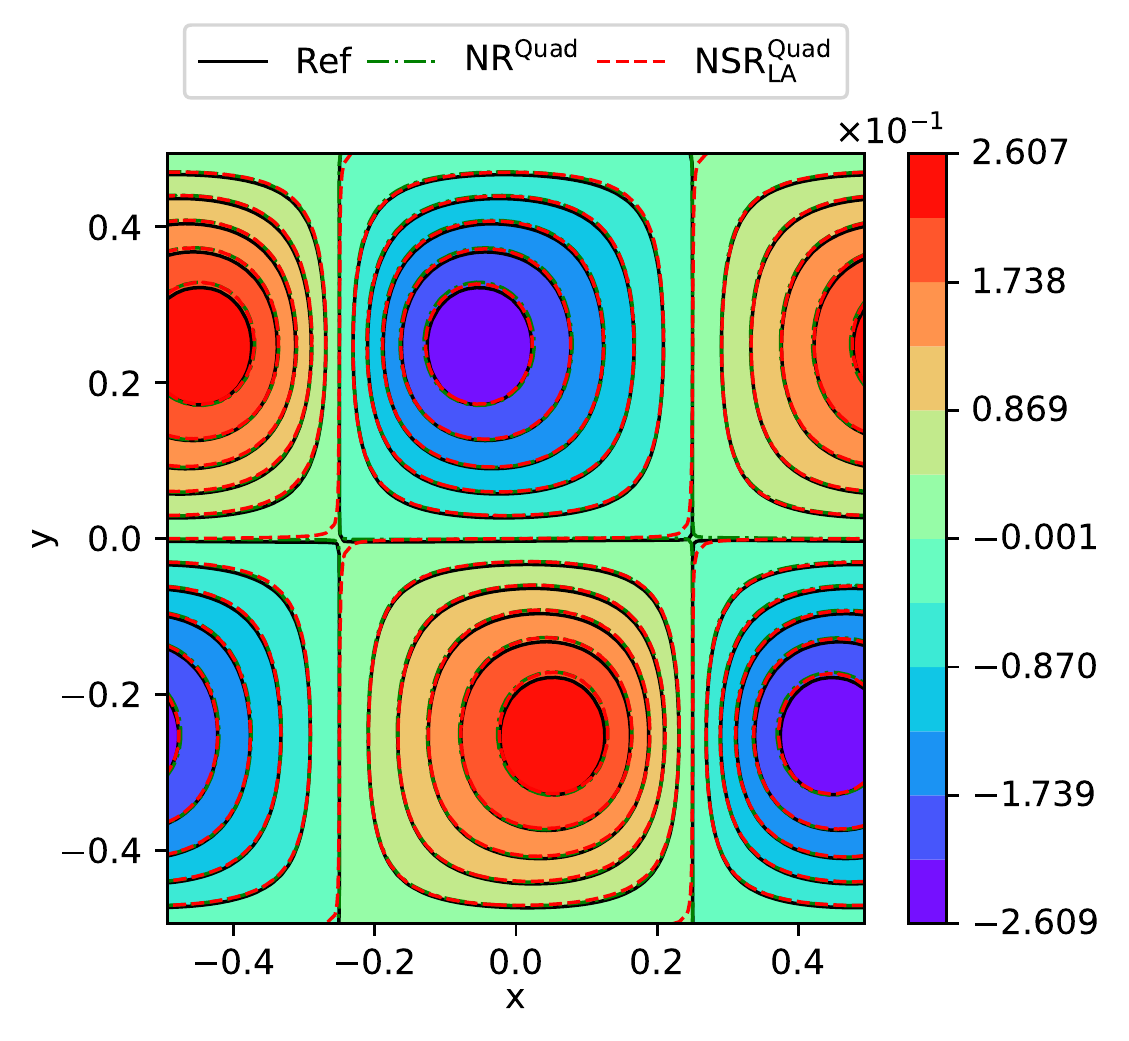}
\includegraphics[width=0.25\linewidth]{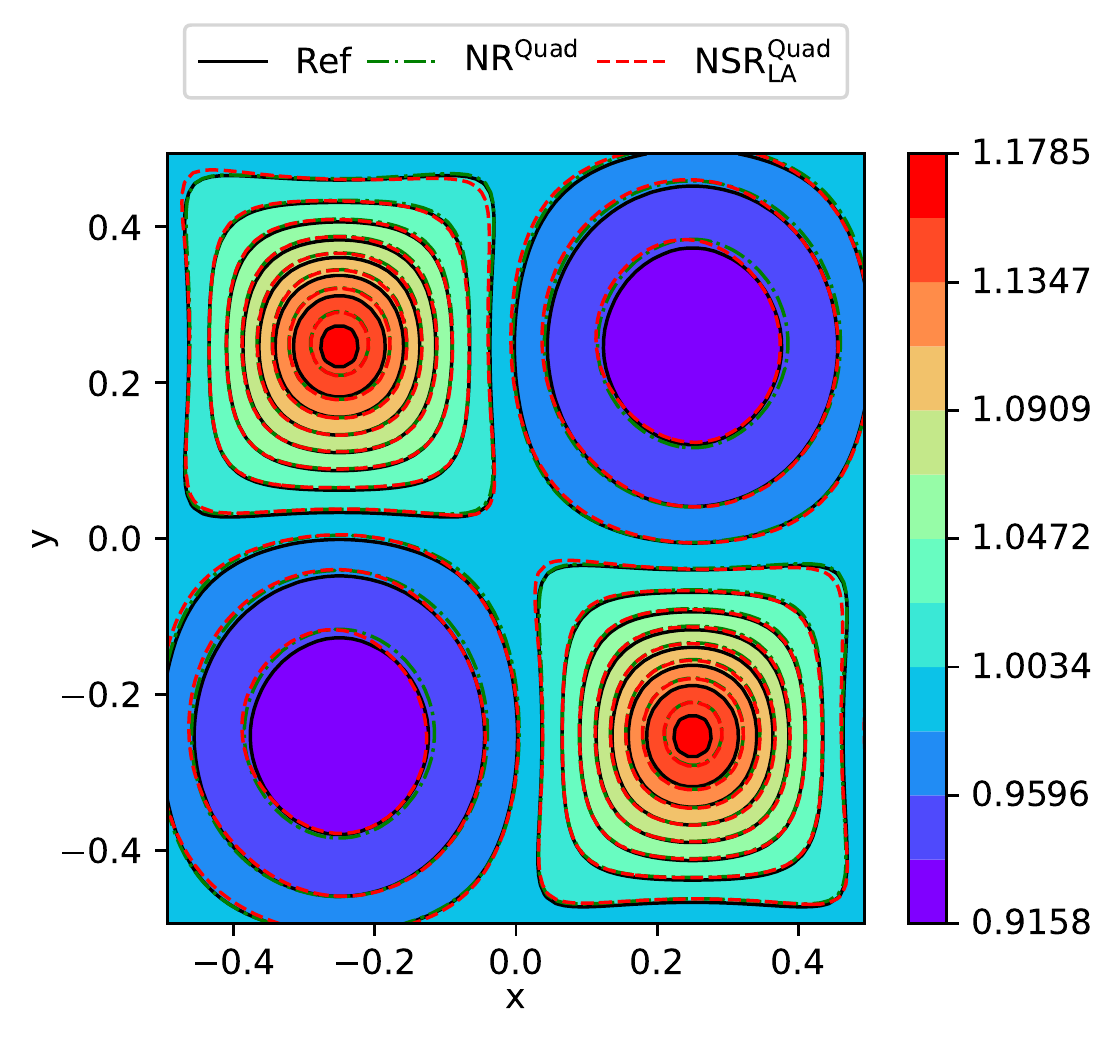}
\caption{(2D case in Sec. \ref{sec:2D}) The numerical solution of the NR/NSR method for $\Kn = 0.01$ at $ = 0.1$, where the three columns are the density $\rho$, the macroscopic velocity $u_1$, and the temperature $T$, respectively. The top row is the solution for the BGK model, and the bottom row is for the quadratic model.} 
\label{fig:wave2d-num-kn001} 
\end{figure}

In this section, the 2D3V problem with continuous initial condition is studied. The initial distribution function is Maxwellian with the macroscopic variables as follows 
\begin{equation}
\label{eq:wave2d-initial}
\rho(x, y)=1+0.5\sin(2\pi x)\sin(2\pi y), \qquad 
\bu(x, y) = \bm{0}, \qquad T(x, y)=1,
\end{equation}
with the computational domain in the spatial space $[-0.5, 0.5]^2$. Here, the periodic boundary condition is utilized, and macroscopic variables such as the density $\rho$, macroscopic velocity $\bu$, and temperature $T$ will evolve periodically as some trigonometric functions. 
\begin{figure}[!hptb]
\centering 
\includegraphics[width=0.25\linewidth]{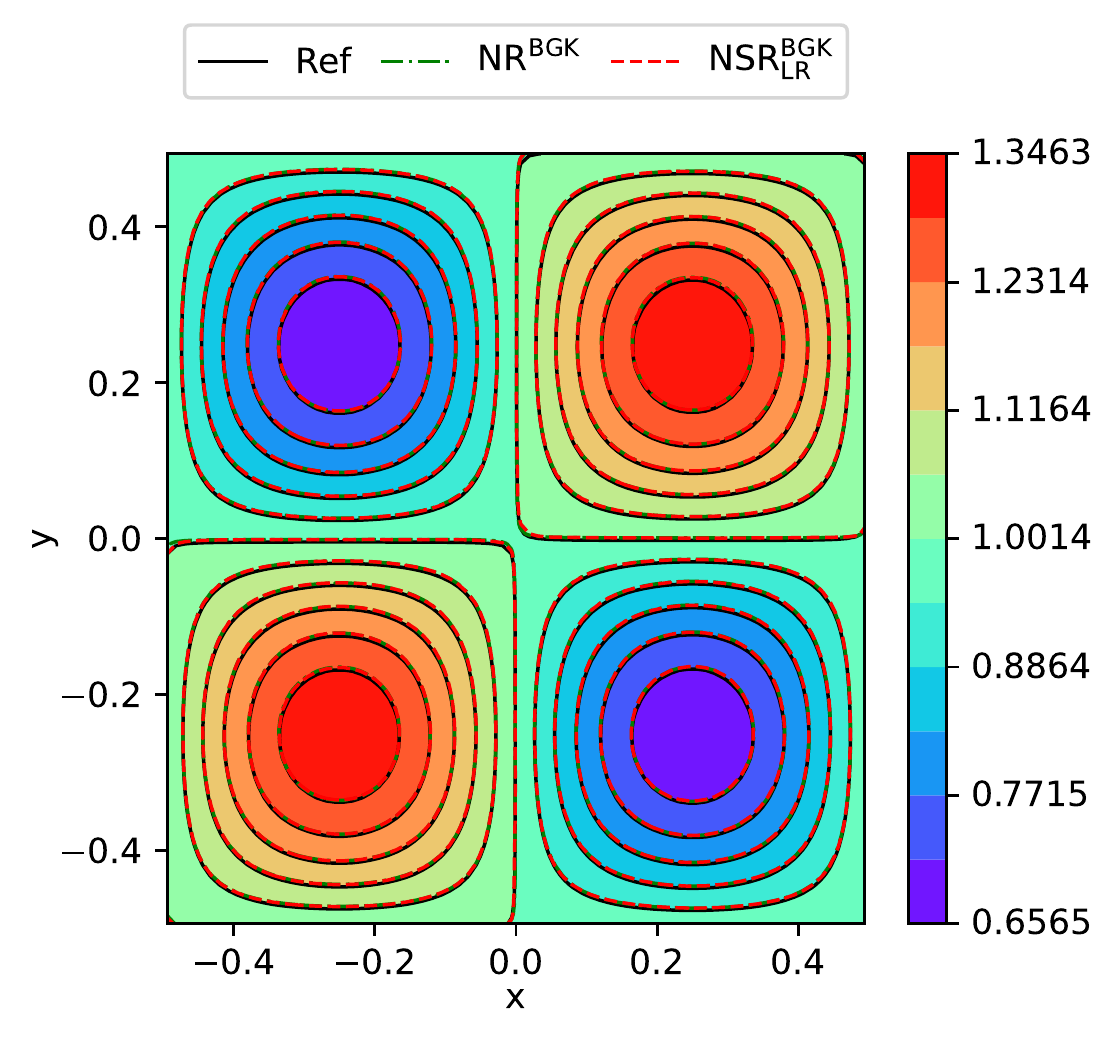}
\includegraphics[width=0.25\linewidth]{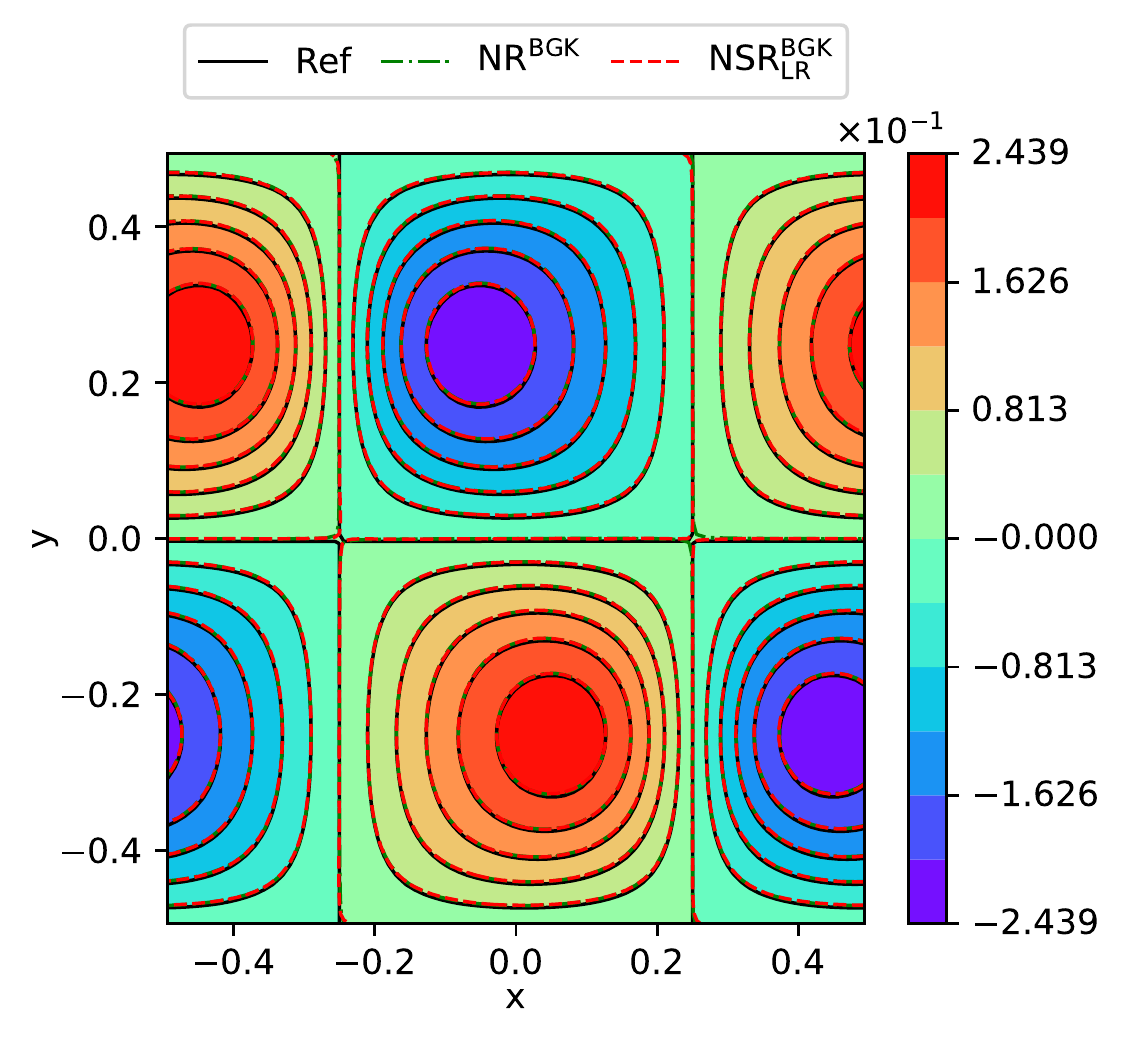}
\includegraphics[width=0.25\linewidth]{fig/exp_2d/D2V3_BGK_Kn0.1_t1_T.pdf}
\includegraphics[width=0.25\linewidth]{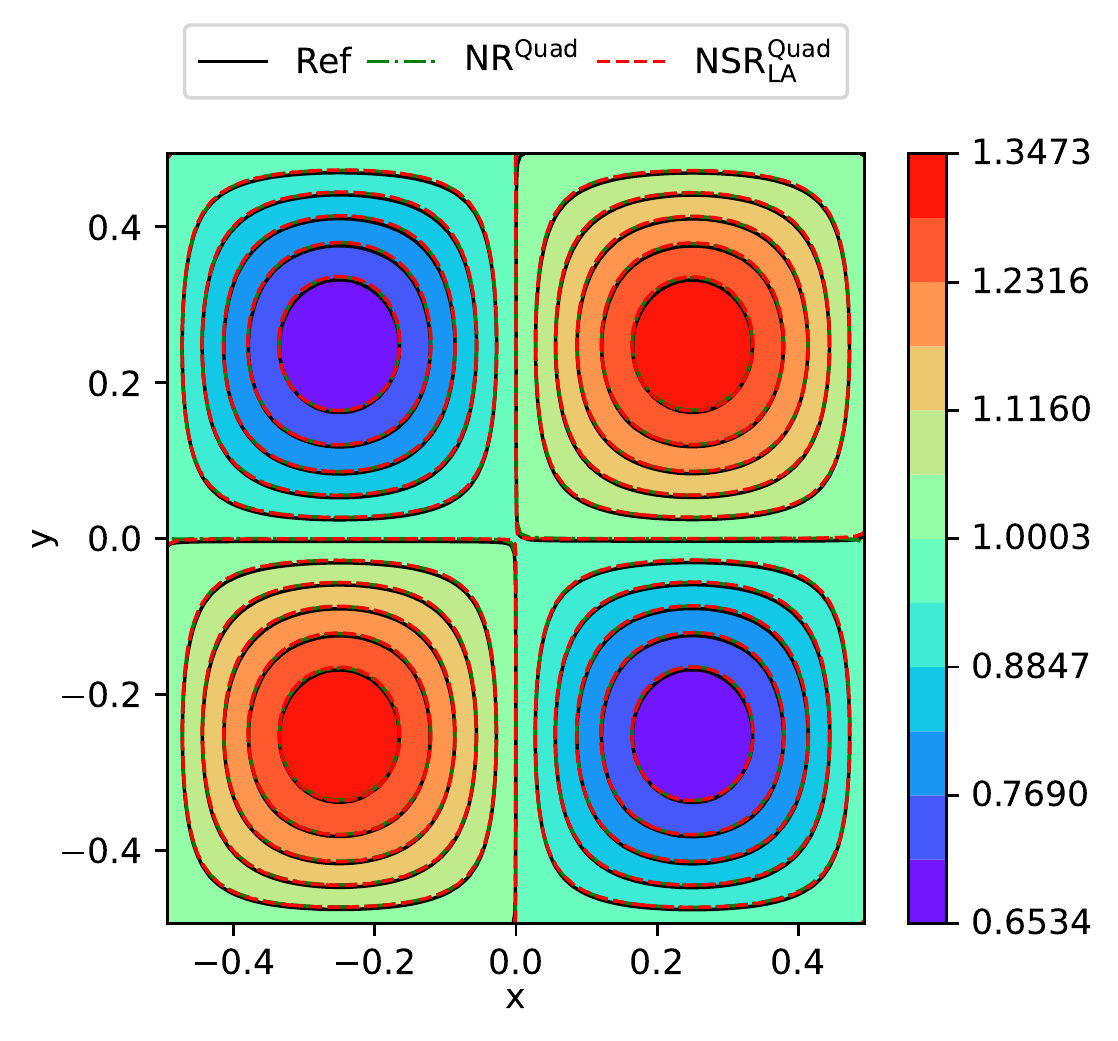}
\includegraphics[width=0.25\linewidth]{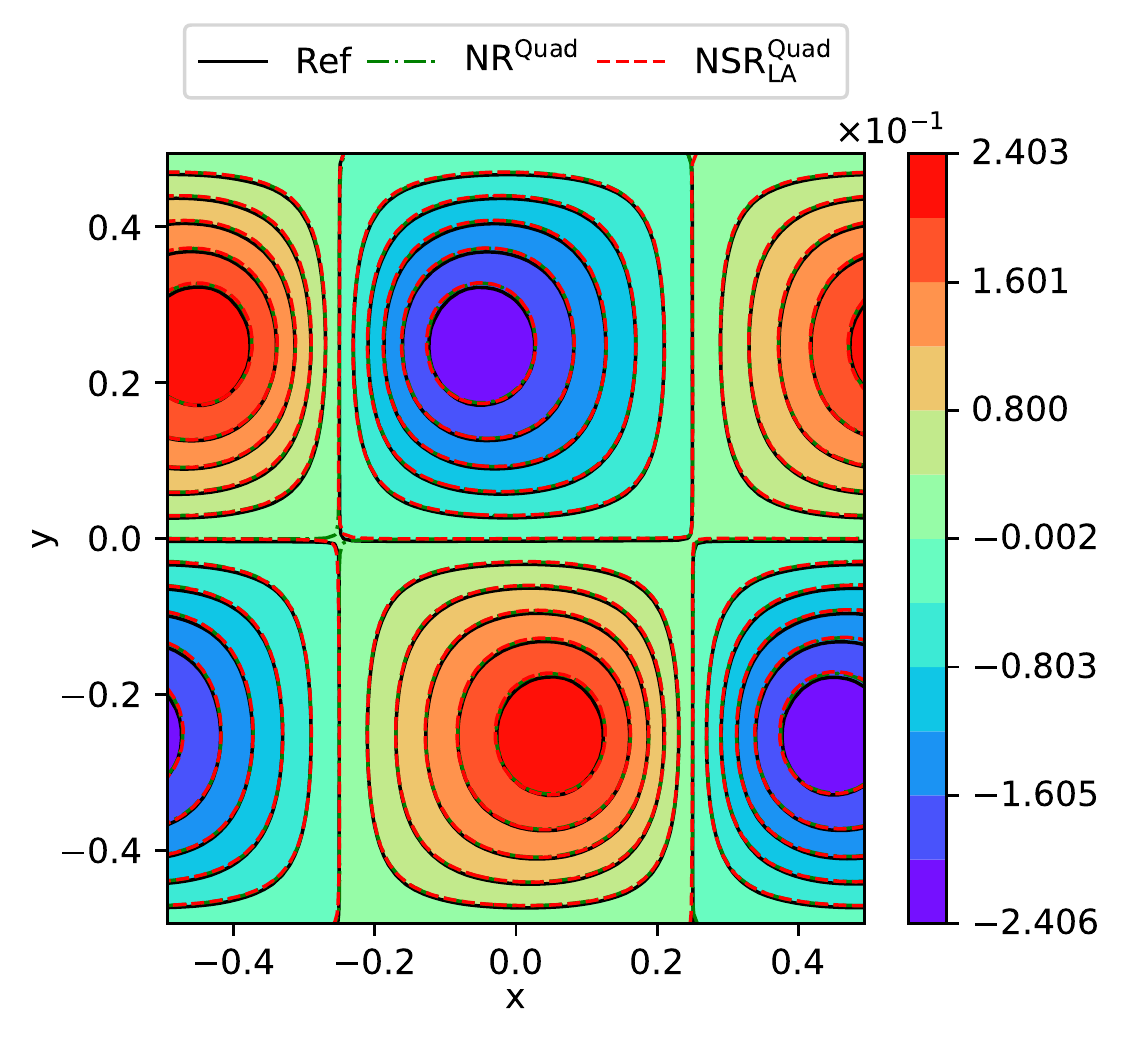}
\includegraphics[width=0.25\linewidth]{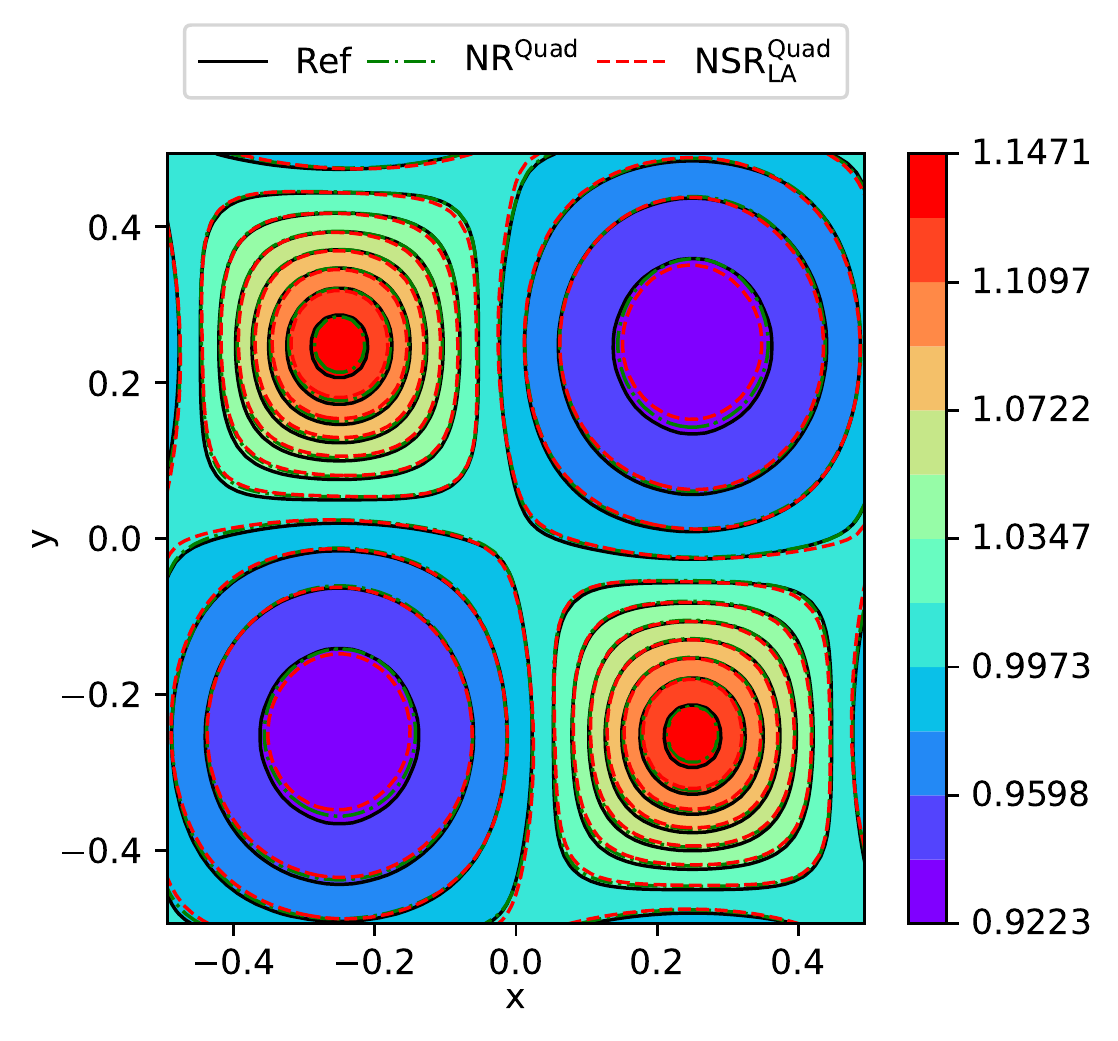}
\caption{(2D case in Sec. \ref{sec:2D}) The numerical solution of the NR/NSR method for $\Kn = 0.1$ at $t= 0.1$, where the three columns are the density $\rho$, the macroscopic velocity $u_1$, and the temperature $T$, respectively. The top row is the solution for the BGK model, and the bottom row is for the quadratic model.} 
\label{fig:wave2d-num-kn01} 
\end{figure}

\begin{figure}[!hptb]
\centering 
\includegraphics[width=0.25\linewidth]{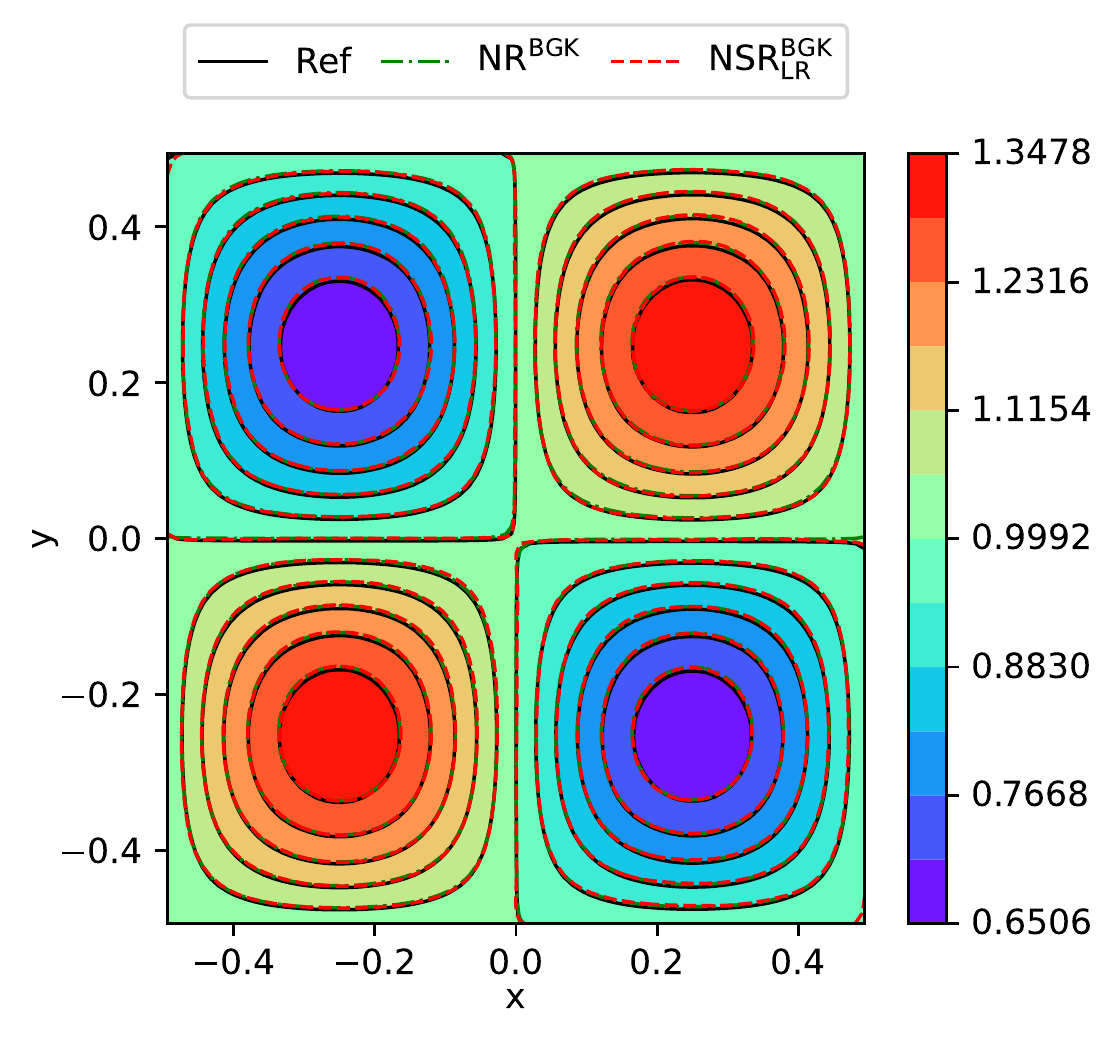}
\includegraphics[width=0.25\linewidth]{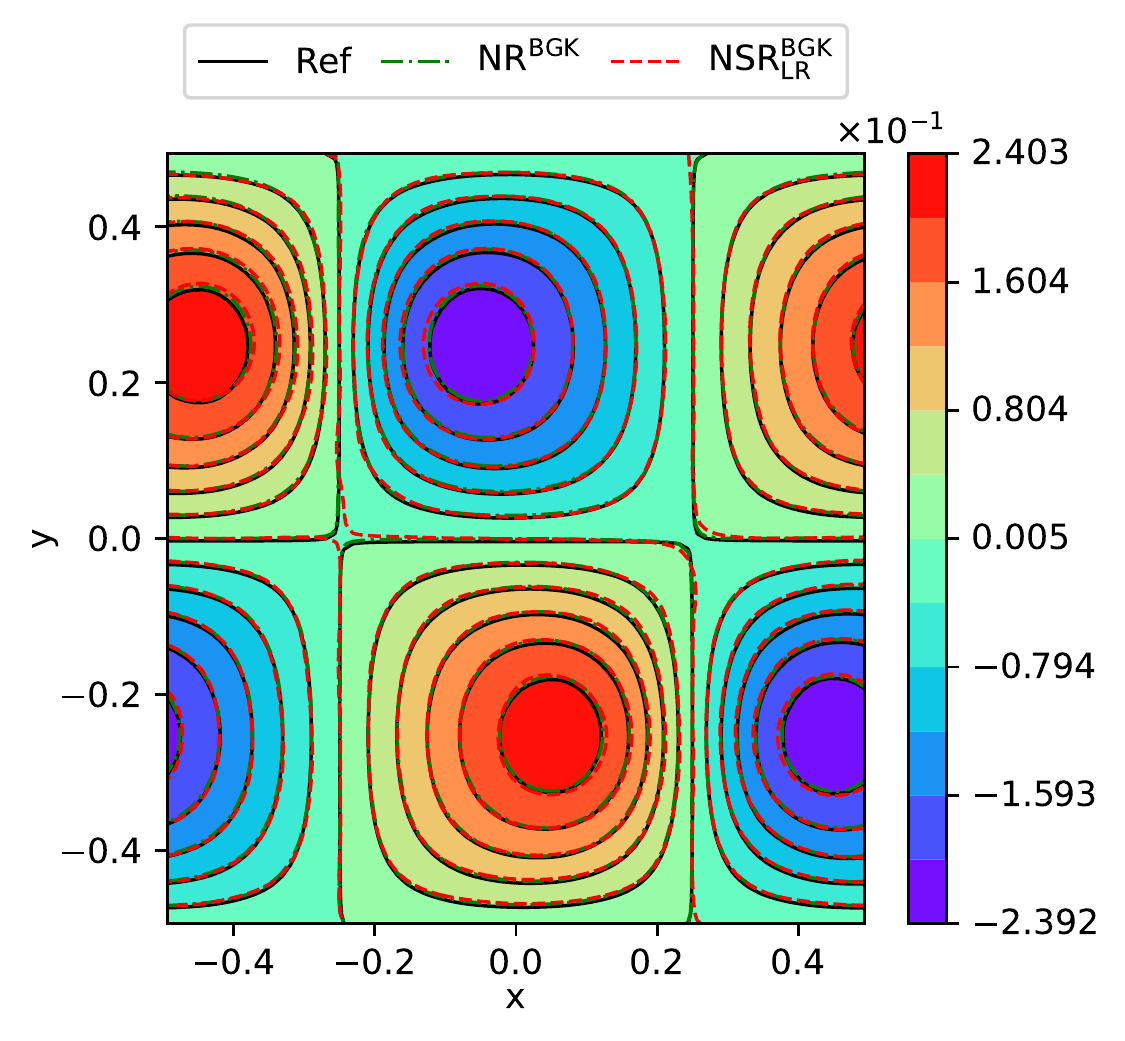}
\includegraphics[width=0.25\linewidth]{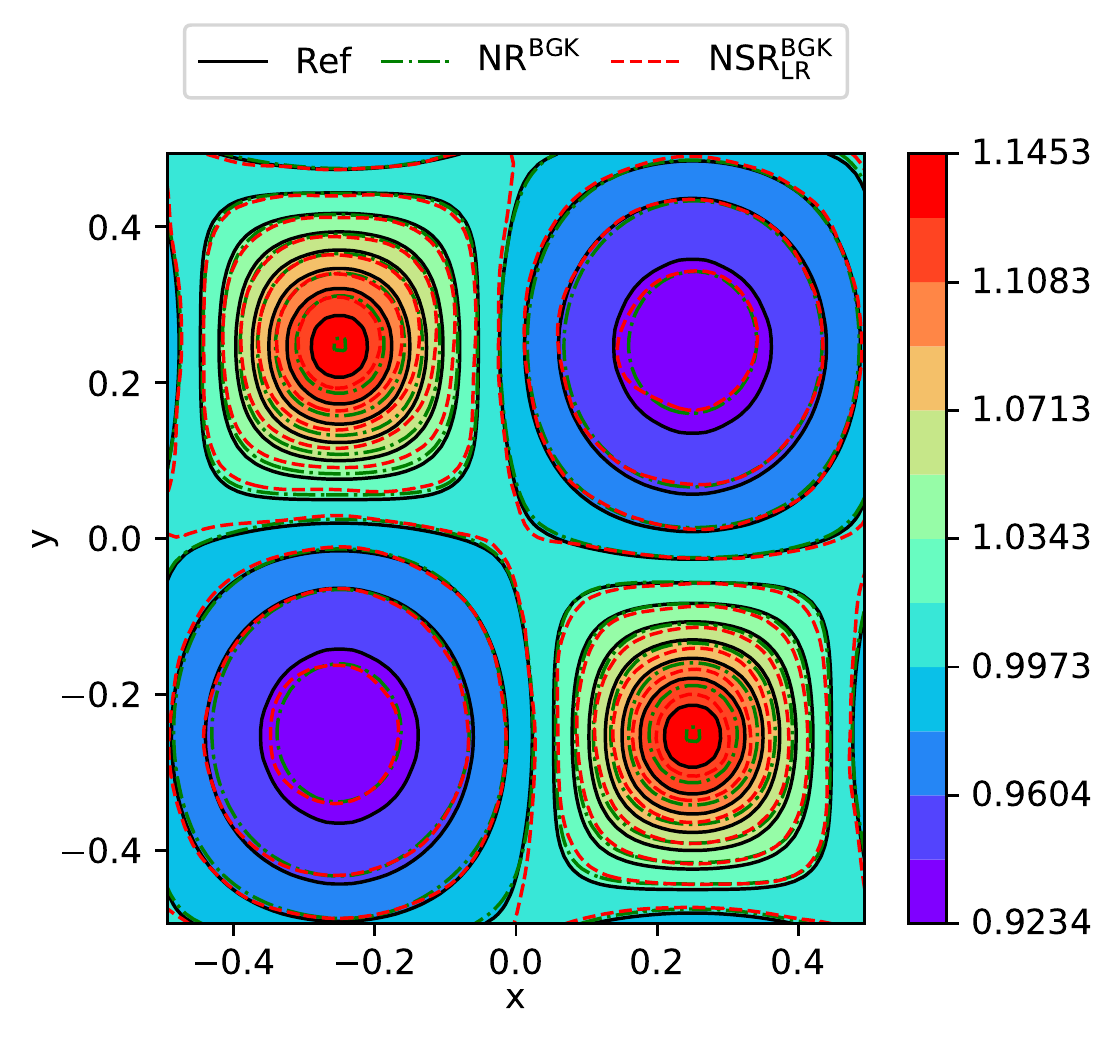}
\includegraphics[width=0.25\linewidth]{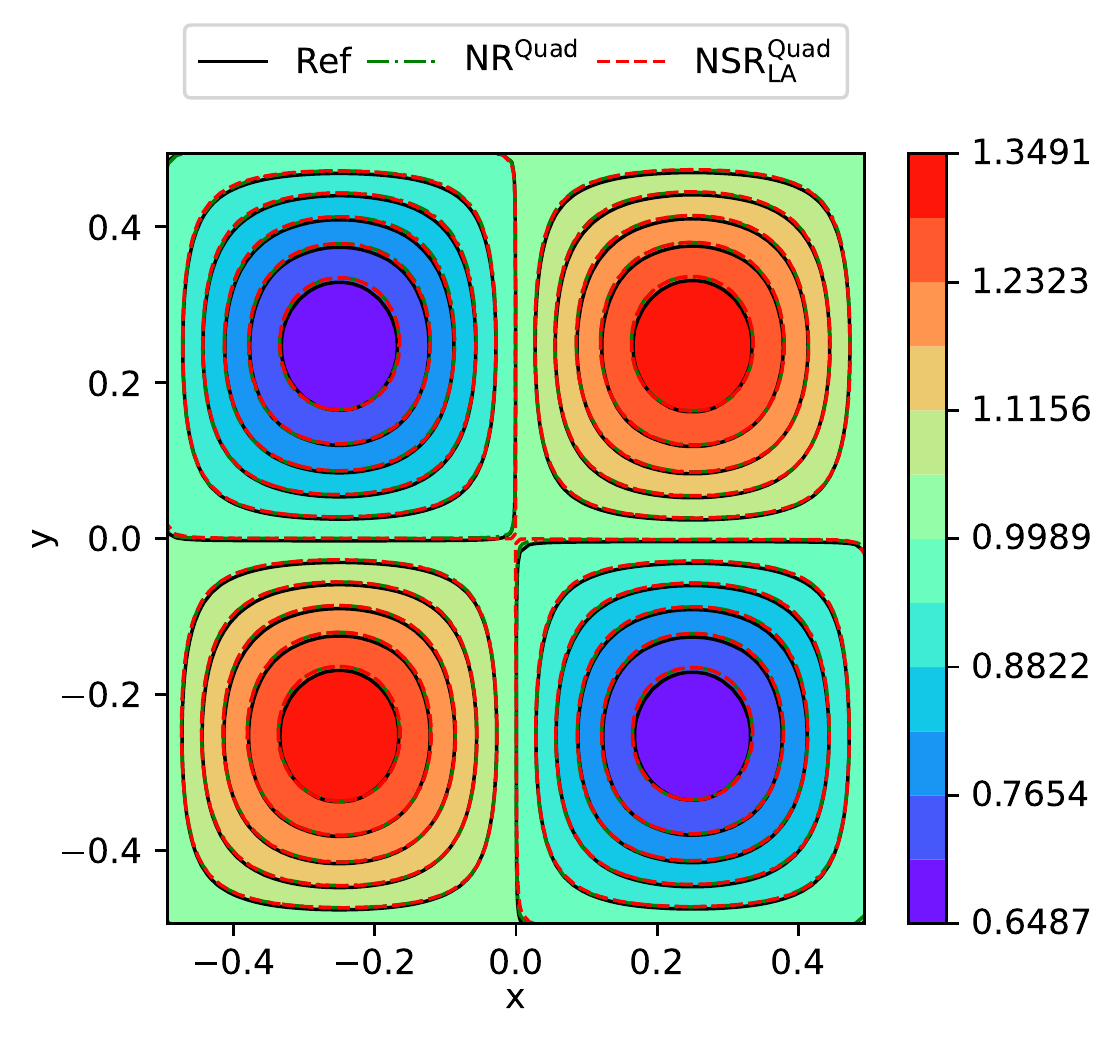}
\includegraphics[width=0.25\linewidth]{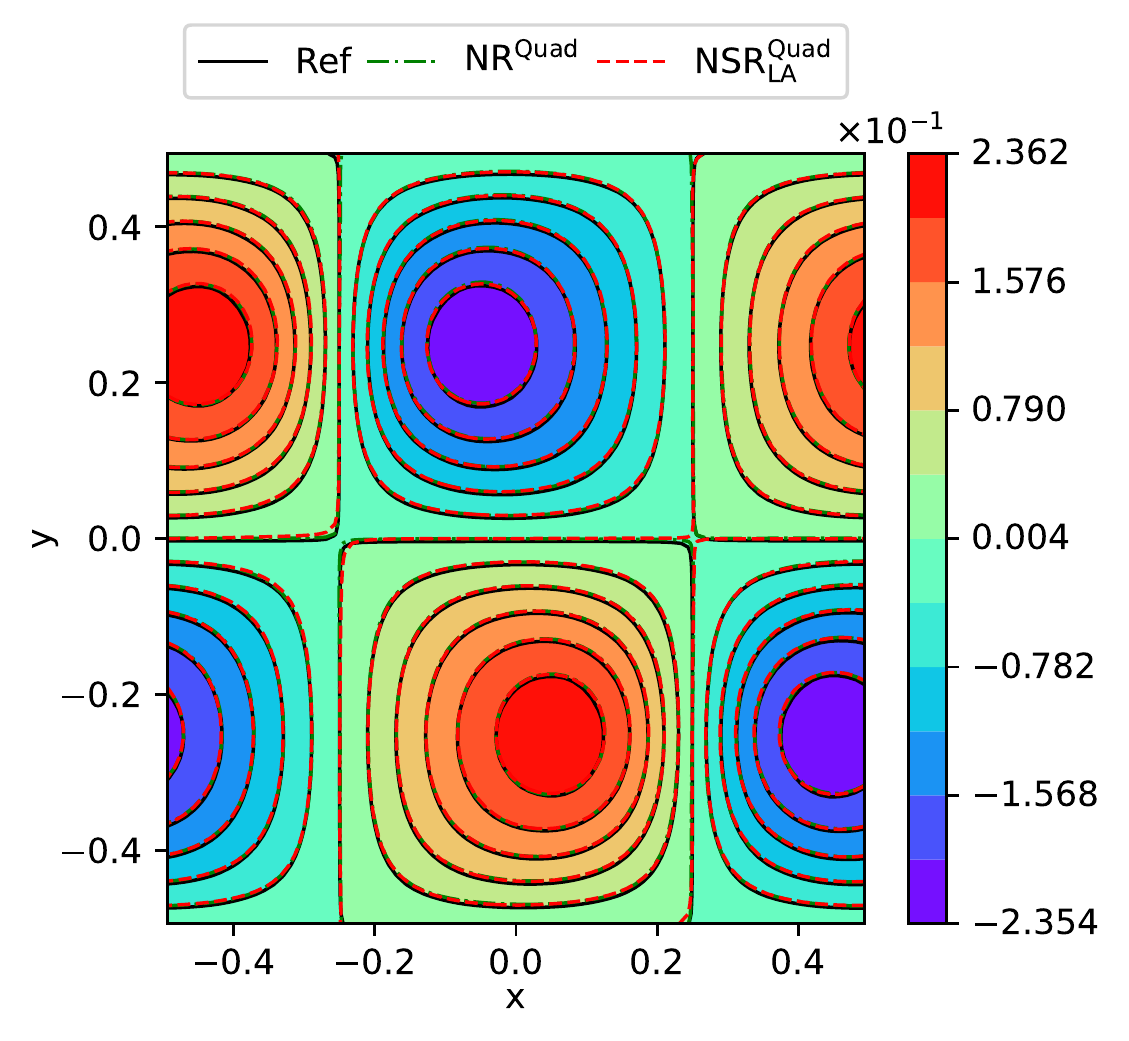}
\includegraphics[width=0.25\linewidth]{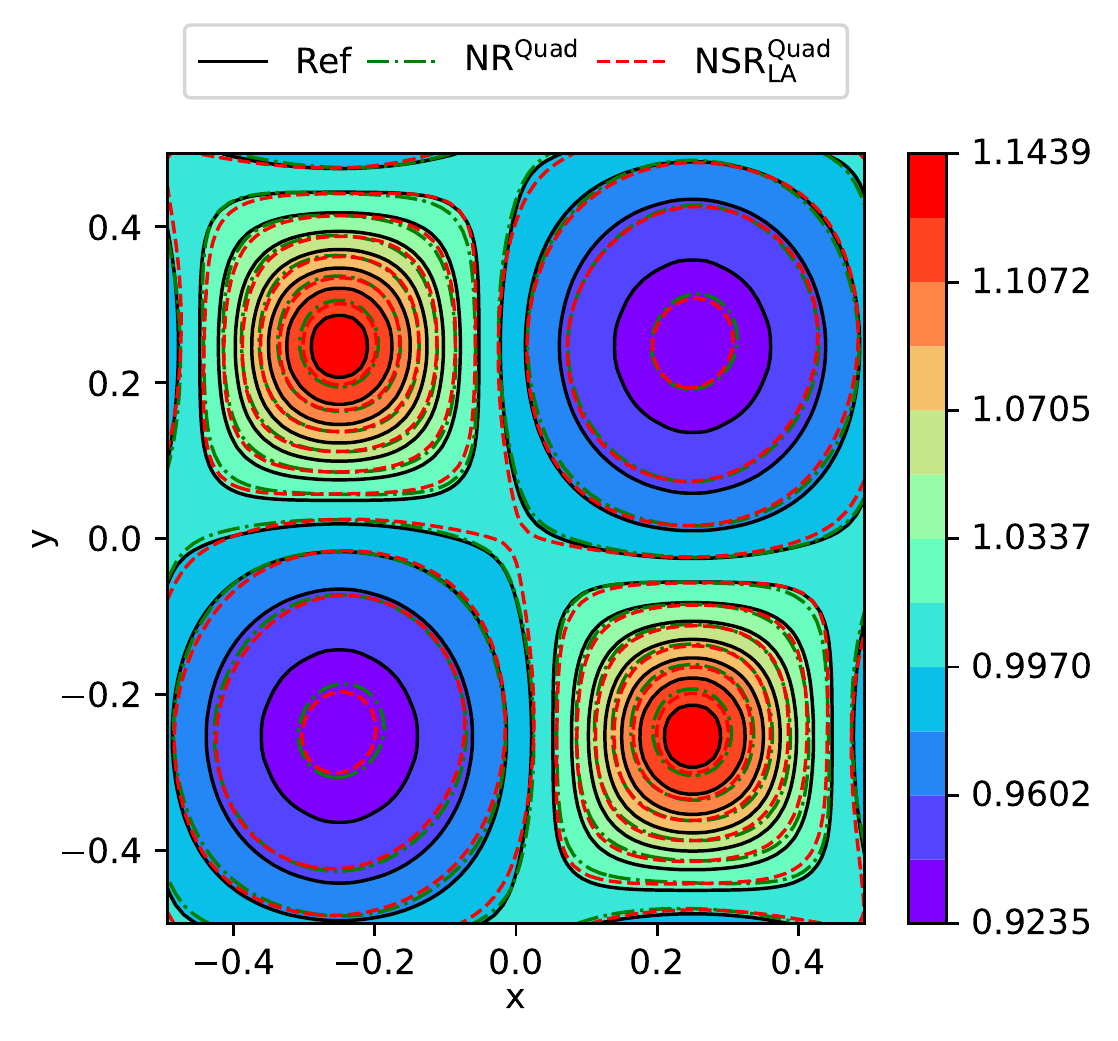}
\caption{(2D case in Sec. \ref{sec:2D}) The numerical solution of the NR/NSR method for $\Kn = 1.0$ at $t = 0.1$, where the three columns are the density $\rho$, the macroscopic velocity $u_1$, and the temperature $T$, respectively. The top row is the solution for the BGK model, and the bottom row is for the quadratic model.} 
\label{fig:wave2d-num-kn10} 
\end{figure}

The BGK and quadratic collision models are tested, where the network  and the computational parameters are the same as in Sec. \ref{subsec:wave1d}, while the sampling number is changed to $N_{\rm IC} = N_{\rm BC} = 500$, and $N_{\rm PDE} = 2000$. The numerical solution with $\nnBGK$, $\nnLR$, $\nnFSM$ and $\nnLA$ for $\Kn = 0.01$ at $t = 0.1$ are shown in Fig. \ref{fig:wave2d-num-kn001}, where the density $\rho$, the macroscopic velocity in the $x$ direction $u_1$, and the temperature $T$ are plotted. For the BGK model and the quadratic model, these three variables all agree well with the reference solution. The reference solution of BGK model is derived by the discrete velocity method, with the spatial mesh $N_x = N_y = 80$, and $24$ grids in each direction of the microscopic velocity space. The reference solution of the quadratic model is obtained by the fast Fourier spectral method, with the spatial mesh $N_x = N_y = 80$, and $24$ modes in each velocity direction. The numerical solution for $\Kn = 0.1$ and $1.0$ at $t = 0.1$ is shown in Fig. \ref{fig:wave2d-num-kn01} and \ref{fig:wave2d-num-kn10}, where the reference solution is obtained with the same parameters as in $\Kn = 0.01$. When $\Kn = 0.1$, we find the numerical solution and the reference are still on top of each other. However, when $\Kn$ is increased to $1.0$, for the density $\rho$, and the macroscopic velocity $u_1$, they match well with the reference solution, but there is a small distance for the temperature $T$. This may be due to that the reference solution is not accurate enough, but this parameter setting has already the maximum memory we can afford. 

To exhibit the numerical error quantitatively, the relative error \eqref{eq:error1} and \eqref{eq:error2} of the four methods with different Knudsen numbers at $t = 0$, and $0.1$ are shown in Tab. \ref{tab:wave2d}. It shows that for the initial data, this error is relatively small, most at the order of $\mO(10^{-4})$, and at $t = 0.1$, this error is increased to $\mO(10^{-3})$, but all at the same order. 

\begin{table}[hpt!]
\centering
\def\arraystretch{1.5}
\scalebox{0.85}{
{\footnotesize
\begin{tabular}{lllllllllll}
\hline
Kn                               & \multicolumn{1}{c}{}    & \multicolumn{3}{c}{0.01}                                                             & \multicolumn{3}{c}{0.1}                                                              & \multicolumn{3}{c}{1.0}                                                        \\ \hline
                                 & \multicolumn{1}{c}{$t$} & \multicolumn{1}{c}{$\rho$} & \multicolumn{1}{c}{$u$} & \multicolumn{1}{c|}{$T$}      & \multicolumn{1}{c}{$\rho$} & \multicolumn{1}{c}{$u$} & \multicolumn{1}{c|}{$T$}      & \multicolumn{1}{c}{$\rho$} & \multicolumn{1}{c}{$u$} & \multicolumn{1}{c}{$T$} \\
\multirow{2}{*}{$\nnBGK$} & 0.0                     & 3.05e-04                   & 1.70e-04                & \multicolumn{1}{l|}{3.83e-04} & 2.00e-04                   & 1.74e-04                & \multicolumn{1}{l|}{2.35e-04} & 4.08e-04                   & 2.28e-04                & 6.14e-04                \\
                                 & 0.1                     & 3.20e-03                   & 2.44e-03                & \multicolumn{1}{l|}{1.39e-03} & 3.29e-03                   & 2.23e-03                & \multicolumn{1}{l|}{1.45e-03} & 3.45e-03                   & 2.25e-03                & 3.11e-03                \\
\multirow{2}{*}{$\nnLR$} & 0.0                     & 2.03e-04                   & 2.20e-04                & \multicolumn{1}{l|}{2.76e-04} & 1.26e-03                   & 4.02e-04                & \multicolumn{1}{l|}{1.33e-03} & 1.27e-03                   & 3.12e-04                & 1.39e-03                \\
                                 & 0.1                     & 4.20e-03                   & 3.10e-03                & \multicolumn{1}{l|}{6.87e-03} & 1.27e-03                   & 3.12e-04               & \multicolumn{1}{l|}{1.39e-03} & 3.55e-03                   & 3.00e-03                & 4.70e-03                \\
\multirow{2}{*}{$\nnFSM$}  & 0.0                     & 3.40e-04                   & 1.41e-04                & \multicolumn{1}{l|}{5.28e-04} & 3.58e-04                   & 2.06e-04                & \multicolumn{1}{l|}{4.82e-04} & 5.20e-04                   & 2.12e-04                & 7.26e-04                \\
                                 & 0.1                     & 3.24e-03                   & 2.88e-03                & \multicolumn{1}{l|}{2.18e-03} & 3.30e-03                   & 2.29e-03                & \multicolumn{1}{l|}{1.30e-03} & 3.62e-03                   & 2.41e-03                & 4.81e-03                \\
\multirow{2}{*}{$\nnLA$} & 0.0                     & 2.85e-04                   & 1.42e-04                & \multicolumn{1}{l|}{2.09e-04} & 2.80e-04                   & 2.21e-04                & \multicolumn{1}{l|}{3.14e-04} & 6.12e-04                   & 2.38e-04                & 5.01e-04                \\
                                 & 0.1                     & 3.22e-03                   & 2.73e-03                & \multicolumn{1}{l|}{2.22e-03} & 3.31e-03                   & 2.50e-03                & \multicolumn{1}{l|}{2.09e-03} & 3.69e-03                   & 2.36e-03                & 5.37e-03                \\ \bottomrule
\end{tabular}}
}
\caption{ (2D case in Sec. \ref{sec:2D}) The relative error between the numerical solution by NR/NSR and the reference solution for the density $\rho$, macroscopic velocity $u_1$ and the temperature $T$
with $\Kn = 0.01, 0.1$ and $1$ at $t = 0$ and $0.1$.
}
\label{tab:wave2d}
\end{table}

\subsection{Transfer learning}
\label{sec:trans}
One limitation of the network-based method to solve PDEs is the slow training speed, and there is some work to improve this, such as the transfer learning \cite{chen2021transfer} and manifold learning \cite{huang2022metaautodecoder}. For the Boltzmann equation, especially for the BGK model, the classical methods such as DVM are always more efficient than the network-based method, even when transfer learning is utilized. But for high dimensional problems, the network-based method is more competitive. 

To explore the efficiency of the NR/NSR method, we studied the computational time utilized for the 2-dimensional problems with transfer learning. Generally speaking, transfer learning is to transfer the network trained on one task to another similar new task as the initial network \cite{chen2021transfer}. Thus, it is expected that the learning process can be speeded up, which is also verified by the numerical experiments. The initial condition for the new 2D3V problem is 
\begin{equation}
\label{eq:ex4_ini}
\rho(x)=1+0.4\sin(2\pi x+0.3))\sin(2\pi (y+0.4)), \qquad {\bm u}(x, y) = 0, \qquad T(x, y) = 1,
\end{equation}
with the computational domain $[-0.5, 0.5]^2$. In the learning process, instead of initializing the neural network randomly, the well-trained network in Sec \ref{sec:2D} is adopted. The other parameters including those for the reference solution are the same as in Sec. \ref{sec:2D}. 

\begin{table}[H]
\centering
\def\arraystretch{1.5}
\scalebox{0.85}{
\footnotesize
\begin{tabular}{@{}llll@{}}
\toprule
Case        & 1  & 2  & 3  \\ 
$N_x$         &  20  &  30  &  40  \\
\bottomrule
Err         &  8.43e-03  &  4.54e-03  &  2.71e-03  \\
Ref-BGK     &  58   &  210  &   595  \\
$\nnBGK$ &  33 &  148  &  -  \\
$\nnLR$  &  32  &   93  &  -  \\
\bottomrule
Err         &  8.45e-03  &  4.55e-03  &  2.72e-03 \\
Ref-FSM     &  300  &  1086  &  2549  \\
$\nnFSM$ &  139  &  406  &  -  \\
$\nnLA$  &  86  &  213  &  -  \\
\bottomrule
\end{tabular}
}
\caption{(Transfer learning in Sec. \ref{sec:trans}) Computational time of the classical methods and the network-based methods to achieve similar accuracy. `` -" indicates that precision could not be achieved.
The time for $\nnLA$ includes that to obtain the collision kernel and to train the neural network.
}
\label{tab:transfer}
\end{table}

We first set the mesh number as $N_x = N_y = 20, 30$ and $40$, and record the error between the numerical solution with different grid numbers using  DVM for BGK model and fast Fourier spectral for the quadratic model. These errors are shown in the third and seventh rows of Tab. \ref{tab:transfer}. Then, the simulations with NR/NSR are carried out to reach the same accuracy. The computational time for the classical methods and the network-based methods are all illustrated in Tab. \ref{tab:transfer}. The reference solver of the classical methods is a Fortran program, running on a server with 72 cores (2 Intel(R) Xeon(R) Gold 6240 CPUs @ 2.60GHz) and 256GB of RAM. NR/NSR are Python programs running on the same server with a single RTX 3090 (with 24Gb of video memory) based on pytorch. Tab. \ref{tab:transfer} exhibits that to arrive at the same accuracy, the computational time of the network-based methods is much shorter compared to the classical methods. Moreover, the neural sparse representation ($\nnLR$ and $\nnLA$) are much faster than the general neural representation methods, especially for the quadratic collision model. These indicate that the neural sparse representation is more efficient for high-dimensional problems, and can be quite promising for solving the 3D3V full Boltzmann equations. However, the network-based methods fail to achieve high precision for the moment, which we will work on in the future.

\section{Conclusion}
\label{sec:conclusion}
The neural network-based approach is utilized to solve the Boltzmann equation. Neural sparse representation for the distribution function is proposed, which is a high-quality ansatz to the Boltzmann equation. The low-rank property of the discrete distribution function is adopted in the BGK model, and a network structure whose output is the CPD factorization of the discrete distribution function is proposed, which effectively reduces the complexity of the network parameters. For the quadratic collision model, the data-driven basis vectors are constructed with the BGK solution through SVD. The quadratic collision term can be approximated with this series of linear basis vectors with much less freedom. Adaptive weight loss function, which includes the initial, boundary conditions and residual of PDE  and the loss from the macroscopic variables,  is designed for the learning process and has greatly improved the approximating efficiency of the network. Numerical examples of the 1D and 2D cases are studied to validate the accuracy and efficiency of these neural representation methods. The effect of transfer learning is studied to show the efficiency of these methods, and more work will be done in the future.


\section*{Acknowledgments}
We thank Dr. Chang Liu from Institute of Applied Physics and Computational Mathematics for the code of the fast Fourier spectral method. This work of Y. Wang is partially supported by the National Natural Science Foundation of China (Grant No. 12171026, U2230402 and 12031013), and Foundation of President of China Academy of Engineering Physics (YZJJZQ2022017).


\bibliographystyle{plain}
\bibliography{article}  

\end{document}